\documentclass[11pt]{article}
\usepackage{amssymb,latexsym,amsmath}
\usepackage{graphicx}

\begin{document}

\newcommand{\bfi}{\bfseries\itshape}

\makeatletter

\@addtoreset{figure}{section}
\def\thefigure{\thesection.\@arabic\c@figure}
\def\fps@figure{h, t}
\@addtoreset{table}{bsection}
\def\thetable{\thesection.\@arabic\c@table}
\def\fps@table{h, t}
\@addtoreset{equation}{section}
\def\theequation{\thesubsection.\arabic{equation}}

\makeatother

\newtheorem{thm}{Theorem}[section]
\newtheorem{prop}[thm]{Proposition}
\newtheorem{lema}[thm]{Lemma}
\newtheorem{cor}[thm]{Corollary}
\newtheorem{defi}[thm]{Definition}
\newtheorem{rk}[thm]{Remark}

\newcommand{\comment}[1]{\par\noindent{\raggedright\texttt{#1}
\par\marginpar{\textsc{Comment}}}}
\newcommand{\todo}[1]{\vspace{5 mm}\par \noindent \marginpar{\textsc{ToDo}}\framebox{\begin{minipage}[c]{0.95 \textwidth}
\tt #1 \end{minipage}}\vspace{5 mm}\par}

\newcommand{\ea}{\mbox{{\bf a}}}
\newcommand{\eu}{\mbox{{\bf u}}}
\newcommand{\ueu}{\underline{\eu}}
\newcommand{\ueo}{\overline{u}}
\newcommand{\oeu}{\overline{\eu}}
\newcommand{\ew}{\mbox{{\bf w}}}
\newcommand{\ef}{\mbox{{\bf f}}}
\newcommand{\eF}{\mbox{{\bf F}}}
\newcommand{\eC}{\mbox{{\bf C}}}
\newcommand{\en}{\mbox{{\bf n}}}
\newcommand{\eT}{\mbox{{\bf T}}}
\newcommand{\eL}{\mbox{{\bf L}}}
\newcommand{\eR}{\mbox{{\bf R}}}
\newcommand{\eV}{\mbox{{\bf V}}}
\newcommand{\eU}{\mbox{{\bf U}}}
\newcommand{\ev}{\mbox{{\bf v}}}
\newcommand{\eve}{\mbox{{\bf e}}}
\newcommand{\uev}{\underline{\ev}}
\newcommand{\eY}{\mbox{{\bf Y}}}
\newcommand{\eK}{\mbox{{\bf K}}}
\newcommand{\eP}{\mbox{{\bf P}}}
\newcommand{\eS}{\mbox{{\bf S}}}
\newcommand{\eJ}{\mbox{{\bf J}}}
\newcommand{\eB}{\mbox{{\bf B}}}
\newcommand{\eH}{\mbox{{\bf H}}}
\newcommand{\leb}{\mathcal{ L}^{n}}
\newcommand{\eI}{\mathcal{ I}}
\newcommand{\eE}{\mathcal{ E}}
\newcommand{\hen}{\mathcal{H}^{n-1}}
\newcommand{\eBV}{\mbox{{\bf BV}}}
\newcommand{\eA}{\mbox{{\bf A}}}
\newcommand{\eSBV}{\mbox{{\bf SBV}}}
\newcommand{\eBD}{\mbox{{\bf BD}}}
\newcommand{\eSBD}{\mbox{{\bf SBD}}}
\newcommand{\ecs}{\mbox{{\bf X}}}
\newcommand{\eg}{\mbox{{\bf g}}}
\newcommand{\paromega}{\partial \Omega}
\newcommand{\gau}{\Gamma_{u}}
\newcommand{\gaf}{\Gamma_{f}}
\newcommand{\sig}{{\bf \sigma}}
\newcommand{\gac}{\Gamma_{\mbox{{\bf c}}}}
\newcommand{\deu}{\dot{\eu}}
\newcommand{\dueu}{\underline{\deu}}
\newcommand{\dev}{\dot{\ev}}
\newcommand{\duev}{\underline{\dev}}
\newcommand{\weak}{\rightharpoonup}
\newcommand{\weakdown}{\rightharpoondown}
\newcommand{\opg}{\stackrel{\mathfrak{g}}{\cdot}}
\newcommand{\opn}{\stackrel{\mathfrak{n}}{\cdot}}
\renewcommand{\contentsname}{ }

\title{Sub-Riemannian geometry and Lie groups \\ 
  Part I \\ 
Seminar Notes, 
DMA-EPFL, 2001}

\author{M. Buliga \\
 \\ 
Institut Bernoulli\\ 
B\^{a}timent MA \\
\'Ecole Polytechnique F\'ed\'erale de Lausanne\\
CH 1015 Lausanne, Switzerland\\
{\footnotesize Marius.Buliga@epfl.ch} \\ 
 \\ 
\and 
and \\ 
 \\ 
Institute of Mathematics, Romanian Academy \\ 
P.O. BOX 1-764, RO 70700\\ 
Bucure\c sti, Romania\\
{\footnotesize Marius.Buliga@imar.ro}}

\date{This version: October 31, 2002}

\maketitle


\newpage

\section*{Introduction}

M. Gromov \cite{gromov}, pages 85--86: 

"{\bf 3.15. Proposition:} Let $(V,g)$ be a Riemannian manifold with 
$g$ continuous. For each $v \in V$ the spaces $\lambda(V,v)$ 
Lipschitz converge as $\lambda \rightarrow \infty$ to the tangent
space $(T_{v} V, 0)$ with its Euclidean metric $g_{v}$. 

{\bf Proof$\left.\right._{+}$:} Start with a $C^{1}$ map 
$(R^{n},0) \rightarrow (V,v)$ whose differential is isometric at
$0$. The $\lambda$-scalings of this provide almost isometries
between large balls in $R^{n}$ and those in $\lambda V$ for $\lambda
\rightarrow \infty$. 

{\bf Remark:} In fact we can {\it define} Riemannian manifolds as
locally compact path metric spaces that satisfy the conclusion of 
Proposition 3.15. "

\vspace{.5cm}

If so, Gromov's remark should apply to any sub-Riemannian manifold.
Why then is the sub-Riemannian case so different from the Riemannian
one? Here is a list of legitimate questions: 

How can one define the manifold structure? Who are the tangent and
cotangent bundles? What is the intrinsic differential calculus? Why
are there  abnormal geodesics if the Hamiltonian formalism on the
cotangent bundle were complete? If the manifold is a compact Lie
group does the tangent bundle carry a natural group structure? What
are differential forms, de Rham cochain, and the variational complex? 
Consider the group of smooth volume preserving transformations. Why
does this group have more invariants than the volume and what is the
interpretation of these invariants? 

The purpose of this working seminar  was to explore as many as
possible open questions from the list above. Special attention has
been payed to the case of a Lie group with a left invariant
distribution.

The seminar, organised by the author and Tudor Ra\c tiu at the 
Mathematics Department, EPFL, started in November 2001.

The paper is by no means self contained. For any unproved result it is indicated the place where a complete proof can be found. The choice of 
the proofs is rather psychological: some of them have a geometrical or 
mixed geometric-analytical meaning 
(like the proof of Hopf-Rinow theorem), 
others help to understand that familiar reasoning applies in apparently 
unfamiliar situations (like in the preparations for the Pansu-Rademacher theorem). Important results (Ball-Box theorem for example) have elementary 
proof in particular situations; this might help to better understand their 
meaning and also their strength.

These pages covers the expository talks given by the author during this seminar. However, this is the first part of three, in preparation,  dedicated to this subject. It covers, with mild modifications, an elementary introduction to the 
field. 

The second part will deal with the applications of the Local-to-Global principle for moment maps arising in  connection with Carnot groups. We prove, for example, that Kostant nonlinear convexity theorem 4.1 \cite{kos} can be seen as the sub-Riemannian version of the linear convexity theorem 8.2. 
{\it op. cit.} The key point is the redefinition of the tangent bundle of a Lie group endowed with a left invariant nonintegrable distribution, using a natural construction based on noncommutative derivative.  

The third part  is devoted 
to the study of some representations of  groups of  bi-Lipschitz maps on Carnot groups. The main idea is that one can formulate an interesting rectifiability theory as a theory of (irreducible) representations of such groups.

{\bf Acknowledgements.} 
I would like 
to express here  my gratitude to Tudor Ra\c tiu for his continuous 
support and for the subtle ways of sharing his wide mathematical views. 
I want also to thank to the participants of this seminar and  generally  
for the warm athmosphere that one can encounter at the Mathematics Department 
of EPFL. I cannot close the aknowledgements without mentioning how much I benefit from the every-day impetus  given by a force of nature near me, which is my wife Claudia.

\newpage

\tableofcontents

\newpage

\section{From metric spaces to Carnot groups}

\subsection{Metric spaces}

\begin{defi} (distance) A function $d: X \times X \rightarrow [0,+
\infty)$ is a distance on $X$ if: 
\begin{enumerate}
\item[(a)] $d(x,y)=0$ if and only if $x=y$.
\item[(b)] for any $x,y$ $d(x,y)=d(y,x)$. 
\item[(c)] for any $x,y,z$ $d(x,z) \leq d(x,y) + d(y,z)$. 
\end{enumerate}
$(X,d)$ is  called a metric space. The open ball with centre $x \in X$ and 
radius $r>0$ is denoted by $B(x,r)$. 
\label{defdis}
\end{defi}

If $d$ ranges in $[0,+\infty]$ then it is called a pseudo-distance. 
The property of two points of being at finite distance is an 
equivalence relation. $d$ is then a distance on each equivalence 
class (leaf).

\begin{defi}
A map between metric spaces $f: X \rightarrow Y$ is Lipschitz if there is a 
positive constant $C$ such that for any 
$x,y \in X$ we have 
$$d_{Y}(f(x), f(y)) \ \leq \ C \ d_{X}(x,y)$$
The least such constant is denoted by $Lip(f)$. 

The dilatation, or metric derivative, of a map $f: X \rightarrow Y$ between metric spaces,  in a point $u \in Y$ is 
$$ dil(f)(u) = \limsup_{\varepsilon \rightarrow 0} \ \sup 
 \left\{ 
\frac{d_{Y}(f(v), f(w))}{d_{X}(v,w)} \ : \ v \not = w \ , \ v,w \in B(u,\varepsilon)
 \right\}$$

The distorsion of a map $f:X \rightarrow Y$ is 
$$dis \ f \ = \ \sup \left\{ \mid d_{Y}(f(y), f(y')) - 
d_{X} (y,y') \mid \mbox{ : } y,y' \in X \right\}$$
\end{defi}
The name "metric derivative" is motivated by the fact that for any  derivable function $f: R \rightarrow R^{n}$ the dilatation is  
$dil(f)(t) \ = \ \mid \dot{f}(t) \mid$.

A curve is a function $f: [a,b] \rightarrow X$. The image of a curve is called path. Length measures paths. Therefore length does not depends on the reparametrisation of the path and it is additive with respect to concatenation of paths. 
 
In a metric space $(X,d)$ one can measure the length of curves in several ways.

The length 
of a  curve with $L^{1}$ dilatation $f: [a,b] \rightarrow X$ is 
$$L(f) = \int_{a}^{b} dil(f)(t) \mbox{ d}t$$

A different way to define a length of a curve is to consider its variation. 

The curve $f$ has bounded
variation if the quantity 
$$ Var(f) \ = \ \sup \left\{ \sum_{i=1}^{n-1} d(f(t_{i}),
f(t_{i+1})) \ \mbox{ : } a \leq t_{1} < ... < t_{n} \leq b
\right\}$$ 
(called variation of $f$) is finite. 

The variation of a curve and the length of a path , as defined
previously, do not agree in general. To see this consider the
following easy example: 
$f: [-1,1] \rightarrow R^{2}$, $f(t) \ = \ (t, \ sign(t))$. 
We have $Var(f) \ = \ 4$ and $L(f([-1,1]) =  2$. Another example: 
the Cantor staircase function is continuous, but
not Lipschitz. It has variation equal to 1 and length of the graph
equal to 2. 

Nevertheless, for 
Lipschitz functions, the two definitions agree. 

\begin{thm}
For each Lipschitz curve $f: [a,b] \rightarrow X$, we have 
$L(f) \ = \ Var(f)$. 
\label{t411amb}
\end{thm}

For the statement and the proof see  the second part of Theorem 4.1.1., Ambrosio \cite{ambrosio}. 

\paragraph{Proof.}
We prove the thesis by double inequality. 
$f$ is continuous therefore $f([a,b])$ is a compact metric space. Let 
$\left\{ x_{n} \mbox{ :  } n \in N \right\}$ be a dense sequence in 
$f([a,b])$. All the functions 
$$t \mapsto \phi_{n}(t) \ = \ d(f(t),x_{n})$$
are Lipschitz and $Lip(\phi_{n}) \leq Lip(f)$, because of the general property: 
$$Lip(f \circ g) \ \leq \ Lip(f) \ Lip(g)$$
if $f,g$ are Lipschitz. In the same way we see that : 
$$dil(f)(t) \ = \ \sup \left\{  dil(\phi_{n})(t)  \mbox{ : } n \in N \right\}$$

We have then, for $t<s$ in $[a,b]$: 
$$d(f(t), f(s)) \ = \ \sup \left\{ \mid d(f(t),x_{n}) - d(f(s),x_{n})\mid \mbox{ : } n \in N \right\} \ \leq \ $$
$$\leq \ \int_{s}^{t} dil(f)(\tau) \mbox{ d}\tau$$
From the definition of the variation we get $$Var(f) \ \leq \ L(f)$$
For the converse inequality let $\varepsilon > 0$ and $n \geq 2$ natural number such that $h = (b-a)/n \ < \varepsilon$. Set $t_{i} \ = \ a + i h$. Then 
$$\frac{1}{h} \int_{a}^{b-\varepsilon} d(f(t + h), f(t)) \mbox{ d}t \ \leq \ 
\frac{1}{h} \sum_{i=0}^{n-2} d(f(\tau + t_{i+1}), f(\tau + t_{i})) \mbox{ d} \tau \ \leq $$ 
$$\leq \ \frac{1}{h} \int_{0}^{h}Var(f) \mbox{ d}\tau \ = \ Var(f)$$
Fatou lemma and definition of dilatation lead us to the inequality 
$$\int_{a}^{b-\varepsilon} dil(f) (t) \mbox{ d}t \ \leq \ Var(f)$$
This finishes the proof because $\varepsilon$ is arbitrary. 
\quad $\blacksquare$

There is a third, more basic way to introduce the length of a curve in a metric
space.

The length of the path $A = f([a,b])$ is by definition the
one-dimensional Hausdorff measure of the path. The definition is the
following: 
$$l(A) \ = \ \lim_{\delta \rightarrow 0}  
 \inf \left\{ \sum_{i \in I} diam \ E_{i}  \mbox{ : } diam \ E_{i} 
< \delta \ , \ \ A \subset \bigcup_{i \in I} E_{i} \right\} $$

The general Hausdorff measure is defined further: 

\begin{defi}
Let $k > 0$ and $(X,d)$ be a metric space. The $k$-Hausdorff measure 
is defined by: 
$$\mathcal{H}^{k}(A)  \ = \ \lim_{\delta \rightarrow 0}  
 \inf \left\{ \sum_{i \in I} \left(diam \ E_{i}\right)^{k}  \mbox{ : } diam \ E_{i} 
< \delta \ , \ \ A \subset \bigcup_{i \in I} E_{i} \right\} $$
Any $k$-Hausdorff measure is an outer measure. The Hausdorff
dimension of a set $A$ is defined by: 
$$\mathcal{H}-dim \ A \ = \ \inf\left\{ k > 0 \mbox{ : }
\mathcal{H}^{k}(A) \ = \ 0 \right\}$$
\end{defi}

Variation of a curve and the Hausdorff measure of the associated path don't coincide. For example take a circle $S^{1}$ in $R^{n}$ parametrised such
that any point is covered two times. Then the variation is two times
the length.

\begin{thm} Suppose that $f: [a,b] \rightarrow X$ is an injective 
 Lipschitz function and $A \ = \ f([a,b])$. Then $l(A) = Var(f)$. 
\label{t441amb}
\end{thm}

For the statement and the proof see  
Theorem 4.4.1., Ambrosio \cite{ambrosio}. 

In order to give the proof of the theorem we need two things. The first is 
the geometrically obvious, but not straightforward to prove in this generality, 
Reparametrisation Theorem (theorem 4.2.1 Ambrosio \cite{ambrosio}). 

\begin{thm}
Any path $A \subset X$ with a Lipschitz parametrisation admits a reparametrisation $f: [a,b] \rightarrow A$ such that $dil(f)(t) = 1$ for almost any $t \in [a,b]$. 
\label{tp}
\end{thm}

The second result that we shall need is Lemma 4.4.1 Ambrosio \cite{ambrosio}. 

\begin{lema}
If $f: [a,b] \rightarrow X$ is continuous then 
$$\mathcal{H}^{1}(f([a,b]) \ \leq \ d(f(a),f(b))$$
\label{led}
\end{lema}

\paragraph{Proof.}
Let us consider the Lipschitz function 
$$\phi: X \rightarrow R \ , \ \ \phi(x) \ = \ d(x, f(a))$$ 
It has the property $Lip(\phi) \leq 1$ therefore by the definition of Hausdorff measure we have 
$$\mathcal{H}^{1}(\phi \circ f ([a,b])) \ \leq \ \mathcal{H}^{1}(f([a,b])$$
On the left hand side of the inequality we have the Hausdorff measure on $R$, 
which coincides with the usual (outer) Lebesgue measure. Moreover 
$\phi \circ f ([a,b]) \ = \ [0,\alpha]$, therefore we obtain
$$\mathcal{H}^{1}(\phi \circ f ([a,b])) \ =  \ \sup \left\{ d(f(t),f(a)) \mbox{ : } t \in [a,b] \right\} \ \geq \ d(f(a), f(b))$$
\quad $\blacksquare$

The proof of theorem \ref{tp} follows. 

\paragraph{Proof.}
It is not restrictive to suppose that $A \ = \ f([a,b])$ can be parametrised by $f: [a,b] \rightarrow A$ such that $dil(f)(t) \ = \ 1$ for all $t \in [a,b]$. Due to  
theorem \ref{t411amb} and again  the reparametrisation theorem, we can choose 
$[a,b] = [0, Var(f)]$. 

For an arbitrary  $\delta > 0$ we choose $n \in N$ such that 
$h = Var(f)/n < \delta$ and we divide the interval $[0,Var(f)]$ in 
intervals $J_{i} \ = \ [ih, (i+1)h]$. The function $f$ has Lipschitz constant 
equal to 1 therefore (see definition of Hausdorff measure and notations therein) 
$$\mathcal{H}^{1}_{\delta} (A) \ \leq \ \sum_{i=0}^{n} diam \ (J_{i}) \ = \ 
Var(f)$$
$\delta$ is arbitrary therefore $\mathcal{H}^{1}(A) \ \leq \ Var(f)$. This is a general inequality which does not use the injectivity hypothesis. 

We prove the converse inequality from injectivity  hypothesis. Let us divide 
the interval $[a,b]$ by $a \leq t_{0} < ... < t_{n} \leq b$. From lemma 
\ref{led} and sub-additivity of Hausdorff measure we have: 
$$\sum_{i = 0}^{n-1} d(f(t_{i}), f(t_{i+1})) \ \leq \ \sum_{i = 0}^{n-1} 
\mathcal{H}^{1}(f([t_{i}), t_{i+1}])) \ \leq \ H^{1}(A)$$ 
The partition of the interval was arbitrary, therefore $Var(f) \leq \mathcal{H}^{1}(A)$. 
\quad $\blacksquare$

In conclusion, in any metric space we can measure length of
Lipschitz curves following one of the recipes from above, or we can
measure length of paths with the one-dimensional Hausdorff measure.
The lengths agree if the path admits a Lipschitz parametrisation.

We shall denote by $l_{d}$ the length functional, defined only on
Lipschitz curves, induced by the distance $d$.

The length induces a new distance $d_{l}$, say on any Lipschitz connected component of the space $(X,d)$. The distance $d_{l}$ is given by: 
$$d_{l}(x,y) \  = \ \inf \ \left\{ l_{d}(f([a,b])) \mbox{ : } f: [a,b] \rightarrow X \ \mbox{ Lipschitz } , \right.$$
$$\left. \ f(a)=x \ , \ f(b) = y \right\}$$

We have therefore  two
operators $d \mapsto l_{d}$ and $l \mapsto d_{l}$. Is one the
inverse of another? The answer is no. This leads to the introduction
of path metric spaces. 

\begin{defi}
A path metric space is a metric space $(X,d)$ such that $d  = d_{l}$. 
\label{dpath}
\end{defi}

In terms of distances there is an easy criterion to decide if a
metric space is path metric (Theorem 1.8., page 6-7, Gromov
\cite{gromov}).

\begin{thm}
A complete metric space is path metric if and only if (a) or (b)
from above is true: 
\begin{enumerate}
\item[(a)] for any $x,y \in X$ and for any $\varepsilon > 0$ there
is $z \in X$ such that 
$$\max \left\{ d(x,z), d(z,y) \right\} \ \leq \  \frac{1}{2} d(x,y)
+ \varepsilon$$
\item[(b)] for any $x,y \in X$ and for any $r_{1}, r_{2} > 0$, 
if $r_{1} + r_{2} \leq d(x,y)$ then 
$$d(B(x,r_{1}) , B(y,r_{2})) \ \leq \ d(x,y) - r_{1} - r_{2}$$
\end{enumerate}
\label{tcrtpath}
\end{thm}

\paragraph{Proof.}
We shall prove only that (a) implies that $X$ is path metric. (b)
implies (a) is straightforward, path metric space implies (b)
likewise. 

Set $\delta = d(x,y)$ and take a sequence $\varepsilon_{k} > 0$, with
finite sum. We shall recursively define a function 
$z = z_{\varepsilon}$ from the dyadic numbers in $[0,1]$ to $X$. 
$z_{\varepsilon}(1/2) \ = \ z_{1/2}$ is a point such that 
 $$\max \left\{ d(x,z_{1/2}), d(z_{1/2},y) \right\} \ \leq \  
\frac{1}{2} \delta (1 + \varepsilon_{1}) $$
Suppose now that all points $z_{\varepsilon} (p/2^{n}) \ = \
z_{p/2^{n}}$ were defined, for $p = 1, ... , 2^{n} -1$. Then 
$z_{2p+1/2^{n+1}}$ is a point such that 
$$\max \left\{ d(z_{p/2^{n}},z_{2p+1/2^{n+1}}), d(z_{2p+1/2^{n+1}}, 
z_{p+1/2^{n}}) \right\} \ \leq \  
\frac{\delta}{2^{n}}  \prod_{k=1}^{n+1} (1 + \varepsilon_{k}) $$
Because $(X,d)$ is a complete metric space it follows that
$z_{\varepsilon}$ can be prolonged to a Lipschitz curve
$c_{\varepsilon}$, defined on the whole interval $[0,1]$, such that 
$c(0) = x$, $c(1) = y$ and 
$$d(x,y) \ \leq \ l(c_{\varepsilon} \ \leq \  d(x,y) \ 
\prod_{k \geq 1} (1 + \varepsilon_{k})$$ 
But the product $\prod_{k \geq 1} (1 + \varepsilon_{k})$ can be made
arbitrarily close to $1$, which proves the thesis. 
\quad $\blacksquare$

Path metric spaces are geometrically more interesting, because one can define 
 geodesics. Let $c$ be any curve; denote the length of the 
restriction of $c$ to the interval $t,t'$ by $l_{c}(t,t')$. 

\begin{defi}
A (local) geodesic is a curve $c: [a,b] \rightarrow X$ with the property that for any $t \in (a,b)$ there is a small $\varepsilon > 0$ such that $c: [t - \varepsilon, t+\varepsilon] \rightarrow X$
is length minimising. A global geodesic is a length minimising curve.
\end{defi}

Therefore in a path metric space a local geodesic has the property  that in the neighbourhood 
of any of it's points the relation 
$$d(c(t),c(t')) = l_{c}(t,t')$$ holds. 
 Any global geodesic is also local geodesic. 

Can one join any two points with a geodesic?

The abstract Hopf-Rinow theorem (Gromov \cite{gromov}, page 9) states that: 

\begin{thm}
If $(X,d)$ is a  connected locally compact path metric space then 
each pair of points can be joined by a global geodesic.
\label{thr} 
\end{thm}

\paragraph{Proof.}
It is sufficient to give the proof for compact path metric spaces. Given the points 
$x,y$, there is a sequence of curves $f_{h}$ joining those points such that 
$l(f_{h}) \leq d(x,y) + 1/h$. The sequence, if parametrised by arclength, is 
equicontinuous; by Arzela-Ascoli theorem one can extract a subsequence (denoted 
also $f_{h}$) which converges uniformly to $f$. By construction the length function is lower semicontinuous hence: 
$$l(f) \leq \liminf_{h \rightarrow \infty} l(f_{h}) \leq d(x,y)$$
Therefore $f$ is a length minimising curve joining $x$ and $y$. 
\quad $\blacksquare$

The following remark of Gromov 1.13(b) \cite{gromov} easily comes from the examination 
of the previous proof. 

\begin{prop}
In a compact path metric space every free homotopy class is represented by a length minimising curve. 
\label{corg}
\end{prop}

\subsection{Local to Global Principle}

Let $(X,d)$ be a connected locally compact and complete path metric space. 

The natural notion of subspace in this category of metric spaces 
is the following one. 

\begin{defi}
A set $C \subset X$ is a closed convex set if the inclusion 
$(C,d_{C}) \subset (X,d)$ is an isometry, where $d_{C}$ is the 
inner path metric distance: 
$$d_{C}(x,y) \  = \ \inf \ \left\{ l_{d}(f([a,b])) \mbox{ : } f: [a,b] \rightarrow X \ \mbox{ Lipschitz } , \right. $$ 
$$ \left. \ f(a)=x \ , \ f(b) = y \ , \ 
f([a,b]) \subset C \  \right\}$$

We shall call a closed set weakly convex if any two points of the set can be joined by a local geodesic lying in the set. 
\end{defi}

A closed set it is therefore convex if any two points $x,y \in C$ can be joined 
by a global geodesic lying in $C$. Note that there might be other 
global geodesics joining two points in $C$ which are not entirely in $C$.

The purpose of the remainder of this section is to state and prove the 
Local-to-Global Principle. The principle has been formulated for the first time 
by Hilgert, Neeb, Plank \cite{hnp}. It has been used to prove convexity theorems; it's most widely known application concerns an alternative proof 
of the Atiyah-Guillemin-Sternberg theorem about convexity of the image of the moment map. We shall prove in this section a general version of the principle, true in path metric spaces. Applications of this principle will be given further, where we will meet also the moment map.

\begin{defi}
A pointed closed convex set is a pair $(C,x)$ with $x \in C \in X$ and 
$C$ closed convex set with the property that  there is a constant $K \geq 1$ 
such that for any $r>0$ and any two points $y,y' \in C \cap B(x,r)$ can be joined by a global geodesic in $C \cap B(x, K r)$. 
\end{defi}

Let $X$ be a connected, compact topological Hausdorff space,  $Y$ a connected 
locally compact path metric space and 
$f: X \rightarrow Y$ with the following properties: 
\begin{enumerate}
\item[LC1.] $f$ is proper and  for any $x \in X$ there is a neighbourhood $V_{x}$ of $x$ such that 
$f^{-1}(f(y)) \cap V_{x}$ is connected for all $y \in V_{x}$. 
\item[LC2.] for any $x \in X$ there is a pointed closed convex  $(C_{x}, f(x))$ 
 and a neighbourhood $V_{x}$ of $x$ such that $f: V_{x} \rightarrow C_{x}$ is open. 
\end{enumerate}

\begin{thm}(Local-to-Global Principle). Under the hypotheses LC1, \\ 
LC2, the image 
$f(X)$ is closed weakly convex in $Y$. 
\label{tltg}
\end{thm}

The proof is inspired by the one from Sleewaegen \cite{sleewa}, section 1.8, Theorem 1.8.1.

\paragraph{Proof.}
The proof splits in two steps, one of topological nature and the other of 
 metric character. 

{\bf Step 1.}
Let us consider on $X$ the equivalence relation $x \approx y$ if $f(x) = f(y)$ and $x,y$ are in the same connected component of $f^{-1}(f(x))$. Denote 
by $\tilde{X}$ the set of equivalence classes, endowed with the quotient topology, and by $\pi: X \rightarrow \tilde{X}$ the canonical projection. 
Let $\tilde{f}: \tilde{X} \rightarrow N$ be the map with the property 
$\tilde{f}\circ \pi \ = \ f$. 

The property  LC2 is hereditary:  it is  invariant to the choice of arbitrary small neighbourhoods $V_{x}$. 
We obtain the following proposition (see Lemmas 1.8.1, 1.8.2. Sleewaegen \cite{sleewa})
using topological reasoning and the hereditarity of LC2: 

\begin{lema}
The space $\tilde{X}$ is Hausdorff, connected and compact. The function 
$\tilde{f}$ is proper, the preimage $\tilde{f}^{-1}(y)$ has a finite number of 
elements and it satisfies the improved version of LC2: for any $\tilde{x} \in \tilde{X}$ there are a neighbourhood $V_{\tilde{x}}$ of $\tilde{x}$, a pointed 
closed convex set  $(C_{\tilde{x}}, \tilde{f}(\tilde{x}))$ and $r_{\tilde{x}} > 0$ such that the restriction 
$\tilde{f} : V_{\tilde{x}} \rightarrow C_{\tilde{x}} \cap B(f(\tilde{x}), r_{\tilde{x}})$ is a homeomorphism on the image.  
\end{lema}

{\bf Step 2.}
Denote by $l$ the length functional in the compact set $f(X)$ and by 
$d_{l}$ the induced length distance. We begin by constructing a distance on 
$\tilde{X}$, which comes from the distance $d_{l}$ and the function $\tilde{f}$. 

For any two points $\tilde{x}, \tilde{y} \in \tilde{X}$ let 
$$\tilde{d}(\tilde{x},\tilde{y}) \ = \ \inf \left\{ l(\tilde{f} \circ \gamma) 
\mbox{ : } \gamma: [a,b] \rightarrow \tilde{X} \mbox{ continuous, } 
\gamma(a) = \tilde{x} \ , \ \gamma(b) = \tilde{y} \right\}$$
We have the following immediate inequalities: 
$$\tilde{d}(\tilde{x},\tilde{y}) \ \geq \ d_{l}(\tilde{f}(\tilde{x}), 
\tilde{f}(\tilde{y})) \ \geq \ d(\tilde{f}(\tilde{x}), 
\tilde{f}(\tilde{y}))$$
In order to check that $d_{l}$ is a distance it suffices to show that 
if $\tilde{d}(\tilde{x},\tilde{y}) \ = \ 0$ then $\tilde{x} = \tilde{y}$. 
Indeed, we have $\tilde{f}(\tilde{x}) \ = \ \tilde{f}(\tilde{y})$ and there are 
two possibilities. If $\tilde{y} \in V_{\tilde{x}}$ (neighbourhood chosen as 
in the improved version of LC2) then $\tilde{y} = \tilde{x}$. Otherwise, if 
$\tilde{y} \not \in V_{\tilde{x}}$ then we have arbitrary short closed  curves 
in $Y$ which start and end at $\tilde{f}(\tilde{x})$ but which go out of 
$C_{\tilde{x}} \cap B(f(\tilde{x}), r_{\tilde{x}})$, which is in contradiction 
with $r_{\tilde{x}} > 0$. 

Use now Gromov Theorem \ref{tcrtpath} to prove that $(\tilde{X},\tilde{d})$ is  
a compact path metric space. Hopf-Rinow Theorem \ref{thr} ensures us that any 
two points $\tilde{x}, \tilde{y} \in \tilde{X}$ can be joined by a global minimising geodesic. Moreover, such a geodesic projects on a global 
$d_{l}$ geodesic in $f(X)$, because of the construction of the distance $\tilde{d}$.

Because $f(X) \ = \tilde{f}(\tilde{X})$ is compact, there is a $K \geq 1$ such that  we can cover $f(X)$  with a finite number of sets of the form $C_{\tilde{x}^{k}} \cap B(\tilde{f}(
\tilde{x}^{k}), r_{k})$ and $C_{\tilde{x}^{k}} \cap B_{N}(\tilde{f}(
\tilde{x}^{k}), K r_{k}) \ \subset \ f(X)$. Then any $d_{l}$ global geodesic 
is made by gluing a finite number of global geodesics in $Y$ and the question 
is if such a curve is local geodesic in the neighbourhood of the gluing points. But the improved LC2 
property ensures that. Indeed, take any two sufficiently closed points 
$\tilde{f}(\tilde{x})$ and $\tilde{f}(\tilde{y})$ on a $d_{l}$ geodesic $\tilde{f}\circ \gamma$, with $\gamma$ global $\tilde{d}$ geodesic. Then 
there is a neighbourhood $V_{\tilde{x}}$ of $\tilde{x}$ such that 
$\tilde{y} \in V_{\tilde{x}}$ and $\tilde{f} \circ \gamma$ 
 is a $Y$ global geodesic in $f(X)$ 
which joins $\tilde{f}(\tilde{x})$ and $\tilde{f}(\tilde{y})$.  
\quad $\blacksquare$

The compactness assumption on the space $X$ is not important. For example if the constant $K$ in the definition of closed convex pointed set can be 
chosen independent of $x$ in LC2 then the compactness assumption is no longer 
needed. 

The classical Local-to-Global Principle is stated for vectorial target space 
$Y$, hence an arbitrary Euclidean space. Such spaces admit unique geodesic 
between any two different points, therefore the conclusion of the principle 
is that the image $f(X)$ is convex. 

It seems hard to check that the function $f$ satisfies LC2. In the case of 
the moment map this is a consequence of the existence of roots for some 
unitary representations. Note here the appearance of group theory. 

\subsection{Distances between  metric spaces} 

The references for this section  are Gromov \cite{gromov}, chapter 3, and Burago \& al. \cite{burago} section 7.4.  There are several definitions of distances between metric spaces. The very fertile idea of introducing such distances belongs to Gromov. 

In order to introduce the Hausdorff distance between metric spaces, recall 
the Hausdorff distance between subsets of a metric space. 

\begin{defi}
For any set $A \subset X$ of a metric space and any $\varepsilon > 0$ set 
the $\varepsilon$ neighbourhood of $A$ to be 
$$A_{\varepsilon} \ = \ \cup_{x \in A} B(x,\varepsilon)$$
The Hausdorff distance between $A,B \subset X$ is defined as
$$d_{H}^{X}(A,B) \ = \ \inf \left\{ \varepsilon > 0 \mbox{ : } A \subset B_{\varepsilon} \ , \ B \subset A_{\varepsilon} \right\}$$
\end{defi} 

By considering all isometric embeddings of two metric spaces $X$, $Y$ into 
an arbitrary metric space $Z$ we obtain the Hausdorff distance between $X$, $Y$ (Gromov \cite{gromov} definition 3.4). 

\begin{defi}
The Hausdorff distance $d_{H}(X,Y)$ between metric spaces $X$ $Y$ is the infimum of the numbers 
$$d_{H}^{Z}(f(X),g(Y))$$
for all isometric embeddings $f: X \rightarrow Z$, $g: Y \rightarrow Z$ in a 
metric space $Z$. 
\end{defi}

If $X$, $Y$ are compact then $d_{H}(X,Y) < + \infty$. Indeed, let $Z$ be the disjoint union of $X,Y$ and $M \ = \ \max \left\{ diam(X) , diam(Y) \right\}$. 
Define the distance on $Z$ to be 
$$d^{Z}(x,y) \ = \ \left\{ \begin{array}{ll} 
d^{X}(x,y) & x, y \in X \\ 
d^{Y}(x,y) & x, y \in Y \\ 
\frac{1}{2} M & \mbox{ otherwise} 
\end{array} \right. $$
Then $d_{H}^{Z}(X,Y) < + \infty$. 

The Hausdorff distance between isometric spaces equals $0$. The converse is also true (Gromov {\it op. cit.} proposition 3.6) in the class of compact metric spaces.

\begin{thm}
If $X,Y$ are compact metric spaces such that $d_{H}(X,Y) = 0$ then 
$X,Y$ are isometric. 
\label{tgro}
\end{thm}

For the proof of the theorem we need the Lipschitz distance ({\it op. cit. definition 3.1}) and a criterion for convergence of metric spaces in the Hausdorff distance ( {\it op. cit. proposition 3.5 }). We shall give the definition of Gromov for Lipschitz distance and the first part of the 
mentioned proposition. 

\begin{defi} 
The Lipschitz distance $d_{L}(X,Y)$ between bi-Lipschitz homeomorphic metric spaces 
$X,Y$ is the infimum of 
$$\mid \log \ dil(f) \mid + \mid \log \ dil(f^{-1}) \mid$$
for all $f: X \rightarrow Y$, bi-Lipschitz homeomorphisms. 
\end{defi}

Obviously, if $d_{L}(X,Y) = 0$ then $X,Y$ are isometric. Indeed, by definition we have a sequence $f_{n}: X \rightarrow Y$ such that $dil (f_{n}) \rightarrow 1$ as $n \rightarrow \infty$. Extract an uniformly convergent subsequence; the limit is an isometry. 

\begin{defi}
A $\varepsilon$-net in the metric space $X$ is a set $P \subset X$ such that 
$P_{\varepsilon} \ = \ X$. The separation of the net $P$ is 
$$sep(P) \ = \ \inf \left\{ d(x,y) \mbox{ : } x \not = y \ , \ x,y \in P \right\}$$
A $\varepsilon$-isometry between $X$ and $Y$ is a function $f: X \rightarrow Y$ such that $dis \ f \ \leq \varepsilon$ and $f(X)$ is a $\varepsilon$ net in 
$Y$. 
\end{defi}

The following proposition gives a connection between  convergence of metric spaces 
in the Hausdorff distance and  convergence of $\varepsilon$ nets in the Lipschitz distance.  

\begin{defi} 
A sequence of metric spaces $(X_{i})$ converges in the sense of Gromov-Hausdorff to the metric space $X$ if $d_{H}(X,X_{i}) \rightarrow 0$ as 
$i \rightarrow \infty$. This means that there is a sequence $\eta_{i} \rightarrow 0$ and isometric embeddings $f_{i}: X \rightarrow Z_{i}$, 
$g_{i} : X_{i} \rightarrow Z_{i}$ such that 
$$d_{H}^{Z_{i}} (f_{i}(X), g_{i}(X_{i})) \ < \ \eta_{i}$$
\label{dhauc}
\end{defi}

\begin{prop}
Let $(X_{i})$ be a sequence of metric spaces converging to $X$ and 
$(\eta_{i})$ as in definition \ref{dhauc}. Then for any $\varepsilon > 0$ and 
any $\varepsilon$-net $P \subset X$  with strictly positive separation there is a sequence $P_{i} \subset X_{i}$ such that 
\begin{enumerate} 
\item[1.] $P_{i}$ is  a $\varepsilon + 2  \eta_{i}$ net in $X_{i}$, 
\item[2.] $d_{L}(P_{i}, P) \rightarrow 0$ as $i \rightarrow \infty$, 
uniformly with respect to $P, P_{i}$. 
\end{enumerate}
\end{prop}

We shall not use further this proposition, given for the sake of completeness. 

The next proposition is corollary 7.3.28 (b), Burago \& al. \cite{burago}. The proof is adapted from Gromov, proof of proposition 3.5 (b). 

\begin{prop}
If there exists a $\varepsilon$-isometry between $X,Y$ then $d_{H}(X,Y) < 2 \varepsilon$. 
\label{pbur}
\end{prop}

\paragraph{Proof.}
Let $f: X \rightarrow Y$ be the $\varepsilon$-isometry. On the disjoint union 
$Z = X \cup Y$ extend the distances $d^{X}$, $d^{Y}$ in the following way. Define the distance between $x \in X$ and $y \in Y$ by 
$$d(x,y) \ = \ \inf\left\{ d^{X}(x,u) + d^{Y}(f(u),y) + \varepsilon \right\}$$
This gives a distance $d^{Z}$ on $Z$. Check that 
$d^{Z}_{H}(X,Y) < 2 \varepsilon$. 
\quad $\blacksquare$

The proof of theorem \ref{tgro} follows as an application of previous propositions.

\subsection{Metric tangent  cones}

The local geometry of a path metric space is described with the help of Gromov-Hausdorff convergence of pointed metric spaces. 

This is definition 3.14 Gromov \cite{gromov}. 
\begin{defi}  
The sequence of pointed metric spaces $(X_{n}, x_{n}, d_{n})$ converges in the sense 
of Gromov-Hausdorff to the pointed space $(X,x, d_{X})$ if 
for any $r>0$, $\varepsilon > 0$ there is $n_{0} \in N$ such that for all 
$n \geq n_{0}$ there exists $f_{n} : B_{n}(x_{n},r) \subset X_{n} \rightarrow X$ such that: 
\begin{enumerate}
\item[(1)] $f_{n}(x_{n}) \ = \ x$, 
\item[(2)] the distorsion of $f_{n}$ is bounded by $\varepsilon$: $dis \ f_{n} \ < \ \varepsilon$, 
\item[(3)]
$B_{X}(x,r) \ \subset \ \left(f_{n}(B_{n}(x_{n},r)) \right)_{\varepsilon}$. 
\end{enumerate}
\end{defi}

The Gromov-Hausdorff limit is defined up to isometry, in the class of compact metric spaces (Proposition 3.6 Gromov \cite{gromov}) or in the class of locally compact cones. Here it is the definition of a cone.

\begin{defi}
A pointed metric space $(X,x_{0})$ is called a cone if for any $\lambda >0$ 
there is a map $\delta_{\lambda}: X \rightarrow X$ such that $\delta_{\lambda}(x_{0}) \ = 
\ x_{0}$ and for any $x,y \in X$ 
$$d(\delta_{\lambda}(x), \delta_{\lambda}(y)) \ = \ \lambda \ d(x,y)$$
Such a map is called a dilatation with center $x_{0}$ and coefficient $\lambda$. 
\end{defi}

The local geometry of a metric space $X$ in the neighbourhood of $x_{0} \in X$ is described by the tangent space to $(X,x_{0})$ (if such object exists). 

\begin{defi}
The tangent space to $(X,x_{0})$ is the Gromov-Hausdorff limit 
$$(T_{x_{0}}, 0, d_{x_{0}})  \ = \ \lim_{\lambda \rightarrow \infty} (X,x_{0}, 
\lambda d)$$
\end{defi}

Three remarks are in order: 
\begin{enumerate}
\item[1.] The tangent space is obviously a cone. We shall see that in a large class of situations is also a group, hence a graded nilpotent group, called 
for short Carnot group.  
\item[2.] The tangent space comes with a metric inside. This space is 
path metric (Proposition 3.8 Gromov \cite{gromov}). 
\item[3.] The tangent cone is defined up to isometry, therefore there is 
no way to use the metric tangent cone definition to construct a tangent bundle. 
For a modification of the definitions in this direction  see Margulis, Mostow \cite{marmos2}.  
\end{enumerate}

One can ask for a classification of metric spaces which admit everywhere the same (up to isometry) tangent cone. Such classification results exist, for example in the category of  Lipschitz $N$ manifold.

\begin{defi}
Let $N$ be a cone with a path metric distance. A Lipschitz $N$ manifold is an 
(equivalence class of, or a maximal) atlas over the cone $N$ such that the 
change of charts is locally bi-Lipschitz. 
\label{dlipm}
\end{defi}

We shall see in section \ref{rigid} a rigidity (i.e. classification result)  property of manifolds. 
Also, we shall prove in section \ref{manif} that a sub-Riemannian manifold (see definition 
\ref{defsr}) which is also a Lie group does not admit  
a manifold structure over the nilpotentisation of the defining distribution, such that the operation and the group exponential to be smooth. We will be forced therefore  to  modify the definition of a manifold in order to show that a group with a left invariant distribution has 
a manifold structure over the nilpotentisation. The modification will lead, surprising but 
natural, to (a generalisation of) the moment map. 

\subsection{Examples of path metric spaces}
\label{exemple}

In this section are given examples of path metric  spaces. These are: 

\begin{enumerate}
\item[-] Riemannian manifolds
\item[-] Finsler manifolds
\item[-] Carnot-Carath\'eodory manifolds
\end{enumerate}

 Let $X$ be a Riemannian manifold and $l(f)$ be 
the length of the curve $f$ with respect to the metric on $X$, if $f$ is 
piecewise $C^{1}$, otherwise  $l(f) = + \infty$. This is, off course, 
{\it a} metric structure on the Riemannian $X$. The class of curves with 
{\it potential finite} length can be enlarged or restricted by using 
analytical arguments.

Change the metric on $X$ by a Finsler metric, i.e. consider a continuous mapping 
$\Delta: TX \rightarrow R$ such that for any $x \in X$, $\lambda \in R$ and 
$v \in T_{x}X$ 
$$\Delta_{x}(\lambda v) = \mid \lambda \mid \Delta_{x}(v)$$
Then for any piecewise $C^{1}$ curve $c: [a,b] \rightarrow X$ define it's length 
by 
$$l(c) = \int_{a}^{b} \Delta_{c(t)}(\dot{c}(t)) \mbox{ d}t$$
This give rise to a (non-oriented) length structure on $X$.

Instead of a Finsler metric $\Delta$ consider an oriented version obtained 
by imposing a milder condition upon $\Delta$: to be positively one-homogeneous. 
That means: for any $x \in X$, $\lambda \geq 0$ and $v \in T_{x}X$ 
$$\Delta_{x}(\lambda v) =  \lambda \Delta_{x}(v)$$
The length structure associated to $\Delta$ is now oriented. 

A particular case is the one of a Finsler-Minkowski metric. Suppose first that 
$(X,g)$ is Riemannian ($g$ is the metric).  Take $Q_{x} \subset 
T_{x}X$ a convex bounded set and denote by $\chi_{x}$ the characteristic function 
$$\chi_{x}(v) \ = \ \left\{ \begin{array}{ll}
0 & \mbox{ if } v \in Q_{x} \\ 
+\infty & \mbox{otherwise}
\end{array} \right.$$
Suppose that $Q_{x}$ contains a ball centered in the origin. Define now the 
Finsler-Minkowski metric  
$$\Delta_{x}(v) \ = \ \sup \left\{ g(v,u) - \chi_{x}(u) \ : \ u \in T_{x}X \right\}$$
that is $\Delta_{x}$ is the polar of $\chi_{x}$. From the hypothesis, there 
are positive constants $c,c'$ such that ($\mid v \mid_{x} = \sqrt{g(v,v)}$) 
$$c \mid v \mid_{x} \leq \Delta_{x}(v) \leq c' \mid v \mid_{x}$$
The length structure associated to $\Delta$ is in  general oriented. It becomes 
non oriented if the set  $Q_{x}$ is 
symmetric with respect to the origin.

In almost all examples  we have a Riemannian manifold and a length 
structure constructed from a map $\Delta$. We shall suppose now that the manifold 
is complete with respect to the Riemannian distance.

\begin{prop}
The Hopf-Rinow theorem holds on a Finsler-Minkowski manifold, provided that 
it is complete with respect to the initial Riemannian distance and that 
$\Delta$ is convex. 
\label{pfm}
\end{prop}

A more particular case is the following: 

\begin{prop}
Let $G$ be a Lie group and $U$ be a left invariant Finsler-Minkowski metric on   
$G$. Then $(G,d_{U})$ is complete and locally compact. 
\label{pmiel}
\end{prop}

\paragraph{Proof.}
The exponential map $\exp$ is a diffeomorphism from $\mathcal{V}(0) \ni U \subset 
g$ to $\exp U$. From the hypothesis we get that $\exp_{|_{U}}$ is $d_{U}$  Lipschitz. The exponential is also open. $G$ is locally a linear group, hence from 
the hypothesis there is a constant $c$ such that 
$$d_{U}(e, \exp v) \geq c \| v\|$$
Therefore the metric topology given by $d_{U}$ is the same as the manifold topology, hence locally compact. 

Take now a Cauchy sequence $(x_{h})$. For a sufficiently small $\varepsilon > 0$ 
there is $N = N(\varepsilon)$ such that $x_{N}^{-1}x_{h}$ lies in a compact neighbourhood of the identity element $e$ for any $h\geq N$. One can extract a convergent subsequence from it and the proof finishes. 
\quad $\blacksquare$ 

As a corollary of propositions \ref{pfm} and \ref{pmiel} we have the following concrete Hopf-Rinow theorem:

\begin{cor}
On a 
Lie group endowed with a convex left (or right) invariant Finsler-Minkowski 
metric any two points can be joined by a geodesic. 
\label{cor1}
\end{cor}

Another class of examples is provided by sub-Riemannian manifolds. 

\begin{defi} 
A sub-Riemannian (SR) manifold is a triple $(M,H, g)$, where $M$ is a 
connected manifold, $H$ is a subbundle of $TM$, named horizontal 
bundle or distribution, and $g$ is a metric (inner-product) on the 
horizontal bundle.  

A horizontal curve is a continuous, almost everywhere
differentiable curve, whose tangent lies in the horizontal bundle. 

The length of a horizontal curve $c: [a,b] \rightarrow M$ is 
$$l(c) \ = \ \int_{a}^{b} g(c(t), \dot{c}(t)) \mbox{ d} t$$

 The SR 
manifold is called a Carnot-Carath\'eodory (CC) space if any two
points can be joined by a finite length horizontal curve. 
\label{defsr}
\end{defi}

A CC space is a path metric space with the Carnot-Carath\'eodory
distance induced by the length $l$: 
$$d(x,y) \ = \ \inf \left\{ l(c) \mbox{ : } c: [a,b] \rightarrow M \
, \ c(a) = x \ , \  c(b) = y \right\} $$

The particular case that will be interesting for us is the following one. 
Consider a connected real Lie group  and a  vector space $V \subset g$, which generates the whole algebra $g$ by 
Lie brackets. Then we shall see that the left  translations of $V$ provide a non-integrable 
distribution $V$ on $G$ (easy form of Chow theorem). For any (left or right invariant) metric defined on $D$ 
we have an associated  Carnot-Carath\'eodory distance.

\section{Carnot groups}

\subsection{Structure of Carnot groups}
\begin{defi}
A Carnot (or stratified nilpotent) group is a
connected simply connected group $N$  with  a distinguished vectorspace 
$V_{1}$ such that the Lie algebra of the group has the  
direct sum decomposition: 
$$n \ = \ \sum_{i=1}^{m} V_{i} \ , \ \ V_{i+1} \ = \ [V_{1},V_{i}]$$
The number $m$ is the step of the group. The number 
$$Q \ = \ \sum_{i=1}^{m} i \ dim V_{i}$$
is called the homogeneous dimension of the group. 
\label{dccgroup}
\end{defi}

Because the group is nilpotent and simply connected, the
exponential mapping is a diffeomorphism. We shall identify the 
group with the algebra, if is not locally otherwise stated.

The structure that we obtain is a set $N$ endowed with a Lie
bracket and a group multiplication operation.

The Baker-Campbell-Hausdorff formula shows that the group operation
is polynomial. It is easy to see that the Lebesgue measure on the
algebra is (by identification) a bi-invariant measure.

We give further examples of such groups:

{\bf (1.)} $R^{n}$ with addition is the only commutative Carnot group. 

{\bf (2.)} The Heisenberg group is the first non-trivial example. 
This is the group $H(n) \ = \ R^{2n} \times R$ with the operation: 
$$(x,\bar{x}) (y,\bar{y}) \ = \ (x+y, \bar{x} + \bar{y} + 
\frac{1}{2} \omega(x,y) )$$
where $\omega$ is the standard symplectic form on $R^{2n}$. 
The Lie bracket is 
$$[(x,\bar{x}), (y,\bar{y})] \ = \ (0,\omega(x,y))$$
The direct sum decomposition of (the algebra of the) group is: 
$$H(n) \ = \ V + Z \ , \ \ V \ = \ R^{2n} \times \left\{ 0 \right\}
\ , \ \ Z \ = \  \left\{ 0 \right\} \times R$$ 
$Z$ is the center of the algebra, the group has step 2 and
homogeneous dimension $2n+2$. 

{\bf (3.)} H-type groups. These are two step nilpotent Lie groups $N$ 
endowed with an inner product $(\cdot , \cdot)$, such that the
following {\it orthogonal} direct sum decomposition occurs: 
$$N \ = \ V + Z$$
$Z$ is the center of the Lie algebra. Define now the function 
$$J : Z \rightarrow End(V) \ , \ \ (J_{z} x, x') \ = \ (z, [x,x'])$$ 
The group $N$ is of H-type if for any $z \in Z$ we have 
$$J_{z} \circ J_{z} \ = \  - \mid z \mid^{2} \  I$$
From the Baker-Campbell-Hausdorff formula we see that the group
operation is 
$$(x,z) (x', z') \ = \ (x + x', z + z' + \frac{1}{2} [x,x'])$$
These groups appear naturally as the nilpotent part in the Iwasawa
decomposition of a semisimple real group of rank one. (see \cite{cdkr})

{\bf (4.)} The last example is the group of $n \times n$ 
upper triangular matrices, which is
nilpotent of step $n-1$.  This example is important because any Carnot group is isomorphic with a subgroup of a group of upper triangular matrices.

Any Carnot group admits a one-parameter family of dilatations. For any 
$\varepsilon > 0$, the associated dilatation is: 
$$ x \ = \ \sum_{i=1}^{m} x_{i} \ \mapsto \ \delta_{\varepsilon} x \
= \ \sum_{i=1}^{m} \varepsilon^{i} x_{i}$$
Any such dilatation is a group morphism and a Lie algebra morphism. 

In fact the class of Carnot groups is characterised by the existence of dilatations. 

\begin{prop}
Suppose that the Lie algebra $\mathfrak{g}$ admits an one parameter group 
$\varepsilon \in (0,+\infty) \mapsto \delta_{\varepsilon}$ of simultaneously 
diagonalisable Lie algebra isomorphisms. Then $\mathfrak{g}$ is the algebra of 
a Carnot group. 
\end{prop}

\paragraph{Proof.}
The hypothesis means that there is a direct sum decomposition of 
$$\mathfrak{g} \ = \ \oplus_{i=1}^{m} V_{i}$$ such that for any 
$\varepsilon > 0$ we have 
$$x \ = \ \sum_{i=1}^{m} x_{i} \ \mapsto \ \delta_{\varepsilon} x \ = \ 
\sum_{i=1}^{m} f_{i}(\varepsilon) x_{i}$$
with $f_{i}$ continuous,  moreover 
$$\delta_{\varepsilon} \circ \delta_{\mu} \ = \ \delta_{\varepsilon \mu}$$
for any $\varepsilon, \mu > 0$ and $\delta_{1} \ = \ id$. From this we get 
$f_{i}(\varepsilon) \ = \ \varepsilon^{\alpha_{i}}$. Each $\delta_{\varepsilon}$ is also Lie algebra morphism, therefore ($x_{i} \in V_{i}$, $x_{j} \in V_{j}$) 
$$ \delta_{\varepsilon}[x_{i}, x_{j}] \ = \ [\delta_{\varepsilon} x_{i}, 
\delta_{\varepsilon} x_{j}] \ = \ \varepsilon^{\alpha_{i} + \alpha_{j}} [x_{i}, x_{j}]$$
We conclude that $\alpha_{i} + \alpha_{j}$ is also an eigenvalue, unless 
$[x_{i}, x_{j}] = 0$. The direct sum decomposition of $\mathfrak{g}$ is finite 
therefore $\alpha_{i} \ = \ i \alpha$ and, more important, $[V_{i}, V_{j}] \ = \ V_{i+j}$. In conclusion $\mathfrak{g}$ is the Lie algebra of a Carnot group 
and $\delta_{\varepsilon}$ are the dilatations, up to a scale factor 
$\alpha$. 
\quad $\blacksquare$

We can always find inner products on $N$ such that the 
decomposition $N \ = \ \sum_{i=1}^{m} V_{i}$ is an orthogonal sum.


We shall endow the group $N$ with a structure of a sub-Riemannian
manifold now. For this take the distribution obtained from left
translates of the space $V_{1}$. The metric on that distribution is
obtained by left translation of the inner product restricted to
$V_{1}$.

If $V_{1}$ Lie generates (the algebra) $N$
then any element $x \in N$ can be written as a product of 
elements from $V_{1}$.  An useful lemma is: 

\begin{lema}
Let $N$ be a Carnot group and $X_{1}, ..., X_{p}$ an orthonormal basis 
for $V_{1}$. Then there is a  
 a natural number $M$ and a function $g: \left\{ 1,...,M \right\} 
\rightarrow \left\{ 1,...,p\right\}$ such that any element 
$x \in N$ can be written as: 
\begin{equation}
x \ = \ \prod_{i = 1}^{M} \exp(t_{i}X_{g(i)})
\label{fp2.4}
\end{equation}
Moreover, if $x$ is sufficiently close (in Euclidean norm) to
$0$ then each $t_{i}$ can be chosen such that $\mid t_{i}\mid \leq C 
\| x \|^{1/m}$
\label{p2.4}
\end{lema}

\paragraph{Proof.}
This is a slight reformulation of Lemma 1.40, Folland, Stein \cite{fostein}. 
We have used Proposition \ref{pest} applied for the homogeneous norm 
$\mid \cdot \mid_{1}$ (see below).
\quad $\blacksquare$

 This means that there is a horizontal curve
joining any two points. Hence the distance 
$$d(x,y) \ = \ \inf \left\{ \int_{a}^{b} \| c^{-1} \dot{c} \| \mbox{ 
d}t \ \mbox{ : } \ c(a) = x , \ c(b) = y , \ c^{-1} \dot{c} \in
V_{1} 
\right\}$$ 
is finite for any two $x,y \in N$. 
The distance is obviously left
invariant.

A Carnot group $N$ can be therefore described in the following way.  

\begin{prop}
$(N,d)$ is a path metric cone and also a group such that left translations 
are isometries. Also, dilatations are group isomorphisms. 
\end{prop}

Associate to the CC distance $d$ the function: 
$$\mid x \mid_{d} \ = \ d(0,x)$$ 
This function has the following properties: 
\begin{enumerate}
\item[(a)] the set $\left\{ x \in G \ : \ \mid x \mid_{d} = 1 \right\}$ 
does not contain $0$.
\item[(b)] $\mid x^{-1} \mid_{d} = \mid x \mid_{d}$ for any $x \in G$. 
\item[(c)] $\mid \delta_{\varepsilon} x \mid_{d} = \varepsilon 
\mid x \mid_{d}$ for 
any $x \in G$ and $\varepsilon > 0$. 
\item[(d)] for any $x,y \in N$ we have 
$$\mid xy \mid_{d} \ \leq \ \mid x \mid_{d} + \mid y \mid_{d}$$
\end{enumerate}

\begin{prop}
Balls with respect to the distance $d$ are open and there is a basis for the 
topology on $N$ (as a manifold) made by $d$ balls. 
\end{prop}

\paragraph{Proof.}
We have to prove that $d$ induces the same topology on $N$ as the manifold topology, which is the standard topology induced by Euclidean distance on 
$N$. 

Set  $d_{R}$ to be  the Riemannian distance induced by a 
metric which extends the metric on $D$ and makes the direct sum 
$N \ = \ \oplus_{i} V_{i}$ orthogonal.  We shall have then the inequality 
$$d_{R}(x,y) \ \leq \ d(x,y)$$
for all $x,y \in N$. 

The function $\mid \cdot \mid_{d}$ is also continuous. Indeed, from 
Lemma \ref{p2.4} we see that 
$$\mid x \mid_{d} \ \leq \ \sum_{i = 1}^{M} \mid t_{i} \mid$$
with $t_{i} \rightarrow 0$ when $\|x \| \rightarrow 0$. 

These two facts prove
the thesis.
\quad $\blacksquare$

The map $\mid \cdot \mid_{d}$ looks like a norm. However it is
intrinsically defined and hard to work with. This is the reason for
the introduction of homogeneous norms.

\begin{defi}
A continuous function from $x \mapsto \mid x \mid$ from $G$ to 
$[0,+\infty)$ is a homogeneous norm if
\begin{enumerate}
\item[(a)] the set $\left\{ x \in G \ : \ \mid x \mid = 1 \right\}$ 
does not contain $0$.
\item[(b)] $\mid x^{-1} \mid = \mid x \mid$ for any $x \in G$. 
\item[(c)] $\mid \delta_{\varepsilon} x \mid = \varepsilon 
\mid x \mid$ for 
any $x \in G$ and $\varepsilon > 0$. 
\end{enumerate}
\label{d2.3}
\end{defi}

Homogeneous norms exist. An example is: 
$$\mid x \mid_{1} \ =  \ \sum_{i=1}^{m} \mid x_{i} \mid^{1/i}$$

\begin{prop}
Any two homogeneous norms are equivalent. Let $\mid \cdot \mid$ be a
homogeneous norm. Then the set 
$\left\{ x \mbox{ : } \mid x \mid = 1 \right\}$ is compact. 

There is a constant $C> 0$ such that for any $x,y \in N$ we have: 
$$ \mid xy \mid \ \leq \ C \ ( \mid x \mid + \mid y \mid ) $$
\label{p11}
\end{prop}

\paragraph{Proof.}
We show first that there are constants $c, C > 0$ such that for any 
$x \in  N$ we have 
$$ c \mid x \mid_{1} \ \leq \  \mid x \mid \ \leq \  C \mid x \mid$$
For this consider the compact set $B = \left\{ x \mbox{ : } \mid x 
\mid_{1} = 1 \right\}$. The norm $\mid \cdot \mid_{1}$ has a minimum
$c$ and a maximum $C$ on the set $B$ ($0<c \leq C$). Let now $x \in
N$, $x \not = 0$. Then $x = \delta_{\varepsilon} y$ with
$\varepsilon \ = \ \mid x \mid_{1}$. The homogeneity of the norms
show that $\mid y \mid_{1} = 1$ and the proof is done. 

A consequence of the equivalence of the norms is that for any norm 
$\mid \cdot \mid$ the set $\left\{ x \mbox{ : } \mid x \mid = 1 
\right\}$ is compact. 

For the third assertion remark that the set in $N^{2}$ of all
$(x,y)$ such that $\mid x \mid  + \mid y \mid \ \leq \ 1$ is compact. 
The function $(x,y) \mapsto \mid x+y \mid$ attains a maximum $C$ on
this set. Use again the homogeneity of the norm to finish the proof. 
\quad $\blacksquare$

Let us denote by $\| \cdot \|$ the Euclidean norm on $N$. The
following estimate is easily obtained using homogeneity, as in the
previous proposition. 

\begin{prop}
Let $\mid \cdot \mid$ be a homogeneous norm. Then there are
constants $c, C > 0$ such that for any $x \in N$, $\mid x \mid < 1$,
we have 
$$c \| x \| \ \leq \ \mid x \mid \ \leq \  C \| x \|^{1/m}$$
\label{pest}
\end{prop}

The balls in CC-distance look roughly like boxes (as it is the case 
with the Euclidean balls). A box is a set 
$$Box(r) \ = \ \left\{ x \ = \ \sum_{i=1}^{m} x_{i} \ \mbox{ : } 
\| x_{i} \| \leq r^{i} \right\}$$

\begin{prop} ("Ball-Box theorem") 
There are positive constants $c,C$ such that for any $r>0$ we have: 
$$Box(c r) \ \subset \ B(0,r) \ \subset \  Box( C r)$$
\label{bbprop}
\end{prop}

\paragraph{Proof.}
The box $Box(r)$ is the ball with radius $r$ with respect to the 
homogeneous norm 
$$\mid x \mid_{\infty} \ = \  \max \left\{ \|x_{i} \|^{1/i}
\right\}$$
By Proposition \ref{p11} the norm $\mid \cdot \mid_{\infty}$ is
equivalent with $\mid \cdot \mid_{d}$, which proves the thesis.
\quad $\blacksquare$

A consequence of the Ball-Box Theorem is 

\begin{thm} (Hopf-Rinow theorem for Carnot groups)
Any two points in a Carnot group can be joined by a geodesic. 
\end{thm}

\paragraph{Proof.}
Take two arbitrary points in $N$. Because the distance is left
invariant we can choose $x = 0$. Because dilatations change the
distances by a constant factor, we can suppose that $y \in B(0,1)$. 

A previous proposition shows that the topology generated by a
homogeneous norm and the topology generated by the Euclidean norm 
are the same. Therefore the space is metrically complete. 
The Ball-Box Theorem implies that the ball $B(0,1)$ is compact. 
We are in the assumptions of the abstract Hopf-Rinow theorem \ref{thr}. 
\quad $\blacksquare$

An easy computation shows that 
$$\mathcal{L}(B(0,\delta_{\varepsilon} r)) \ = \ \varepsilon^{Q} 
\mathcal{L} (B(0,r))$$

As a consequence the Lebesgue measure is absolutely continuous with 
respect to the Hausdorff measure $\mathcal{H}^{Q}$. Because of the
invariance with respect to the group operation it follows that the
Lebesgue measure is a multiple of the mentioned Hausdorff measure. 

The following theorem is Theorem 2. from Mitchell \cite{mit}, in the particular case of Carnot groups. 

\begin{thm}
The ball $B(0,1)$ has Hausdorff dimension $Q$. 
\end{thm}

\paragraph{Proof.}
We know that the volume of a ball with radius $\varepsilon$ is 
$c \varepsilon^{Q}$. 

Consider a maximal filling of $B(0,1)$ with balls of radius 
$\varepsilon$. There are $N_{\varepsilon}$ such balls in the
filling; an upper bound for this number is: 
$$ N_{\varepsilon} \ \leq \  1/\varepsilon^{Q}$$
The set of concentric balls of radius $2\varepsilon$ cover 
$B(0,1)$; each of these balls has diameter smaller than $4
\varepsilon$, so the Hausdorff $\alpha$ measure of $B(0,1)$ is 
smaller than
$$\lim_{\varepsilon \rightarrow 0} N_{\varepsilon}
(2\varepsilon)^{\alpha}$$ 
which is $0$ if $\alpha > Q$. Therefore the Hausdorff dimension is 
smaller than $Q$. 

Conversely, given any covering of $B(0,1)$ by sets of diameter $\leq
\varepsilon$, there is an associated covering with balls of the same
diameter; the number $M_{\varepsilon}$ of this balls has a lower 
bound: 
$$M_{\varepsilon} \ \geq \ 1/\varepsilon^{Q}$$
thus there is a lower bound 
$$\sum_{cover} \varepsilon^{\alpha} \ \geq  \ \varepsilon^{\alpha} / 
\varepsilon^{Q}$$
which shows that if $\alpha < Q$ then $\mathcal{H}^{\alpha}(B(0,1))
= \infty$. Therefore the Hausdorff dimension of the ball is greater
than $Q$. 
\quad $\blacksquare$


We collect the important facts discovered until now: 
let $N$ be a Carnot group endowed with the left invariant distribution
generated by $V_{1}$ and with an Euclidean norm on $V_{1}$.

\begin{enumerate}
\item[(a)] If $V_{1}$ Lie-generates the whole Lie algebra of $N$
then any two points can be joined by a horizontal path. 
\item[(b)] The metric  topology of $N$ is the same as Euclidean 
topology. 
\item[(c)] The ball $B(0,r)$ looks roughly like the box 
$\left\{ x \ = \ \sum_{i=1}^{m} x_{i} \ \mbox{ : } 
\| x_{i} \| \leq r^{i} \right\}$.
\item[(d)] the Hausdorff measure $\mathcal{H}^{Q}$ is group
invariant and the Hausdorff dimension of a ball is $Q$.
\item[(e)] the Hopf-Rinow Theorem applies. 
\item[(f)] there is a one-parameter group of dilatations, where a 
dilatation is an isomorphism $\delta_{\varepsilon}$ of $N$ which 
transforms the distance $d$ in $\varepsilon d$. 
\end{enumerate}

\subsection{Pansu differentiability}

A Carnot group has it's own concept of differentiability, introduced by Pansu \cite{pansu}.

In Euclidean spaces, given $f: R^{n} \rightarrow R^{m}$ and 
a fixed point $x \in R^{n}$, one considers the difference function: 
$$X \in B(0,1) \subset R^{n} \  \mapsto \ \frac{f(x+ tX) - f(x)}{t} \in R^{m}$$
The convergence of the difference function as $t \rightarrow 0$ in the uniform 
convergence gives rise to the concept of differentiability in it's classical sense. The same convergence, but in measure, leads to approximate differentiability. 
Other topologies might be considered (see Vodop'yanov \cite{vodopis}).

In the frame of Carnot groups the difference function can be written using only dilatations and the group operation. Indeed, for any function between Carnot groups 
$f: G \rightarrow P$,  for  any fixed point $x \in G$ and $\varepsilon >0$  the finite difference function is defined by the formula: 
$$X \in B(1) \subset G \  \mapsto \ \delta_{\varepsilon}^{-1} \left(f(x)^{-1}f\left( 
x \delta_{\varepsilon}X\right) \right) \in P$$
In the expression of the finite difference function enters $\delta_{\varepsilon}^{-1}$ and $\delta_{\varepsilon}$, which are dilatations in $P$, respectively $G$. 

Pansu's differentiability is obtained from uniform convergence of the difference 
function when $\varepsilon \rightarrow 0$.
 
The derivative of a function $f: G \rightarrow P$ is linear in the sense 
explained further.  For simplicity we shall consider only the case $G=P$. In this way we don't have to use a heavy notation for the dilatations. 

\begin{defi}
Let $N$ be a Carnot group. The function 
$F:N \rightarrow N$ is linear if 
\begin{enumerate}
\item[(a)] $F$ is a {\it group} morphism, 
\item[(b)] for any $\varepsilon > 0$ $F \circ \delta_{\varepsilon} \
= \ \delta_{\varepsilon} \circ F$. 
\end{enumerate}
We shall denote by $HL(N)$ the group of invertible linear maps  of 
$N$, called the  linear group of $N$. 
\label{dlin}
\end{defi}

The condition (b) means that $F$, seen as an algebra morphism, 
preserves the grading of $N$. 

The definition of Pansu differentiability follows: 

\begin{defi}
Let $f: N \rightarrow N$ and $x \in N$. We say that $f$ is 
(Pansu) differentiable in the point $x$ if there is a linear 
function $Df(x): N \rightarrow N$ such that 
$$\sup \left\{ d(F_{\varepsilon}(y), Df(x)y) \ \mbox{ : } \ y \in B(0,1)
\right\}$$
converges to $0$ when $\varepsilon \rightarrow 0$. The functions $F_{\varepsilon}$
are the finite difference functions, defined by 
$$F_{t} (y) \ = \ \delta_{t}^{-1} \left( f(x)^{-1} f(x
\delta_{t}y)\right)$$
\end{defi}

The definition says that $f$ is differentiable at $x$ if the 
sequence of finite differences $F_{t}$ uniformly converges to a
linear map when $t$ tends to $0$. 

We are interested to see how this differential looks like. For 
any $f: N \rightarrow N$, $x,y \in N$ $Df(x)y$ means
$$Df(x)y \ =  \ \lim_{t \rightarrow 0} \delta_{t}^{-1} 
\left( f(x)^{-1} f(x \delta_{t} y) \right)$$
provided that the limit exists.

\begin{prop}
Let $f: N \rightarrow N$, $y, z \in N$ such that: 
\begin{enumerate}
\item[(i)] $Df(x)y$, $Df(x)z$ exist for any $x \in N$. 
\item[(ii)] The map $x \mapsto Df(x)z$ is continuous. 
\item[(iii)] $x \mapsto 
\delta_{t}^{-1} \left( f(x)^{-1} f(x \delta_{t}z)\right)$
converges uniformly to $Df(x)z$. 
\end{enumerate}
Then for any $a,b > 0$ and any $w \ = \ \delta_{a} x \delta_{b} y$
the limit $Df(x)w$ exists and we have: 
$$Df(x) \delta_{a} x \delta_{b} y \ = \ \delta_{a} Df(x)y \ 
\delta_{b} Df(x) z$$
\label{ppansuprep}
\end{prop}

\paragraph{Proof.}
The proof is standard. Remark that if $Df(x)y$ exists then for any 
$a>0$ the limit $Df(x) \delta_{a} y$ exists and 
$$Df(x) \delta_{a} y \ = \ \delta_{a} Df(x) y$$
It is not restrictive therefore to suppose that 
$w = yz$, such that $\mid y \mid_{d} \ =  \ 1$. 

We write: 
$$\delta_{t}^{-1} \left( f(x)^{-1} f(x \delta_{t}w) \right) \ = \ 
(1)(2)(3)$$
where 
$$(1) \ = \ \delta_{t}^{-1} \left( f(x)^{-1} f(x \delta_{t}y)
\right)$$
$$(2) \ = \ \delta_{t}^{-1} \left( f(x \delta_{t}y)^{-1} 
f(x \delta_{t}y \delta_{t} z) \right) \left( Df(x \delta_{t}y) z
\right)^{-1}$$
$$(3) \ = \ Df(x \delta_{t}y) z$$
When $t$ tends to $0$ (i) implies that $(1)$ tends to $Df(x)y$, 
(ii) implies that $(3)$   tends to $Df(x) z$ and (iii) implies that 
$(2)$ goes to $0$. The proof is finished. 
\quad $\blacksquare$

If $f$ is Lipschitz then the previous proposition holds almost
everywhere. In  Proposition 3.2, Pansu \cite{pansu} there is a more general
statement: 

\begin{prop}
Let $f: N \rightarrow N$ have finite dilatation and suppose that for 
almost any $x \in N$ the limits $Df(x)y$ and $Df(x) z$ exist. Then for
almost any $x$ and for any $w \ = \ \delta_{a} y \ \delta_{b} z$ 
the limit $Df(x) w$ exists and we have: 
$$Df(x) \delta_{a} x \delta_{b} y \ = \ \delta_{a} Df(x)y \ 
\delta_{b} Df(x) z$$
\label{ppansu}
\end{prop}

\paragraph{Proof.}
The idea is to use  Proposition \ref{ppansuprep}. We shall suppose 
first that $f$ is Lipschitz. Then any finite difference function
$F_{t}$ is also Lipschitz. 

Recall the general Egorov theorem: 

\begin{lema} (Egorov Theorem)
Let $\mu$ be a finite measure on $X$ and $(f_{n})_{n}$ a sequence of
measurable functions which converges $\mu$ almost everywhere (a.e.) 
to $f$. Then for any $\varepsilon> 0$ there exists a measurable set 
$X_{\varepsilon}$ such that $\mu(X \setminus X_{\varepsilon}) <
\varepsilon$ and $f_{n}$ converges uniformly to $f$ on
$X_{\varepsilon}$. 
\end{lema}

\paragraph{Proof.}
It is not restrictive to suppose that $f_{n}$ converges pointwise to
$f$. Define, for any $q,p \in \mathbf{N}$, the set: 
$$X_{q,p} \ = \ \left\{ x \in X \ \mbox{ : } \  \mid f_{n}(x) - f(x) 
\mid \ < \ 1/p \ , \ \forall n \geq q \right\}$$
For fixed $p$ the sequence of sets $X_{q,p}$ is increasing and at
the limit it fills the space: 
$$\bigcup_{q \in \mathbf{N}} X_{q,p} \ = \ X$$
Therefore for any $\varepsilon > 0$ and any $p$ there is a $q(p)$ 
such that $\mu(X \setminus X_{q(p) p}) \ < \ \varepsilon / 2^{p}$. 

Define then $$X_{\varepsilon} \ = \ \bigcap_{p \in N} X_{q(p) p}$$ 
and check that it satisfies the conclusion. 
\quad $\blacksquare$

We can improve the conclusion of Egorov Theorem by claiming that if 
$\mu$ is a Borel measure then each  $X_{\varepsilon}$ can be chosen 
to be open. 

In our situation we can take $X$ to be a ball in  the group $N$ and 
$\mu$ to be the $Q$ Hausdorff measure, which is Borel. Then we are
able to apply Proposition \ref{ppansuprep} on $X_{\varepsilon}$.
When we tend $\varepsilon$ to $0$ we obtain the claim (for $f$
Lipschitz).  
\quad $\blacksquare$

A consequence of this proposition is   Corollaire 3.3, Pansu \cite{pansu}.  

\begin{cor}
If $f: N \rightarrow N$ has finite dilatation and $X_{1}, ... , X_{K}$ is 
a basis for the distribution $V_{1}$ such that 
$$s \mapsto f(x \delta_{s} X_{i})$$ is differentiable in $s=0$ for
almost any $x \in N$, then $f$ is differentiable almost everywhere 
and the differential is linear, i.e.  if 
$$y \ = \ \prod_{i = 1}^{M} \delta_{t_{i}} X_{g(i)}$$
then 
$$Df(x) y \ = \ \prod_{i = 1}^{M} \delta_{t_{i}} X_{g(i)}(f)(x)$$
\label{corr1}
\end{cor}

With the definition of (Pansu) derivative at hand, it is natural to introduce the class of $C^{1}$ maps from $N$ to $M$. 

\begin{defi}
Let $N,M$ be  Carnot groups. A function $f: N \rightarrow M$ is of 
class $C^{1}$ if  it is  Pansu derivable everywhere and the derivative is continuous as a function $Df : N \times N \rightarrow M \times M$. Here 
$N \times N$ and $M \times M$ are endowed with the product topology. 
\end{defi}

\begin{rk}
For example left translations $L_{x}: N \rightarrow N$, $L_{x}(y) \ = \ xy$ 
are $C^{1}$ but right translations $R_{x}: 
N \rightarrow N$, $L_{x}(y) \ = \ yx$ are not even derivable. 
\end{rk}

We can introduce the class of $C^{1}$ $N$ manifolds. 

\begin{defi}
Let $N$ be a Carnot group. A $C^{1}$ $N$ manifold is an 
(equivalence class of, or a maximal) atlas over the group $N$ such that the 
change of charts is locally bi-Lipschitz. 
\label{dcum}
\end{defi}

\subsection{Rademacher theorem}

A very important result is 
theorem 2, Pansu \cite{pansu}, which  contains the 
Rademacher theorem for Carnot groups. 

\begin{thm}
Let $f: M \rightarrow N$ be a Lipschitz  function between Carnot groups. 
Then $f$ is differentiable almost everywhere.
\label{ppansuu}
\end{thm}

The proof is based on the corollary \ref{corr1} and the technique of development of a curve, Pansu sections 4.3 - 4.6 \cite{pansu}.

Let $c: [0,1] \rightarrow N$ be a Lipschitz curve such that $c(0) = 0$.  To 
any division $\Sigma: 0=t_{0} < .... < t_{n} = 1$ of the interval $[0,1]$ is a associated the element $\sigma_{\Sigma} \in N$ (the algebra) given by: 
$$\sigma_{\Sigma} \ = \ \sum_{k=0}^{n} c(t_{k})^{-1}c(t_{k+1})$$
Lemma (18) Pansu \cite{pansu1}  implies the existence of a constant $C>0$ such that 
\begin{equation}
\| \sigma_{\Sigma} \ - \ c(1) \| \ \leq \ C \left(  \sum_{k=0}^{n} 
\|c(t_{k})^{-1}c(t_{k+1})\| \right)^{2}
\label{pu1}
\end{equation}

Take now a finer division $\Sigma'$ and look at the interval 
$[t_{k}, t_{k+1}]$ divided further by $\Sigma'$ like this: 
$$t_{k} = t'_{l} < ... < t'_{m}=t_{k+1}$$
The estimate \eqref{pu1} applied for each interval $[t_{k}, t_{k+1}]$ lead 
us to the inequality: 
$$\| \sigma_{\Sigma'} \ - \ \sigma_{\Sigma} \| \ \leq \ \sum_{k=0}^{n} d(c(t_{k}), c(t_{k+1}))^{2}$$
The curve has finite length (being Lipschitz) therefore the right hand side 
of the previous inequality tends to 0 with the norm of the division. Set 
$$\sigma(s) \ = \ \lim_{\| \Sigma \| \rightarrow 0} \sigma_{\Sigma}(s)$$
where $\sigma_{\Sigma}(s)$ is relative to the curve $c$ restricted to the interval $[0,s]$. The curve such defined is called the development of 
the curve $c$. It is easy to see that $\sigma$ has the same length 
(measured with the Euclidean norm $\| \cdot \|$) as $c$. If we parametrise $c$ 
with the length then we have the estimate: 
\begin{equation}
\|\sigma(s) \- \ c(s) \| \ \leq \ C s^{2}
\label{impes}
\end{equation}
This shows that $\sigma$ is a Lipschitz curve (with respect to the Euclidean distance). Indeed, prove first that $\sigma$ has finite dilatation in almost 
any point, using \eqref{impes} and the fact that $c$ is Lipschitz. Then show that the dilatation is majorised by the Lipschitz constant of $c$. 
By the classical Rademacher theorem $\sigma$ is almost everywhere derivable. 

\begin{rk}
The development of a curve can be done in an arbitrary Lie connected group, 
endowed with a left invariant distribution which generates the algebra. 
One should add some logarithms, because in the case of Carnot groups, we have 
identified the group with the algebra. The inequality \eqref{impes} still holds. 
\label{impr}
\end{rk}

Conversely, given a curve $\sigma$ in the algebra $N$, we can perform the inverse operation to development (called multiplicative integral by Pansu; 
we shall call it "lift"). Indeed, to any division $\Sigma$ of the interval 
$[0,s]$ we associate 
the point 
$$c_{\Sigma}(s) \ = \ \prod_{k=0}^{n} \left( \sigma(t_{k+1}) \ - \ \sigma(t_{k}) \right)$$
Define then 
$$c(s) \ = \ \lim_{\| \Sigma \| \rightarrow 0} c_{\Sigma}(s)$$
Remark that if $\sigma([0,1]) \subset V_{1}$ and it is almost everywhere differentiable then $c$ is a horizontal curve. 

Fix $s \in [0,1)$ and parametrise $c$ by the length. Apply the inequality 
\eqref{impes} to the curve $t \mapsto c(s)^{-1}c(s+t)$. We get the fact that 
the vertical part of $\sigma(s+t) - \sigma(s)$ is controlled by 
$s^{2}$. Therefore, if $c$ is Lipschitz then $\sigma$ is included in $V_{1}$. 

Denote the $i$ multiple bracket $[x,[x,...[x,y]...]$ by $[x,y]_{i}$. For any 
division $\Sigma$ of the interval $[0,s]$ set:
$$A^{i}_{\Sigma}(s) \ = \ \sum_{k=0}^{n} [\sigma(t_{k}),\sigma(t_{k+1})]_{i}$$
As before, one can show that when the norm of the division tends to 0, 
$A^{i}_{\Sigma}$ converges. Denote by 
$$A^{i}(s) \ = \ \int_{0}^{s} [\sigma, d\sigma]_{i} \ = \  \lim_{\| \Sigma \| \rightarrow 0} A^{i}_{\Sigma}(s)$$
the i-area function. The estimate corresponding to \eqref{impes} 
is 
\begin{equation}
\| A^{i}(s) \| \ \leq \ C s^{i+1}
\label{impesai}
\end{equation}

What is the signification of $A^{i}(s)$? The answer is simple, based on the well known formula of derivation of left translations 
in a group. Define, in a neighbourhood of $0$ in the Lie algebra 
$\mathfrak{g}$ of the Lie group $G$, the operation 
$$X \opg Y \ = \ \log_{G}\left( \exp_{G}(X) \exp_{G}(Y) \right)$$
The left translation by $X$ is the function $L_{X}(Y) \ = \ X \opg Y$. It is known that 
$$D L_{X}(0) (Z) \ = \ \sum_{i=0} \frac{1}{(i+1)!} [X,Z]_{i}$$
It follows that if the group $G$ is Carnot then $[X,Z]_{i}$ measures the infinitesimal variation of the $V_{i}$ component of $L_{X}(Y)$, for 
$Y = 0$. Otherwise said, 
$$\frac{d}{dt}\sigma_{i}(t) \ = \  \frac{1}{(i+1)!} \lim_{\varepsilon \rightarrow 0} \varepsilon^{-i} \ \int_{t}^{t + \varepsilon} [\sigma, d\sigma]_{i}$$
Because $\sigma$ is Lipschitz, the left hand side exists for all $i$ for a.e. 
$t$. The estimate \eqref{impesai} tells us that the right hand side equals $0$. This implies the following proposition (Pansu, 4.1 \cite{pansu}). 

\begin{prop}
If $c$ is Lipschitz  then $c$ is differentiable almost everywhere. 
\label{pr41}
\end{prop}

\paragraph{Proof of theorem \ref{ppansuu}.}
The proposition \ref{pr41} implies that we are in the hypothesis of 
corollary \ref{corr1}. 
\quad $\blacksquare$

Pansu-Rademacher theorem has been improved for Lipschitz functions 
$f: A \subset M \rightarrow N$, where M, N are Carnot groups and $A$ 
is just measurable, in Vodop'yanov, Ukhlov \cite{voduk}. Their technique differs from Pansu. It resembles  with the one used in Margulis, Mostow \cite{marmos1}, where 
it is proven that any quasi-conformal map from a Carnot-Carath\'eodory manifold 
to another is a.e. differentiable, without using Rademacher theorem. Same result for quasi-conformal maps on Carnot groups has been first proved by Koranyi, Reimann \cite{kore}. 
Finally  Magnani \cite{magnani} reproved a.e. differentiability of Lipschitz functions on Carnot groups, defined on measurable sets,  continuing Pansu technique.  

For a review of other connected results, such as approximate differentiability, 
differentiability in Sobolev topology, see the excellent Vodop'yanov \cite{vodopis}.

\subsection{Area formulas}

The subject of this subsection is the change of variable formula. In order to prove it one needs more involved (though basic) knowledge on metric measure 
spaces. We shall just sketch the proofs, delaying for the second part of these 
notes the  rigorous proofs. 

\begin{defi}
A metric measure (or mm) space is a metric space $(X,d)$ endowed with a Borel measure 
$\mu$ satisfying the doubling condition, i.e. there is a positive constant $C$ such that 
$$\mu(B(x,2r)) \ \leq \ C \mu(B(x,r))$$
for any $x \in X$, $r>0$. 
\end{defi}

A Carnot group is a metric measure space, endowed with the invariant volume measure (which is the Lebesgue measure). Indeed, this is a consequence of the existence of dilatations. 

In any mm space the Vitali covering theorem holds: 

\begin{thm}
(Vitali) Let $(X,d,\mu)$ be a mm space, $A \subset X$ and $\mathcal{F}$ a family of closed sets in $X$ which covers $A$. If for any $x \in A$ and 
$\varepsilon > 0$ there exists $V \in \mathcal{F}$ such that $x \in V$ and 
$diam \ V \ < \varepsilon$ then one can extract from $\mathcal{F}$ a countable 
disjoint family $(V_{i})_{i \in N}$ such that 
$$\mu( A \setminus \cup_{i \in N} V_{i}) \ = \ 0$$
\end{thm}

Vitali covering theorem is the only technical ingredient needed to prove the following result (compare with Margulis, Mostow \cite{marmos1}, lemma 2.3). But before this we need  a definition.

\begin{defi}
Let $(X,d,\mu)$ be a mm space and $\sigma$ a measure on $X$. The lower spherical density of the measure $\sigma$ at the point 
$x \in X$ is 
$$\frac{d \sigma}{d \mu}^{-} (x) \ = \ \liminf_{\varepsilon \rightarrow 0} 
\frac{\sigma(B(x,\varepsilon))}{\mu(B(x,\varepsilon))}$$

If $(Y,\sigma)$ is  measure space and $f: X \rightarrow Y$ is measurable 
then the pull-back of $\sigma$ by $f$ is the measure $f^{*}(\sigma)$ on 
$X$ defined by: 
$$f^{*}(\sigma) (A) \ = \ \sigma(f(A))$$

The Jacobian of $f$ in $x \in X$ is the spherical density of $f^{*}(\sigma)$ in 
$x$: 
$$J f(x) \ = \  \frac{d f^{*}(\sigma)}{d \mu}^{-} (x) \ = \ \liminf_{\varepsilon \rightarrow 0} 
\frac{\sigma(f(B(x,\varepsilon)))}{\mu(B(x,\varepsilon))}$$
\end{defi}

\begin{prop}
Let $(X,d,\mu)$ be a mm space, $(Y,\sigma)$ a finite measure space and $f: X \rightarrow Y$ be an injective measurable map. We have then the inequality: 
\begin{equation}
\int_{X} J f(x) \mbox{ d}\mu \ \leq \ \sigma(f(X))
\label{pe23}
\end{equation}
\end{prop}

In order to establish the change of variables (or area formula) for 
$f$ bi-Lipschitz between Carnot groups one has to prove first that in 
\eqref{pe23} we have equality. This is done by showing that 
$\mu(Z) = 0$, where $Z$ is the set of zeroes of $J f$. Further, one has to use the Rademacher theorem \ref{ppansuu} in order to show that the Jacobian equals the (modulus of the) determinant of the derivative, almost everywhere. Indeed, the idea is to 
use the uniform convergence in the definition of Pansu derivative. We arrive to the change of variables theorem. 

\begin{thm}
Let $f: A \subset M \rightarrow N$ be a bi-Lipschitz function from a measurable 
subset $A$ of a Carnot group $M$ to another Carnot group $N$. We have then: 
\begin{equation}
\int_{A} \mid \det Df(x) \mid \mbox{ d}x \ = \ vol \ f(A)
\end{equation}
\label{cvf}
\end{thm}

A corollary of the area formula is the fact that bi-Lipschitz maps have Lusin property. 

\begin{cor}
 Let $f: A \subset M \rightarrow N$ be a bi-Lipschitz function from a measurable 
subset $A$ of a Carnot group $M$ to another Carnot group $N$. If $B \subset A$ 
is a set of $0$ measure then $f(B)$ has $0$ measure.
\label{clus} 
\end{cor}

\subsection{Rigidity phenomenon}
\label{rigid}

As an application to the previous sections we shall prove the following well known rigidity result.

\begin{thm}
Let X,Y be Lipschitz manifolds over the Carnot groups M, respectively N. 
If there is no subgroup of $N$ isomorphic with $M$ then there is no bi-Lipschitz embedding of $X$ in $Y$. 
\label{tr1}
\end{thm}

\paragraph{Proof.}
Suppose that such an embedding exists. This implies (by choosing charts for 
$X$, $Y$) that there is a bi-Lipschitz  embedding $f: A \subset M \rightarrow N$ with $A$ open. Rademacher theorem \ref{ppansuu} and corollary \ref{clus} imply that for almost any $x \in A$ $Df(x)$ exists and it is an injective morphism from $M$ to $N$. This contradicts the hypothesis. 
\quad $\blacksquare$

For example there is no bi-Lipschitz function from an open set in a Heisenberg group to an Euclidean space $R^{n}$. Otherwise stated, the Heisenberg group 
(or any other non-commutative Carnot group), is purely unrectifiable. This 
calls for an intrinsic definition of rectifiability. The subject will be treated in the second part of these notes, where we shall propose the following point of view: rectifiability is a theory of irreducible representations of 
the group of bi-Lipschitz homeomorphisms. 

Let us give another example of rigidity, connected to Sussmann \cite{sus},  
section 8,  
example of abundance of abnormal geodesics on four-dimensional sub-Riemannian 
manifolds. 

Consider a connected Lie group $G$ with Lie algebra $\mathfrak{g}$ generated by $X_{1}$, $X_{2}$, $X_{3}$, $X_{4}$, such that 
the following conditions hold: 
\begin{enumerate}
\item[(i)] $X_{3} \ = \ [X_{1},X_{2}]$, $X_{4} \ = \ [X_{1}, X_{3}]$, 
\item[(ii)] $[X_{2},X_{3}] \ = \ \alpha X_{1} \ + \ \beta X_{2} \ + \ \gamma 
X_{3}$ such that $\beta \not = 0$
\end{enumerate}
Close examination shows that there is a three-dimensional family of Lie algebras satisfying these bracket conditions. Set $D \ = \ span \ \left\{ X_{1}, 
X_{2} \right\}$ and endow $G$ with the left invariant distribution generated by $D$ and with the metric obtained by left translation of a metric on $D$ such that $X_{1}, X_{2}$ are orthonormal. Sussmann proves that 
$t \mapsto \exp_{G}(tX_{2})$ is an abnormal geodesic (see the second part 
of these notes or Sussmann \cite{sus}, or Montgomery \cite{montgo} and the references therein). 

In this section we shall look to the nilpotentisation of the group $G$ 
(section \ref{nilpo}). The algebra $\mathfrak{g}$ admits the filtration 
$$V^{1} \ = \ span \ \left\{X_{1},X_{2}\right\} \ \subset \ V^{2} \ = \ 
span \ \left\{X_{1},X_{2}, X_{3}\right\} \ \subset \ \mathfrak{g}$$
The nilpotentisation is $\mathfrak{g}$ endowed with the following bracket: 
$$[X_{1},X_{2}]_{N} \ = \ X_{3} \ , \ \ [X_{1}, X_{3}]_{N} \ = \ X_{4}$$
all the other brackets being zero. The dilatations in this Carnot algebra 
are 
$$\delta_{\varepsilon} (\sum a_{i}X_{i}) \ = \ \varepsilon a_{1}X_{1} \ + \ 
\varepsilon a_{2}X_{2} \ + \ \varepsilon^{2} a_{3}X_{3} \ + \ \varepsilon^{3} a_{4}X_{4}$$
Denote by $N$ the Carnot group which has the algebra $(\mathfrak{g}, [ \cdot , \cdot]_{N})$. The group operation on $N$ will be denoted by $\opn$. 

The linear group $HL(N)$ is the class 
of linear invertible maps which commute with the dilatations and preserve the 
nilpotent bracket. We give without proof the description of this group 
(for a similar proof see proposition \ref{p1}). 

\begin{prop}
In the basis $\left\{ X_{1}, ... , X_{4}\right\}$ 
the linear group $HL(N)$ is made by all linear transformations represented by the matrices of the form 
$$\left( \begin{array}{cccc}
a_{11} & 0 & 0 & 0 \\ 
a_{21} & a_{22} & 0 & 0 \\ 
0 & 0 & a_{11}a_{22} & 0 \\ 
0 & 0 & 0 & a_{11}^{2}a_{22}
\end{array} \right)$$
for all $a_{11}, a_{22} \not = 0$. 
\label{pni}
\end{prop}

Set $H \ = \ span \ \left\{ X_{1}, X_{3}, X_{4} \right\}$. This is a group 
isomorphic with the Heisenberg group $H(1)$, with respect to the nilpotent bracket. Define now the family of horizontal curves 
$$c_{z}(t) \ = \ z \opn (tX_{2})$$
for any $z \in H$. Any such curve is transversal to $H$. 

\begin{prop}
Let $\phi: \mathfrak{g} \rightarrow \mathfrak{g}$ be a bi-Lipschitz map 
on the Carnot group $N$. Then for almost every $z \in H$ (with respect to the Lebesgue measure on $H$) there is $z' \in H$ such that 
$$\phi(c_{z}(R)) \ \subset  \ c_{z'}(R)$$
Otherwise stated the family of curves $c_{z}$ is preserved by $N$ bi-Lipschitz maps. 
\end{prop}

\paragraph{Proof.}
Let $M$ be the family of all $z \in H$ such that $\phi$ is not a.e. differentiable in the points $c_{z}(t)$. Then $M$ has to be negligible in 
$H$, otherwise we contradict Rademacher theorem \ref{ppansuu}. 
Take now $z \in H \setminus M$. Then for almost any $t \in R$ $\phi$ is differentiable in $c_{z}(t)$ and $\phi \circ c_{z}$ is differentiable in $t$. 
From the chain rule it follows that 
$$\frac{d}{dt}\left( \phi \circ c_{z} \right) (t) \ = \ D \phi(c_{z}(t)) 
\frac{d}{dt} c_{z}(t)$$
(all derivatives are Pansu). But 
$$\frac{d}{dt} c_{z}(t) \ = X_{2}$$
and from proposition \ref{pni} it follows that 
$$\frac{d}{dt}\left( \phi \circ c_{z} \right) (t) \ = \ \lambda(z,t) X_{2}$$
Indeed, $D\phi(x) \ \in HL(N)$ and the direction $X_{2}$ is preserved by any 
linear $F \in HL(N)$. 
Because the curve $\phi \circ c_{z}$ is almost everywhere tangent to 
$X_{2}$ it follows that it is a curve in the family $c_{z}$
\quad $\blacksquare$

Let us turn back to theorem \ref{tr1} and remark that it can be improved. 
With the notation used before $H$ is a subgroup of $N$ and $H$ is isomorphic 
with the Heisenberg group $H(1)$. Nevertheless the isomorphism does not 
commute with dilatations, or better, there is no isomorphism $f: H(1) \rightarrow N$ commuting with dilatations: 
$$f \circ \delta^{H(1)}_{\varepsilon} \ = \ \delta^{N}_{\varepsilon}$$
Therefore there can be no bi-Lipschitz embedding of $H(1)$ in $N$. This is because linear maps are not only group morphisms, but also commute with dilatations. 

\begin{defi}
Let   $N$, $M$ be Carnot groups. We write $N \leq M$ if there is an injective 
group morphism $f: N \rightarrow M$ which commutes with dilatations. 
\end{defi}

Theorem \ref{tr1} improves like this: 

\begin{thm}
Let X,Y be Lipschitz manifolds over the Carnot groups M, respectively N. 
If  $N \not \leq M$ then there is no bi-Lipschitz embedding of $X$ in $Y$. 
\label{tr2}
\end{thm} 

Rigidity in the sense of this section manifests in subtler ways. The purpose of 
Pansu paper \cite{pansu} was to extend a result of Mostow \cite{mosto}, called 
Mostow rigidity. Although it is straightforward now to explain what Mostow rigidity means and how it can be proven, it is beyond the purposes of these notes.

\section{The Heisenberg group} 
\label{shei}

Let us rest a little bit and look closer to  an example. The Heisenberg group 
is the most simple non commutative Carnot group. We shall apply the achievements of the previous chapter to this group. 

\subsection{The group}
The Heisenberg group $H(n) = R^{2n+1}$ is a 2-step  nilpotent group with the 
operation: 
$$(x,\bar{x})  (y,\bar{y}) = (x + y, \bar{x} + \bar{y} + \frac{1}{2} \omega(x,y))$$
where $\omega$ is the standard symplectic form on $R^{2n}$. We shall identify 
the Lie algebra with the Lie group. The bracket is 
$$[(x,\bar{x}),(y,\bar{y})] = (0, \omega(x,y))$$
The Heisenberg algebra is generated by 
$$V = R^{2n} \times \left\{ 0 \right\}$$ 
and we have the relations $V + [V,V] = H(n)$, $\left\{0\right\} \times R \ = \ [V,V] \ = \ Z(H(n))$.

The dilatations on $H(n)$ are 
$$\delta_{\varepsilon} (x,\bar{x}) = (\varepsilon x , \varepsilon^{2} \bar{x})$$

We shall denote by $HL(H(n))$ 
the group of invertible linear transformations and by $SL(H(n))$ the 
subgroup of volume preserving ones. 

\begin{prop}
We have the isomorphisms $$HL(H(n)) \approx CSp(n) \ , \ \ SL(H(n)) \approx  Sp(n)$$
\label{p1}
\end{prop}

\paragraph{Proof.}
By direct computation. We are looking first for the algebra isomorphisms of 
$H(n)$. Let the matrix 
$$\left( \begin{array}{cc} 
          A & b \\
          c & a 
         \end{array} \right)$$
represent such a morphism, with $A \in gl(2n, R)$, $b,c \in R^{2n}$ and 
$a \in R$. The bracket preserving condition reads: for any $(x,\bar{x}), 
(y,\bar{y}) \in H(n)$ we have
$$(0,\omega(A x + \bar{x} b, A y + \bar{y} b)) = ( \omega(x,y) b, a \omega(x,y))$$
We find therefore $b = 0$ and $\omega(Ax, Ay) = a \omega(x,y)$, so $A \in 
CSp(n)$ and $a \geq 0$, $a^{n} = \det A$. 

The preservation of the grading gives $c=0$. The volume preserving condition 
means $a^{n+1} = 1$ hence $a= 1$ and $A \in  Sp(n)$. 
\quad $\blacksquare$

\paragraph{Get acquainted with Pansu differential}
{\bf Derivative of a curve:} 
Let us see which are the smooth (i.e. derivable) curves. Consider 
$\tilde{c}: [0,1] \rightarrow H(n)$, $t \in (0,1)$ and 
$\varepsilon > 0$ sufficiently small. Then the finite difference 
function associated to $c,t, \varepsilon$ is 
$$C_{\varepsilon}(t)(z) \ = \
\delta_{\varepsilon}^{-1} \left( \varepsilon\tilde{c}(t)^{-1}
\tilde{c}(t + \varepsilon z) \right)$$ 
After a short computation we obtain: 
$$C_{\varepsilon}(t)(z) \ = \ \left(\frac{c(t + \varepsilon z) -
c(t)}{\varepsilon} , \frac{\bar{c}(t+\varepsilon z) -
\bar{c}(t)}{\varepsilon^{2}} - \frac{1}{2} \omega(c(t), 
\frac{c(t + \varepsilon z) -
c(t)}{\varepsilon^{2}}) \right)$$
When $\varepsilon \rightarrow 0$ we see that the finite difference
function converges if: 
$$\dot{\bar{c}}(t) \ = \ \frac{1}{2} \omega(c(t), \dot{c}(t))$$
Hence the curve has to be horizontal; in this case we see that 
$$D c(t) z \ = \ z ( \dot{c}(t), 0)$$
This is almost the tangent to the curve. The tangent is obtained 
by taking $z = 1$ and the left translation of $Dc(t) 1$ by $c(t)$.

The horizontality condition implies that,  given a curve 
$t \mapsto c(t) \in R^{2n}$, there is
only one horizontal curve $t \mapsto (c(t),\bar{c}(t))$, such that 
$\bar{c}(0) = 0$. This curve is called the lift of $c$. 

{\bf Derivative of a functional:} 
Take now $f: H(n) \rightarrow R$ and compute its Pansu derivative. 
The finite difference function is 
$$F_{\varepsilon}(x,  \bar{x})(y,\bar{y}) \ = \ 
\left( f(x + \varepsilon y, \bar{x} + \frac{\varepsilon}{2} \omega(x,y) +
\varepsilon^{2} \bar{y}) - f
(x,\bar{x}) \right)/\varepsilon$$
Suppose that $f$  is classically derivable. Then it is Pansu derivable  
and this derivative has the expression: 
$$Df(x , \bar{x}) (y,\bar{y}) \ = \ \frac{\partial f}{\partial
x}(x,\bar{x}) y
+ \frac{1}{2} \omega(x,y) \ \frac{\partial f}{\partial
\bar{x}}(x,\bar{x}) $$
Not any Pansu derivable functional is derivable in the classical sense. As an
example, check that the square of any homogeneous norm is Pansu derivable everywhere, 
but not derivable everywhere in the classical sense.

\subsection{Lifts of symplectic diffeomorphisms}

In this section we are interested in the 
group of volume preserving diffeomorphisms of $H(n)$, with certain classical regularity. We establish connections between volume preserving 
diffeomorphisms of $H(n)$ and symplectomorphisms of $R^{2n}$.

\paragraph{Volume preserving diffeomorphisms}

\begin{defi}
$Diff^{2}(H(n),vol)$ is  the group of volume preserving 
diffeomorphisms $\tilde{\phi}$ of $H(n)$ such that 
$\tilde{\phi}$ and it's inverse have (classical) regularity $C^{2}$. 
In the same way we define $Sympl^{2}(R^{2n})$ to be the group of 
$C^{2}$ symplectomorphisms of $R^{2n}$. 
\end{defi}

\begin{thm}
We have the isomorphism of groups
$$Diff^{2}(H(n),vol) \ \approx \ Sympl^{2}(R^{2n}) \times R$$
given by the mapping 
$$\tilde{f} = (f,\bar{f}) \  \in \ Diff^{2}(H(n),vol) \ \mapsto \ \left( 
f \in Sympl^{2}(R^{2n}) , \bar{f}(0,0) \right)$$
The inverse of this isomorphism has the expression
$$\left( f \in Sympl^{2}(R^{2n}) , a \in R \right) \ \mapsto  \  \tilde{f} = (f,\bar{f}) \  \in \ Diff^{2}(H(n),vol)$$ 
$$\tilde{f}(x,\bar{x}) \ = \ (f(x), \ \bar{x} + F(x))$$
where $F(0)= a$ and $dF \ = \ f^{*} \lambda \ - \ \lambda$. 
\label{t1}
\end{thm}

\paragraph{Proof.}
Let $\tilde{f} = (f,\bar{f}) : H(n) \rightarrow H(n)$ be an element of 
the group \\ 
$Diff^{2}(H(n), vol)$. 
We shall compute: 
$$D \tilde{f} ((x,\bar{x})) (y,\bar{y}) \ = \ 
\lim_{\varepsilon \rightarrow 0} \delta_{\varepsilon^{-1}}  \left( 
\left(\tilde{f}(x,\bar{x})\right)^{-1}  \tilde{f}\left((x,\bar{x})  \delta_{\varepsilon} (y,\bar{y})\right)\right) $$
We know that $D \tilde{f}(x,\bar{x})$ has to be a linear mapping. 

After a short computation we see that we have to pass to the limit 
$\varepsilon \rightarrow 0$ in the following expressions (representing the two 
components of $D \tilde{f} ((x,\bar{x})) (y,\bar{y})$): 
\begin{equation}
\frac{1}{\varepsilon} \left( f\left(x+ \varepsilon y, \bar{x} + \varepsilon^{2} 
\bar{y} + \frac{\varepsilon}{2} \omega(x,y)\right) - f(x,\bar{x}) \right)
\label{exp1}
\end{equation}
\begin{equation}
\frac{1}{\varepsilon^{2}} \left( \bar{f}\left(x+ \varepsilon y, \bar{x} + \varepsilon^{2} 
\bar{y} + \frac{\varepsilon}{2} \omega(x,y)\right) - \bar{f}(x,\bar{x}) - 
\right.  
\label{exp2}
\end{equation}
$$\left. - 
\frac{1}{2}\omega\left(f(x,\bar{x}), f\left(x+ \varepsilon y, \bar{x} + \varepsilon^{2} 
\bar{y} + \frac{\varepsilon}{2} \omega(x,y)\right)\right)\right)
$$

The first component \eqref{exp1} tends to 
$$\frac{\partial f}{\partial x} (x,\bar{x}) y + \frac{1}{2} \frac{\partial f}{\partial \bar{x}} (x,\bar{x}) \omega(x,y)$$
The terms of order $\varepsilon$ must cancel in the second component \eqref{exp2}. We obtain the cancellation 
condition (we shall omit from now on the argument $(x,\bar{x})$ from all functions): 
\begin{equation}
\frac{1}{2} \omega(x,y) \frac{\partial \bar{f}}{\partial \bar{x}} - 
\frac{1}{2} \omega(f, \frac{\partial f}{\partial x} y) - 
\frac{1}{4} \omega(x,y) \omega(f, \frac{\partial f}{\partial \bar{x}}) + 
\frac{\partial \bar{f}}{\partial x} \cdot y \ = \ 0
\label{cancel}
\end{equation}
The second component tends to 
$$\frac{\partial \bar{f}}{\partial \bar{x}} \bar{y} - \frac{1}{2} \omega(f, 
\frac{\partial f}{\partial \bar{x}}) \bar{y}$$
The group morphism $D \tilde{f}(x,\bar{x})$ 
is represented by the matrix: 
\begin{equation}
d \tilde{f}(x,\bar{x}) \ = \ \left( \begin{array}{cc}
\frac{\partial f}{\partial x} + \frac{1}{2} \frac{\partial f}{\partial \bar{x}} 
\otimes Jx & 0 \\
0 & \frac{\partial \bar{f}}{\partial \bar{x}} - \frac{1}{2} \omega(f, 
\frac{\partial f}{\partial \bar{x}}) 
       \end{array} \right)
\label{tang}
\end{equation}
We shall remember now that $\tilde{f}$ is volume preserving. According to 
proposition \ref{p1}, this means: 
\begin{equation}
\frac{\partial f}{\partial x} + \frac{1}{2} \frac{\partial f}{\partial \bar{x}} 
\otimes Jx \ \in Sp(n) 
\label{c1}
\end{equation}
\begin{equation}
\frac{\partial \bar{f}}{\partial \bar{x}} - \frac{1}{2} \omega(f, 
\frac{\partial f}{\partial \bar{x}}) = 1 
\label{c2}
\end{equation}
The cancellation condition \eqref{cancel} and relation \eqref{c2} give
\begin{equation}
\frac{\partial \bar{f}}{\partial x} y \ = \  \frac{1}{2} \omega(f,\frac{\partial 
f}{\partial x} y ) \ - \ \frac{1}{2} \omega(x,y)
\label{c3}
\end{equation}

These conditions describe completely the class of volume preserving diffeomorphisms 
of $H(n)$. Conditions \eqref{c2} and \eqref{c3} are in fact differential equations 
for the function $\bar{f}$ when $f$ is given. However, there is a compatibility 
condition in terms of $f$ which has to be fulfilled for  \eqref{c3} to have 
a solution $\bar{f}$. Let us look closer to \eqref{c3}. We can see the symplectic 
form $\omega$ as a closed 2-form. Let $\lambda$ be a 1-form such that 
$d \lambda = \omega$. If we take the (regular) differential with respect 
to $x$ in \eqref{c3} we quickly obtain the compatibility condition
\begin{equation}
\frac{\partial f}{\partial x} \ \in \ Sp(n)
\label{c4}
\end{equation}
and \eqref{c3} takes the form: 
\begin{equation}
 \ d \bar{f} \ = \ f^{*} \lambda \ - \ \lambda
\label{c5}
\end{equation}
(all functions seen as functions of $x$ only).

Conditions \eqref{c4} and \eqref{c1} imply: there is a scalar function 
$\mu = \mu(x,\bar{x})$ such that 
$$\frac{\partial f}{\partial \bar{x}} \ = \ \mu \ Jx $$
Let us see what we have until now: 
\begin{equation}
\frac{\partial f}{\partial x} \ \in \ Sp(n) 
\label{cc1}
\end{equation}
\begin{equation}
\frac{\partial \bar{f}}{\partial x} \ = \ \frac{1}{2} \left[ \left( 
\frac{\partial f}{\partial x}\right)^{T} J f \ - \ J x \right]
\label{cc2}
\end{equation}
\begin{equation}
\frac{\partial \bar{f}}{\partial \bar{x}} \ = \ 1 + \frac{1}{2} \omega(f, 
\frac{\partial f}{\partial \bar{x}}) 
\label{cc3}
\end{equation}
\begin{equation}
\frac{\partial f}{\partial \bar{x}} \ = \ \mu \ Jx
\label{cc4}
\end{equation}
Now, differentiate \eqref{cc2} with respect to $\bar{x}$ and use \eqref{cc4}. In the same time 
differentiate \eqref{cc3} with respect to $x$. From the equality 
$$\frac{\partial^{2} \bar{f}}{\partial x \partial \bar{x}} \ = \ 
\frac{\partial^{2} \bar{f}}{\partial \bar{x} \partial x}$$ 
we shall obtain by straightforward computation $\mu = 0$. 
\quad $\blacksquare$

\subsection{Hamilton's equations}

For any element $\phi \in Sympl^{2}(R^{2n})$ we define the lift $\tilde{\phi}$ to be  the image 
of $(\phi, 0)$ by the isomorphism described in the theorem \ref{t1}.

Let $A \subset R^{2n}$ be a set. $Sympl^{2}(A)_{c}$ is the group of symplectomorphisms with regularity $C^{2}$, 
with compact support in $A$.

\begin{defi}
For any flow $t \mapsto \phi_{t} \in Sympl^{2}(A)_{c}$ 
denote by $\phi^{h}(\cdot, x))$ the  horizontal flow in 
$H(n)$ obtained by the lift of all curves $t \mapsto \phi(t,x)$ and 
by $\tilde{\phi}(\cdot, t)$ the flow obtained by the lift of all 
$\phi_{t}$. The vertical flow is defined by the expression 
\begin{equation}
\phi^{v} \ = \ \tilde{\phi}^{-1} \circ \phi^{h}
\label{hameq}
\end{equation}
\label{dhameq}
\end{defi}

Relation \eqref{hameq} can be seen as Hamilton equation. 

\begin{prop}
Let $t \in [0,1] \mapsto \phi^{v}_{t}$ be a curve of diffeomorphisms 
of $H(n)$ satisfying the equation: 
\begin{equation}
\frac{d}{dt} \phi^{v}_{t}(x,\bar{x}) \ = \ (0, H(t,x)) \ \ , \ \ 
\phi^{v}_{0} \ = \ id_{H(n)}
\label{ham}
\end{equation}
Then the flow $t \mapsto \phi_{t}$ which satisfies \eqref{hameq} and 
$\phi_{0} \ = \ id_{R^{2n}}$ is the Hamiltonian flow generated by
$H$. 

Conversely, for any Hamiltonian flow $t \mapsto \phi_{t}$, generated 
by $H$, the vertical flow $ t \mapsto \phi^{v}_{t}$ satisfies the
equation \eqref{ham}. 
\label{pham}
\end{prop}

\paragraph{Proof.}
Write the lifts $\tilde{\phi}_{t}$ and $\phi^{h}_{t}$, compute then 
the differential of the quantity 
$\dot{\tilde{\phi}}_{t} - \dot{\phi}^{h}_{t}$ and show that it
equals the differential of $H$. 
\quad $\blacksquare$

\paragraph{Flows of volume preserving diffeomorphisms}

We want to know if there is any nontrivial smooth (according to Pansu differentiability) 
flow of volume preserving diffeomorphisms. 

\begin{prop}
Suppose that $t \mapsto \tilde{\phi}_{t} \in Diff^{2}(H(n),vol)$ is a flow such that 
\begin{enumerate}
\item[-] is $C^{2}$ in the classical sense with respect to $(x,t)$, 
\item[-] is horizontal, that is $t \mapsto \tilde{\phi}_{t}(x)$ is a horizontal curve 
for any $x$. 
\end{enumerate}
Then the flow is constant. 
\label{pho}
\end{prop}

\paragraph{Proof.}
By direct computation, involving second order derivatives. Indeed, let 
$\tilde{\phi}_{t}(x,\bar{x}) \ = \ (\phi_{t}(x), \bar{x} + F_{t}(x))$. 
From the condition $\tilde{\phi}_{t} \in Diff^{2}(H(n),vol)$ we obtain 
$$\frac{\partial F_{t}}{\partial x} y \ = \  \frac{1}{2} \omega(\phi_{t}(x),\frac{\partial 
\phi_{t}}{\partial x}(x) y ) \ - \ \frac{1}{2} \omega(x,y)$$
and from the hypothesis that $t \mapsto \tilde{\phi}_{t}(x)$ is a horizontal curve 
for any $x$ we get 
$$\frac{d F_{t}}{dt}(x) \ = \ \frac{1}{2} \omega(\phi_{t}(x), \dot{\phi}_{t}(x))$$
Equal now the derivative of the RHS of the first relation with respect to $t$ with the derivative of the RHS of the second relation with respect to $x$. We get the equality, for any $y \in R^{2n}$: 
$$ 0  \ = \ \frac{1}{2} \omega(\frac{ \partial \phi_{t}}{\partial x}(x) y, \dot{\phi}_{t}(x))$$
therefore $\dot{\phi}_{t}(x) \ = \ 0$. 
\quad $\blacksquare$

One should expect such a result to be true, based on two remarks. The first, general remark: take a 
flow of left translations in a Carnot group, that is a flow $t \mapsto \phi_{t}(x) \ = \ x_{t} x$. 
We can see directly that each $\phi_{t}$ is smooth, because the distribution is 
left invariant. But the flow is not horizontal, because the distribution is not 
right invariant.  The second, particular remark: any  flow which satisfies the hypothesis 
of proposition \ref{pho} corresponds to a Hamiltonian flow with null Hamiltonian function, hence the flow is constant. 

At a first glance it is  disappointing to see that the group of volume preserving 
diffeomorphisms contains no smooth paths according to the intrinsic calculus 
on Carnot groups. But this makes the richness of such groups of homeomorphisms, as we shall see.

\subsection{Volume preserving bi-Lipschitz maps}

We shall work with the following groups of 
homeomorphisms. 

\begin{defi}
The group $Hom(H(n), vol, Lip)$ is formed by all  
locally bi-Lipschitz, volume preserving homeomorphisms of $H(n)$, which have the form: 
$$\tilde{\phi}(x,\bar{x}) \ = \ (\phi(x), \bar{x} + F(x))$$
 
The group $Sympl(R^{2n}, Lip)$ of locally bi-Lipschitz symplectomorphisms of $R^{2n}$, in the sense that for a.e. $x \in R^{2n}$ the derivative 
$D\phi(x)$ (which exists by classical Rademacher theorem) is symplectic. 

Given $A \subset R^{2n}$, we denote by $Hom(H(n), vol, Lip)(A)$ the group 
of maps $\tilde{\phi}$ which belong to $Hom(H(n), vol, Lip)$ such that $\phi$ has compact support in $A$ (i.e. it differs from identity on a compact set relative to $A$). 

The group  $Sympl(R^{2n}, Lip)(A)$ is defined in an analogous way.  
\end{defi}

Remark that any element $\tilde{\phi} \in  Hom(H(n), vol, Lip)$ preserves 
the "vertical" left invariant distribution  
$$(x \bar{x}) \ \mapsto \ (x,\bar{x}) Z(H(n))$$ 
for a.e. $(x,\bar{x}) \in H(n)$.

We shall prove that any locally bi-Lipschitz volume preserving homeomorphism 
 of $H(n)$ belongs to $Hom(H(n), vol, Lip)$.

\begin{thm}
Take any   $\tilde{\phi}$ locally bi-Lipschitz volume preserving homeomorphism 
 of $H(n)$. Then $\tilde{\phi}$ has the form: 
$$\tilde{\phi}(x,\bar{x}) \ = \ (\phi(x), \bar{x} \ + \ F(x))$$ 
Moreover $$\phi \in Sympl(R^{2n}, Lip)$$
  $F: R^{2n} \rightarrow R$ is 
Lipschitz and for almost any point  $(x,\bar{x}) \in H(n)$ we have:
$$DF(x) y \ = \ \frac{1}{2} \omega(\phi(x), D\phi(x)y) \ - \ \frac{1}{2}\omega(x,y)$$ 
\label{pn1}
\end{thm}

\paragraph{Proof.}
Set $\tilde{\phi}(\tilde{x}) \ = \ (\phi(\tilde{x}), \bar{\phi}(\tilde{x}))$. 
Also denote by 
$$\mid (y,\bar{y}) \mid^{2} \ = \ \mid y\mid^{2} \ + \ \mid\bar{y}\mid$$

By theorem \ref{ppansuu} $\tilde{\phi}$ is almost everywhere derivable and 
the derivative can be written in the particular form: 
$$\left( \begin{array}{cc}
A & 0 \\
0 & 1
       \end{array} \right)$$ with  $A \ = \ A(x,\bar{x}) \  \in Sp(n,R)$.

For a.e. $x \in \mathbb{R}^{2n}$ the function $\tilde{\phi}$ is derivable. The derivative has the form: 
$$D\tilde{\phi}(x,\bar{x}) (y,\bar{y}) \ = \ (A y , \bar{y})$$ 
where by definition 
$$D\tilde{\phi}(x,\bar{x})(y,\bar{y}) \ = \ \lim_{\varepsilon 
\rightarrow 0} \delta_{\varepsilon}^{-1} \left( \tilde{\phi}(x,\bar{x})^{-1} \ \phi((x,\bar{x})\delta_{\varepsilon}(y,\bar{y}))\right)$$
Let us write down what Pansu derivability means.

 For the 
$\mathbb{R}^{2n}$ component we have: the limit  
\begin{equation}
\lim_{\varepsilon \rightarrow 0} \mid 
\frac{1}{\varepsilon} [ \phi((x,\bar{x})\delta_{\varepsilon}(y,\bar{y})) - \phi(x,\bar{x})] - A \mid \ = \ 0
\label{nip1}
\end{equation}
is uniform with respect to $\tilde{y} \in B(0,1)$. We shall rescale all for 
the ball $B(0,R) \in H(n)$. Relation \eqref{nip1} means: for any 
$\lambda > 0$ there exists $\varepsilon_{0}(\lambda)>0$  such that for any 
$R>0$ and any $\varepsilon > \varepsilon_{0}/R$ and $\tilde{y} \in H(n)$ such that 
$\mid \tilde{y} \mid < R$  we have 
\begin{equation}
\mid \frac{1}{\varepsilon} \left( \phi(\tilde{x}\delta_{\varepsilon}\tilde{y}) 
\ - \ \phi(\tilde{x}) \right) \ - \ A y \mid \ < \ \lambda R
\label{nip11}
\end{equation}
Suppose that $\tilde{x} \ = \ (x ,\bar{x})$ satisfies $\mid x \mid \leq R$ with 
$R>1$. 
Fix now $\lambda > 0$ and consider arbitrary  $y \in R^{2n}$, $\mid y \mid < 1$ and $\mu>0$. Set 
$$\bar{y}_{\mu} \ = \ \frac{1}{2\mu} \omega(y,x)$$
There exists a constant $C > 0$ such that 
 $$\mid (y,\bar{y}_{\mu}\mid^{2} \ \leq \ 1 \ + \frac{C}{\mu} (1+R)^{2}$$
For $\varepsilon' < \varepsilon_{0}(\lambda)$ choose 
$$\varepsilon \ = \ \mu \ = \ \frac{1}{4} \left(  
\sqrt{C^{2}(1+R)^{4} + 4 \varepsilon'} \- \ C(1+R)^{2} \right)$$
Then the inequality \eqref{nip11} is true for $R$ replaced 
by $\left( 1 \ + \frac{C}{\mu} (1+R)^{2}\right)^{1/2}$, $\tilde{y}= (y, \bar{y}_{\mu})$, with arbitrary $\mid y \mid < 1$. We get: for any $\lambda > 0$ 
there is $\varepsilon(\lambda) = \varepsilon_{0}(\lambda)/2 >0$ such that for any $\varepsilon < \varepsilon(\lambda)$ and any $\mid y \mid < 1$ we have 
$$\mid \frac{1}{\varepsilon} \left( \phi(x+\varepsilon y, \bar{x}) - \phi(x,\bar{x}) \right) \ - A y \mid \ < \ 2 \lambda$$
Therefore $x \mapsto \phi(x,\bar{x})$ is derivable a.e. with (locally) bounded derivative. It is therefore (locally) Lipschitz. 

Let us now return to \eqref{nip1} and choose $\bar{y} = 0$. We get: 
the map 
$$y \mapsto \phi(x+y, \bar{x} + \frac{1}{2}\omega(x,y))$$
is derivable in $y = 0$ and the derivative equals $A$. But we have just proven that the derivative of $\phi(\cdot, \bar{x})$ in $x$ equals $A$. We get: 
$\phi$ is constant with respect to $\bar{x}$. 

We have proven until now that 
$$\tilde{\phi}(x,\bar{x}) \ = \ (\phi(x), \bar{\phi}(x,\bar{x}))$$
with $\phi \in Sympl(R^{2n}, Lip)$.

Let us use this information in the definition of the center component of the Pansu derivative. We have: 
\begin{equation}
\mid \frac{1}{\varepsilon^{2}} \left( \bar{\phi}(\tilde{x} \delta_{\varepsilon} 
\tilde{y}) - \bar{\phi}(\tilde{x}) - \frac{1}{2} \omega( \phi(x), \phi(x+ \varepsilon y)) \right) - \bar{y} - 
\label{vuvu}
\end{equation}
$$ - \frac{1}{2} \omega(\frac{1}{\varepsilon} (\phi(x+\varepsilon y) - \phi(x)), 
Ay) \mid \ \rightarrow 0$$
as $\varepsilon \rightarrow 0$, uniformly with respect to $\mid \tilde{y} \mid < 1$. From previously proved facts we obtain that \eqref{vuvu} is equivalent 
to 
$$\mid \frac{1}{\varepsilon^{2}} \left( \bar{\phi}(\tilde{x} \delta_{\varepsilon} 
\tilde{y}) - \bar{\phi}(\tilde{x}) - \frac{1}{2} \omega( \phi(x), \phi(x+ \varepsilon y)- \phi(x)) \right) - \bar{y}\mid \ \rightarrow 0$$
as $\varepsilon \rightarrow 0$, uniformly with respect to $\tilde{y}$. 
We quickly get that 
$\bar{\phi}(x,\bar{x}) \ = \ \bar{x} \ + \ F(x)$.
 
All in all the Pansu derivability (center component) reads  
\begin{equation}
\lim_{\varepsilon \rightarrow 0} 
\mid \frac{1}{\varepsilon^{2}} [F(x + \varepsilon y) - F(x)] + \frac{1}{2\varepsilon}\omega(x,y) - 
\frac{1}{2} \omega(\phi(x), \frac{1}{\varepsilon^{2}} [ \phi(x + \varepsilon y) - \phi(x)]) \mid \ = \ 0
\label{nip2}
\end{equation}
This means that for almost any $x \in \mathbb{R}^{2n}$ 
\begin{enumerate}
\item[1] the function 
$\phi$ is derivable and the derivative is equal to $A \ = \ A(x)$,
\item[2] the function $F$ is derivable and is connected to $\phi$ by the 
relation from the conclusion of the proposition. 
\end{enumerate}

We have seen that   
$\tilde{\phi}$ is locally Lipschitz implies that $x \in \mathbb{R}^{2n} 
\mapsto A(x)$ is locally bounded, therefore $\phi$ is locally Lipschitz. 
In the same way we obtain that $F$ is locally Lipschitz. 
\quad $\blacksquare$

The  flows in the group $Hom(H(n), vol, Lip)$ are defined further. 

\begin{defi}

A flow in the group $Hom(H(n),vol, Lip)$ is a curve  
$$t \mapsto \tilde{\phi}_{t} \in Hom(H(n), vol, Lip)$$ such that for 
a.e. $x \in R^{2n}$ the curve 
$t \mapsto \phi_{t}(x) \in R^{2n}$ is (locally) Lipschitz. 

For any flow we can define the horizontal lift of this flow like this: a.e. curve 
$t \mapsto \phi_{t}(x)$ lifts to a horizontal curve 
$t \mapsto \phi_{t}^{h}(x,0)$. Define then 
$$\phi_{t}^{h}(x,\bar{x}) \ = \ \phi_{t}^{h}(x,0) (0,\bar{x})$$
\label{dflo}
\end{defi}

If the flow is smooth (in the classical sense) and $\tilde{\phi}_{t} 
\in Diff^{2}(H(n),vol)$ then this lift is the same as the one described in 
definition \ref{dhameq}. We can define now the vertical flow by the formula
\eqref{hameq}, that is
$$\phi^{v} \ = \ \tilde{\phi}^{-1} \circ \phi^{h}$$

There is an  analog of proposition \ref{pho}. In the proof we shall need 
lemma \ref{pglnr}, which comes after.

\begin{prop}
Let $t \mapsto \tilde{\phi}_{t} \in Hom(H(n), vol, Lip)$ be a curve such that  the function  $t \mapsto 
\Phi(\tilde{x},t) = (\tilde{\phi}_{t}(\tilde{x}), t)$ is locally Lipschitz from $H(n) \times R$ to itself. 
Then $t \mapsto \tilde{\phi}_{t}$ is a constant curve. 
\label{tn2}
\end{prop}

\paragraph{Proof.}
By Rademacher theorem \ref{ppansu} for the group $H(n) \times R$ we obtain that 
$\Phi$ is almost everywhere derivable. Use now lemma \ref{pglnr} to deduce the claim. 
\quad $\blacksquare$

A short preparation is needed in order to state the lemma \ref{pglnr}. Let $N$ be a noncommutative Carnot group. 
We shall  look at the  group $N \times R$ with the group operation defined component wise. 
This is also a Carnot group. Indeed, consider the family of dilatations 
$$\delta_{\varepsilon}(x,t) \ = \ (\delta_{\varepsilon}(x), \varepsilon t)$$
which gives to $N \times R$ the structure of a CC group. The left invariant distribution 
on the group which generates the distance is (the left translation of)  $W_{1} = V_{1} 
\times R$.

\begin{lema}
Let $N$ be a noncommutative Carnot group which admits the orthogonal decomposition 
$$N \ = \ V_{1} + [N,N]$$ and satisfies the condition 
$$V_{1} \cap Z(N) \ = \ 0$$  The group of linear transformations of $N\times R$ is then 
$$HL(N \times R) \ = \ \left\{ \left( 
\begin{array}{cc}
A & 0 \\
c & d
       \end{array} \right) \ \mbox{ : } A \in HL(N) \ , \ c \in V_{1} \ , \ d \in R 
\right\}$$
\label{pglnr}
\end{lema}

\paragraph{Proof.}
We shall proceed as in the proof of proposition \ref{p1}. We are looking first at the Lie algebra isomorphisms of $N \times R$, with general form 
$$\left( 
\begin{array}{cc}
A & b \\
c & d
       \end{array} \right) $$
We obtain the conditions: 
\begin{enumerate}
\item[(i)] $c$ orthogonal on $[N,N]$, 
\item[(ii)] $b$ commutes with the image of $A$:  [b,Ay] = 0, for any $y \in N$, 
\item[(iii)] $A$ is an algebra isomorphism of $N$. 
\end{enumerate}
From (ii), (iii) we deduce that $b$ is in the center  of $N$ and from (i) we see that $c \in V_{1}$. 

We want now  the isomorphism to commute with dilatations. This condition gives: 
\begin{enumerate}
\item[(iv)] $b \in V_{1}$, 
\item[(v)] $A$ commutes with the dilatations of $N$. 
\end{enumerate}
(iii) and (v) imply that $A \in HL(N)$ and (iv) that $b = 0$. 
\quad $\blacksquare$

\paragraph{Hamiltonian diffeomorphisms: more structure}

In this section we look closer to the structure of the group of volume preserving homeomorphisms of the Heisenberg group.

Let $Hom^{h}(H(n),vol, Lip)(A)$ be the group of time one homeomorphisms 
$t \mapsto \phi^{h}$, for all curves $t \mapsto \phi_{t} \in 
Sympl(R^{2n},Lip)(A)$ such that $(x,t) \mapsto (\phi_{t}(x),t)$ is locally 
Lipschitz.

The elements of this group are also volume preserving, but they are not 
smooth with respect to the Pansu derivative.

\begin{defi}
The group  $Hom(H(n), vol)(A)$ contains all maps $\tilde{\phi}$ which have 
a.e. the  form: 
$$\tilde{\phi}(x,\bar{x}) \ = \  (\phi(x), \bar{x} + F(x))$$
where $\phi \in Sympl(R^{2n},Lip)(A)$ and $F: R^{2n} \rightarrow R$ 
is locally Lipschitz and constant outside a compact set included in the closure of $A$ (for short: with compact support in $A$). 
\label{pn2}
\end{defi}

This group contains  three privileged subgroups: 
\begin{enumerate}
\item[-] $Hom(H(n),vol,Lip)(A)$, 
\item[-] $Hom^{h}(H(n),vol, Lip)(A)$
\item[-]  and $Hom^{v}(H(n), vol, Lip)(A)$.
\end{enumerate} 
The last is the group of vertical homeomorphisms, any of which has the form: 
$$\phi^{v}(x,\bar{x}) \ = \ (x, \bar{x} + F(x))$$
with $F$ locally Lipschitz, with compact support in $A$.

Take a one parameter subgroup $t \mapsto \tilde{\phi}_{t} \in Hom(H(n),vol,Lip)(A)$, in the sense of the definition \ref{dflo}. We know that it cannot be smooth as a curve in in the group $Hom(H(n),vol,Lip)(A)$, but we also know that there are vertical and horizontal flows 
$t \mapsto \phi_{t}^{v} \ , \ \phi_{t}^{h}$ such that 
we have the decomposition $\tilde{\phi}_{t} \circ \phi_{t}^{v} \ =  \ \phi_{t}^{h}$.  
Unfortunately none of the flows $t \mapsto \phi_{t}^{v} \ , \ \phi_{t}^{h}$ are one parameter groups.

There is more structure here that it seems. Consider the 
class 
$$HAM(H(n))(A) = Hom(H(n), vol,Lip)(A) \times Hom^{v}(H(n),vol, Lip)(A)$$
For any pair in this class we shall use the notation 
$(\tilde{\phi}, \phi^{v})$ 
This class forms a group with the (semidirect product) operation: 
$$(\tilde{\phi}, \phi^{v}) (\tilde{\psi}, \psi^{v}) \ = \ 
( \tilde{\phi} \circ \tilde{\psi}, \phi^{v} \circ \tilde{\phi} \circ 
\psi^{v} \circ \tilde{\phi}^{-1})$$

\begin{prop}
If $t \mapsto \tilde{\phi}_{t} \in Hom(H(n),vol,Lip)(A)$ is an one parameter group then 
$t \mapsto (\tilde{\phi_{t}}, \phi_{t}^{v}) \in HAM(H(n))$ is an one parameter group. 
\label{popg}
\end{prop}

\paragraph{Proof.}
The check is left to the reader. Use definition and proposition 1, chapter 5, Hofer \& 
Zehnder \cite{hozen}, page 144. 
\quad $\blacksquare$

We can (indirectly) put a distribution on the group $HAM(H(n))(A)$ by specifying the class of horizontal curves. Such a
 curve $t \mapsto (\tilde{\phi_{t}}, \phi^{v}_{t})$ in the group
$HAM(H(n))(A)$ projects in the first component $t \mapsto \tilde{\phi}_{t}$ to a flow in $Hom(H(n),vol, Lip)(A)$.  Define for any  $(\tilde{\phi},
\phi^{v}) \in HAM(H(n))(A)$ the function 
$$\phi^{h} \ = \ \tilde{\phi} \circ \phi^{v}$$
and  say   that 
the  curve $$t \mapsto (\tilde{\phi}_{t}, \phi_{t}^{v}) \in HAM(H(n))(A)$$  
is horizontal if 
$$t \mapsto \phi^{h}_{t}(x,\bar{x})$$ is horizontal for any 
$(x,\bar{x}) \in H(n)$.

We introduce the following length function for horizontal curves: 
$$L \left(t \mapsto (\tilde{\phi}_{t}, \phi_{t}^{v})\right) \ = \ 
\int_{0}^{1} \| \dot{\phi_{t}^{v}} \|_{L^{\infty}(A)} \mbox{ d}t $$
With the help of the length we are able to endow the group $HAM(H(n))(A)$ as a 
path metric space. The distance is defined by: 
$$dist\left((\tilde{\phi}_{I}, \phi^{v}_{I}), (\tilde{\phi}_{II}, \phi^{v}_{II})\right) 
\ = \ \inf  L\left(t \mapsto (\tilde{\phi}_{t}, \phi_{t}^{v})\right)$$
over all horizontal curves  $t \mapsto (\tilde{\phi}_{t}, \phi_{t}^{v})$ such that 
$$(\tilde{\phi}_{0}, \phi_{0}^{v}) \ = \ (\tilde{\phi}_{I}, \phi_{I}^{v})$$
$$(\tilde{\phi}_{1}, \phi_{1}^{v}) \ = \ (\tilde{\phi}_{II}, \phi_{II}^{v})$$
The distance is not defined for any two points in $HAM(H(n))(A)$. In principle 
the distance can be degenerated.

The group $HAM(H(n))(A)$ acts on $AZ$ (where $Z$ is the center of 
$H(n)$) by: 
$$(\tilde{\phi}, \phi^{v}) (x, \bar{x}) \ = \ \phi^{v} \circ \tilde{\phi} (x, \bar{x})$$

Remark  that the group 
$$Hom^{v}(H(n), vol, Lip)(A) \ \equiv \left\{ id \right\} \times Hom^{v}(H(n),vol, Lip)(A)$$
 is normal in $HAM(H(n))(A)$. We denote by $$HAM(H(n))(A)/Hom^{v}(H(n),vol, Lip)(A)$$ the factor group. This is isomorphic with the group $Sympl(R^{2n},Lip)(A)$. 

Consider now the natural action of $Hom^{v}(H(n),vol, Lip)(A)$ on $AZ$. The space of orbits $AZ/Hom^{v}(H(n),vol,Lip)(A)$ is nothing but 
$A$, with the Euclidean metric. 

Factorize now the action of $HAM(H(n))(A)$ on $AZ$, by the group 
$$Hom^{v}(H(n),vol,Lip)(A)$$ We obtain the action of  
$HAM(H(n))(A)/Hom^{v}(H(n),vol,Lip)(A)$ on  $$AZ/Hom^{v}(H(n),vol,Lip)(A)$$  Elementary examination shows this is  the action of the symplectomorphisms group on $R^{2n}$.

All in all we have the following (for the definition of the Hofer distance see next section): 

\begin{thm}
The action of $HAM(H(n))(A)$ on $AZ$ descends after reduction with the group $Hom^{v}(H(n), vol,Lip)$ 
to the action of symplectomorphisms group with compact support in $A$.
 
The distance $dist$ on $HAM(H(n))$ descends to the Hofer distance on the connected 
component of identity of the symplectomorphisms group. 
\label{tmai}
\end{thm}

\subsection{Symplectomorphisms, capacities and Hofer distance}

Symplectic capacities are invariants under the action of the symplectomorphisms group. 
Hofer geometry is the geometry of the group of Hamiltonian diffeomorphisms, with respect 
to the Hofer distance. For an introduction into the subject see Hofer, Zehnder
\cite{hozen} chapters 2,3 and 5, and Polterovich \cite{polte}, chapters 1,2. 

 A symplectic capacity is a map which associates to any symplectic 
manifold $(M,\omega)$ a number $c(M,\omega) \in [0,+\infty]$. Symplectic capacities 
are special cases of conformal symplectic invariants, described by: 
\begin{enumerate}
\item[A1.] Monotonicity: $c(M,\omega) \leq c(N,\tau)$ if there is a symplectic 
embedding  from $M$ to $N$, 
\item[A2.] Conformality: $c(M,\varepsilon \omega) = \mid \varepsilon \mid c(M,\omega)$ 
for any $\alpha \in R$, $\alpha \not = 0$. 
\end{enumerate}

We can see a conformal symplectic invariant from another point of view. Take 
a symplectic manifold $(M,\omega)$ and consider the invariant defined over the class 
of Borel sets $B(M)$,  (seen as embedded submanifolds). In the particular case of 
$R^{2n}$ with the standard symplectic form, an invariant is a function 
$c: B(R^{2n}) \rightarrow [0,+\infty]$ such that: 
\begin{enumerate}
\item[B1.] Monotonicity: $c(M) \leq c(N)$ if there is a symplectomorphism $\phi$ such 
that $\phi(M) \subset N$,
\item[B2.] Conformality: $c(\varepsilon M ) = \varepsilon^{2} c(M)$ for any 
$\varepsilon \in R$. 
\end{enumerate} 

An invariant is nontrivial if it takes finite values on sets with infinite volume, 
like cylinders: 
$$Z(R) = \left\{ x \in R^{2n} \mbox{ : } x^{2}_{1} + x_{2}^{2} < R \right\}$$

There exist highly nontrivial invariants, as the following theorem shows: 

\begin{thm} (Gromov's squeezing theorem) The ball $B(r)$ can be symplectically embedded in the cylinder $Z(R)$   if and only if $r \leq R$. 
\end{thm}

This theorem permits to define the invariant: 
$$c(A) \ = \ \sup \left\{ R^{2} \mbox{ : } \exists \phi(B(R)) \subset A \right\}$$
called Gromov's capacity. 

Another important invariant is Hofer-Zehnder capacity. In order to introduce this 
we need the notion of a Hamiltonian flow.  

A flow of symplectomorphisms $t \mapsto \phi_{t}$ is Hamiltonian if there is 
a function $H: M \times R \rightarrow R$ such that for any time $t$ and place $x$ we have 
$$\omega( \dot{\phi}_{t}(x), v) \ = \ dH(\phi_{t}(x),t) v$$
for any $v \in T_{\phi_{t}(x)}M$. 

Let $H(R^{2n})$ be the set of compactly supported Hamiltonians. Given 
a set $A \subset R^{2n}$, the class of admissible Hamiltonians is $H(A)$, made by all  compactly supported maps in $A$ such that the generated  Hamiltonian flow does not have closed orbits of periods smaller than 1. Then the Hofer-Zehnder capacity 
is defined by: 
$$hz(A) \ = \ \sup \left\{ \| H \|_{\infty} \mbox{ : } H \in H(A) \right\}$$

Let us denote by $Ham(A)$ the class of Hamiltonian diffeomorphisms compactly supported in $A$. A Hamiltonian diffeomorphism is the time one value of a Hamiltonian flow. In the case which interest us, that is $R^{2n}$, $Ham(A)$ is the connected component of the identity in the group of compactly supported symplectomorphisms. 

A curve of Hamiltonian diffeomorphisms (with compact support) is a Hamiltonian flow. 
For any such curve $t \mapsto c(t)$ we shall denote by $t\mapsto H_{c}(t, \cdot)$ the associated Hamiltonian  function (with compact support). 

On the group of Hamiltonian diffeomorphisms there is a bi-invariant distance introduced by Hofer. This is given by the expression: 
\begin{equation}
d_{H}(\phi,\psi) \ = \ \inf \left\{ \int_{0}^{1} \| H_{c}(t) \|_{\infty, R^{2n}} \mbox{ d}t \mbox{ : } c: [0,1] \rightarrow Ham(R^{2n}) \right\}
\label{hodis}
\end{equation}

It is easy to check that $d$ is indeed bi-invariant and it satisfies the triangle property. It is a deep result that $d$ is non-degenerate, that is 
$d(id, \phi) \ = \ 0 $ implies $\phi = \ id$. 

With the help of the Hofer distance one can define another symplectic invariant, called displacement energy. For a set $A \subset R^{2n}$ the displacement energy is: 
$$de(A) \ = \ \inf \left\{ d_{H}(id, \phi) \ \mbox{ : } \phi \in Ham(R^{2n}) \ \ , \ \phi(A) \cap A = \emptyset \right\}$$

A displacement energy can be associated to any group of transformations endowed with a bi-invariant distance (see Eliashberg \& Polterovich \cite{elia}, section 1.3). The fact that the displacement energy is a nontrivial invariant is equivalent with the 
non-degeneracy of the Hofer distance. Section 2.  Hofer \& Zehnder \cite{hozen} 
is dedicated to the implications of the existence of a non-trivial capacity. 	All in all the non-degeneracy of the Hofer distance (proved for example in the Section 5. Hofer \& Zehnder \cite{hozen}) is the cornerstone 
of symplectic rigidity theory.

\subsection{Hausdorff dimension and Hofer distance}

We shall give in this section a metric proof of the non-degeneracy of the 
Hofer distance. 

A flow of symplectomorphisms $t \mapsto \phi_{t}$  with compact support in $A$ is Hamiltonian if it is the projection onto $R^{2n}$ of a flow in 
$Hom(H(n),vol, Lip)(A)$.   

Consider such a flow which joins identity $id$ with $\phi$. Take the lift 
of the flow $t \in [0,1] \mapsto \tilde{\phi}_{t} \in Hom(H(n),vol, Lip)(A)$. 

\begin{prop}
The curve $t \mapsto \tilde{\phi}_{t}(x,0)$ has Hausdorff dimension $2$ and measure 
\begin{equation}
\mathcal{H}^{2} \left( t \mapsto \tilde{\phi}_{t}(x,\bar{x})\right) \ = \ \int_{0}^{1} \mid H_{t}(\phi_{t}(x)) \mid \mbox{ d}t
\label{aham}
\end{equation}
\label{ppham}
\end{prop}

\paragraph{Proof.}
The curve is not horizontal. The tangent has a vertical part (equal to the Hamiltonian). Use the definition of the (spherical) Hausdorff measure 
$\mathcal{H}^{2}$ and  the Ball-Box theorem \ref{bbprop} to obtain the formula \eqref{aham}. 
\quad $\blacksquare$

We shall prove now that the Hofer distance \eqref{hodis} is non-degenerate. 
For given $\phi$ with compact support in $A$ and generating function 
$F$, look at the one parameter family: 
$$\tilde{\phi}^{a}(x,\bar{x})  \ = \ (\phi(x), \bar{x} + F(x) + a)$$
The volume of the cylinder 
$$\left\{ (x,z) \mbox{ : } x \in A , \ z \mbox{ between } 0 \mbox{ and } 
F(x)+a \right\}$$
attains the minimum   $V(\phi,A)$ for an $a_{0} \in R$. 

Given an arbitrary flow $t \in [0,1] \mapsto \tilde{\phi}_{t} \in Hom(H(n),vol, Lip)(A)$ such that $\tilde{\phi}_{0} \ = \ id$ and 
$\tilde{\phi}_{1} \ = \ \tilde{\phi}^{a}$, we have the following inequality: 
\begin{equation}
V(\phi,A) \ \leq \ \int_{A} \mathcal{H}^{2}(\left( t \mapsto \tilde{\phi}_{t}(x,0)\right) \mbox{ d}x 
\label{ineq1}
\end{equation}
Indeed, the family of curves $t \mapsto \tilde{\phi}_{t}(x,0)$ provides a foliation of the set 
$$\left\{ \tilde{\phi}_{t}(x,0) \mbox{ : } x \in A \right\}$$ 
The volume of this set is lesser than the RHS of inequality \eqref{ineq1}: 
\begin{equation}
 C \  Vol \left( \left\{ \tilde{\phi}_{t}(x,0) \mbox{ : } x \in A \right\}\right) \ \leq \ \int_{A} \mathcal{H}^{2} \left( t \mapsto \tilde{\phi}_{t}(x,0)\right) \mbox{ d}x
\label{cla}
\end{equation}
where $C>0$ is an universal positive constant (that is it depends only on the 
algebraic structure of $H(n)$). 
This is a statement which is local in nature; it is simply saying that 
the area of the parallelogram defined by the vectors $\mathbf{a}$, 
$\mathbf{b}$ is smaller than the product of the norms $\|\mathbf{a}\| 
\|\mathbf{b}\|$. We give further a  proof of this inequality.  

Let $\tilde{A} \subset H(n)$ be a set with finite Hausdorff measure $Q = 2n+2$  which contains a set 
$\tilde{B} \ = \  B \times \left\{ 0 \right\}$  such that there is a flow 
$$t \in [0,1] \mapsto \phi_{t}: H(n) \rightarrow H(n)$$ 
with the properties: 
\begin{enumerate}
\item[(a)] $\phi_{0}(x) \ = \ x$ for any  $x \in B$,  
\item[(b)] for almost any $t$ $\phi_{t}$ preserves the Hausdorff measure $Q$,
\item[(c)] for any $t \in [0,1]$ $\phi_{t}$ is bi-Lipschitz and the Lipschitz 
constants of $\phi_{t}$, $\phi_{t}^{-1}$ are upper bounded by a positive constant $Const$,  
\item[(d)] $\tilde{A} \ = \ \left\{ \phi_{t}(x) \mbox{ : } t \in [0,1] \ , \ x \in B \right\}$, 
\item[(e)] for $\mathcal{H}^{2n}$ almost any $x \in B$ the fiber $f_{x} \ = \ \left\{ \phi_{t}(x) 
\mbox{ : } t \in [0,1] \right\}$ has    finite Hausdorff dimension $P= 2$. 
\end{enumerate}
 Then 
$$C \ \mathcal{H}^{Q}(\tilde{A}) \ \leq \ \int_{B} \mathcal{H}^{2}(f_{x}) \mbox{ d}(x)$$

Indeed, we start from the isodiametric inequality in the Heisenberg group: 
$$\left(  diam \ \ \tilde{D} \right)^{Q} \ \geq \ c \ \mathcal{H}^{Q}(\tilde{D})$$
We can rescale the Hausdorff measure $Q$ in order to have $c=1$. Consider also the projection  
$\pi: H(n) \rightarrow R^{2n}$, $\pi(x,\bar{x}) \ = \ x$. We can rescale the 
Lebesgue measure on $R^{2n}$ such that the volume of the projection of the 
ball $B^{CC}(\tilde{x}, 1)$ equals $1$. Here we have used the notation: $B^{CC}(\tilde{x},R)$ is the ball in the Carnot-Carath\'eodory distance in $H(n)$, and $B^{E}(x,R)$ is the Euclidean ball of radius $R$ in $R^{2n}$. 

Fix $\varepsilon \ 0$ small and cover $B$ with Euclidean balls of radius  $\varepsilon$ $B^{E}(x_{j}, \varepsilon)$, $j \in J_{\varepsilon}$, in an efficient way, i.e. such that $J_{\varepsilon} \varepsilon^{2n}$ tends to 
$vol(B)$ when $\varepsilon \rightarrow 0$. We want now to cover $\tilde{A}$ with sets of the form 
$$\tilde{\phi}_{t_{k}}(B^{CC}((x_{j},0), \varepsilon))$$
where $k$ depends on $x_{j}$. How many sets do we need?  

From the isodiametric inequality and the fact that $\tilde{\phi}_{t}$ is 
volume preserving,  we see that 
$$diam \ \tilde{\phi}_{t}(B^{CC}((x_{j},0), \varepsilon))  \geq \varepsilon$$
Put this together with the fact that the fiber $f_{x_{j}}$ has Hausdorff dimension $P=2$ and get that we need about   $\mathcal{H}^{P}(f_{x_{j}})/\varepsilon^{P}$ such sets to cover the fiber $f_{x_{j}}$. The sum of volumes of the cover is greater then the volume of $A$. Tend $\varepsilon$ to $0$ and obtain the claimed inequality \eqref{cla}. The constant $C$ comes from the rescalings of the measures. Without rescalling we have: 
$$C \ = \ \frac{vol(B^{CC}(0,1))}{vol(B^{E}(0,1)}$$

We go back now to the main stream of the proof and remark that 
by continuity of the flow the volume of the set 
$$\left\{ \tilde{\phi}_{t}(x,0) \mbox{ : } x \in A \right\}$$ 
is  greater than $V(\phi,A)$. 

For any curve of the flow we have the uniform obvious estimate: 
\begin{equation}
\mathcal{H}^{2}(\left( t \mapsto \tilde{\phi}_{t}(x,0)\right) \ \leq \ 
\int_{0}^{1} \| H_{t}(\phi_{t}(x)) \|_{A,\infty} \mbox{ d}t
\label{ineq2}
\end{equation}
Put together inequalities \eqref{ineq1}, \eqref{ineq2} and get 
\begin{equation}
C \ V(\phi,A) \ \leq \ \int_{0}^{1} \| H_{t}(\phi_{t}(x)) \|_{A,\infty} \mbox{ d}t
\label{inef}
\end{equation}

 Use the definition of the Hofer distance \eqref{hodis} to obtain the inequality 
\begin{equation}
C \ V(\phi,A) \ \leq \ vol(A) \ d_{H}(id, \phi)
\label{inef1}
\end{equation}
This proves the non-degeneracy of the Hofer distance, because if the 
RHS of \eqref{inef1} equals 0 then $V(\phi,A)$ is 0, which means that the generating function of $\phi$ is almost everywhere constant, therefore 
$\phi$ is the identity everywhere in $R^{2n}$. 

We close with the translation of the inequality \eqref{inef} in symplectic terms. 

\begin{prop}
Let $\phi$ be a Hamiltonian diffeomorphism with compact support in 
$A$ and $F$ its generating function, that is $dF \ = \ \phi^{*} \lambda - 
\lambda$, where $d\lambda \ = \omega$. Consider a Hamiltonian flow 
$t \mapsto H_{t}$, with compact support in $A$, such that the time one 
map equals $\phi$. Then the following inequality holds: 
$$ C \ \inf \left\{ \int_{A} \mid F(x) - c \mid \mbox{ d}x \mbox{ : } c \in 
R \right\} \ \leq \ vol(A) \ \int_{0}^{1} \| H_{t}\|_{\infty,A} \mbox{ d}t$$
with $C >0$ an universal constant equal to 
$$C \ = \ \frac{vol(B^{CC}(0,1))}{vol(B^{E}(0,1)}$$
\end{prop}

\subsection{Invariants of volume preserving maps}

Theorem \ref{pn1} gives the strucure of a volume preserving locally bi-Lipschitz homeomorphism. We have denoted by $Hom(H(n), vol , Lip)$ the class 
of these maps. Any $\tilde{\phi} \in Hom(H(n), vol , Lip)$ can be written as 
$$\tilde{\phi}(x,\bar{x}) \ = \ (\phi(x), \bar{x} + F(x))$$ where 
$\phi \in Sympl(R^{2n}, Lip)$ is a locally bi-Lipschitz symplectomorphism of 
$R^{2n}$. For a set $B \subset R^{2n}$ we denote by $vol(B)$ it's Lebesgue 
measure. For a set $C \subset R$ $l(C)$ will be it's length (one dimensional Lebesgue measure). An immediate consequence of the structure theorem \ref{pn1} is given further.

\begin{cor}
Let $\tilde{A} \subset H(n)$ be a $\mathcal{H}^{2n+2}$ measurable set such that 
$\mathcal{H}^{2n}(A) \ < \ + \infty$, where $A$ denotes the orthogonal projection of $\tilde{A}$ on $R^{2n}$. Then for any $\tilde{\phi} \in 
Hom(H(n), vol , Lip)$ we have: 
$$w(\tilde{A}) \ = \ \frac{\mathcal{H}^{2n+2}(\tilde{A})}{vol (A)} \ = \ \frac{\mathcal{H}^{2n+2}(\tilde{\phi}(\tilde{A}))}{vol (\phi(A))}$$
We call $w(\tilde{A})$ the width of $\tilde{A}$.  

More general, for any $i \geq 1$ define the $i$-height of $\tilde{A}$ to be 
$$\left(h_{i}(\tilde{A})\right)^{i} \ = \ \frac{1}{vol(A)^{i}} \int_{A} 
\left( l(\left\{\bar{x} \in R \mbox{ : } (x,\bar{x}) \in \tilde{A} \right\}) \right)^{1/i} \mbox{ d}x$$
Then $h_{i}(\tilde{A}) \ = \ h_{i}(\tilde{\phi}(\tilde{A}))$. 
 \label{cev}
\end{cor}

\paragraph{Proof.}
Any symplectomorphism is volume preserving. The thesis follows from this and the structure theorem \ref{pn1}. 
\quad $\blacksquare$

Suppose $\tilde{A}$ is open, bounded and classically smooth. In the case 
of the Heisenberg group $H(1)$ Pansu \cite{paninis} proved the following 
isoperimetric inequality: 
$$\mathcal{H}^{4}(\tilde{A}) \ \leq \  C \left( \mathcal{H}^{3}(\partial \tilde{A} 
\right)^{\frac{4}{3}}$$

The isoperimetric inequality is a very important subject in analysis in metric spaces. It has been intensively studied in  it's disguised form (Sobolev inequalities) in many papers. To have a grasp of what is hidden here the reader 
can consult Varopoulos 
\& al. \cite{varsaco}, Hajlasz, Koskela \cite{hako} and  Garofalo, Nhieu \cite{ganh}.  
In order to give a decent exposition of  this subject we would need a good notion of perimeter. In this first part we don't enter much into the measure theoretic aspects though.   

Let us work with objects which have classical regularity, as in Gromov \cite{gromo} section 2.3. There  we find a more general isoperimetric inequality, true in particular in any Carnot group, provided that $\tilde{A}$ is open, bounded, with (classically) smooth boundary. If we denote by $Q$ the homogeneous dimension of the Carnot group $N$ (which equals its Hausdorff dimension) then Gromov's isoperimetric inequality is: 
$$\mathcal{H}^{Q}(\tilde{A}) \ \leq \ C \left( \mathcal{H}^{Q-1}(\partial \tilde{A}) \right)^{\frac{Q}{Q-1}}$$
Recall that for Heisenberg group $H(n)$ the homogeneous dimension is 
$Q = 2n+2$. 

We should  use further the group $Diff^{2}(H(n),vol)$ instead of  $Hom(H(n), vol , Lip)$, because we work with classical regularity  . We kindly ask the reader to believe that  a good definition of the perimeter of the set 
$\tilde{A}$ allows us to use the group $Hom(H(n), vol , Lip)$. We shall denote 
the perimeter function by $Per$. In "smooth" situations $Per(\partial \tilde{A} ) \ = \ \mathcal{H}^{Q-1}(\partial \tilde{A})$.

We can define then the isoperimeter of $\tilde{A}$ to be: 
$$isop(\tilde{A})^{2n+1} \ = \ \inf  \frac{\left( Per \ \partial \tilde{\phi}(\tilde{A} ) \right)^{2n+2}}{\left(\mathcal{H}^{2n}(A) \right)^{2n+1}}$$
with respect to all $\tilde{\phi} \in Hom(H(n), vol, Lip)$. This is another 
invariant of the group of volume preserving locally bi-Lipschitz maps. 

Yet another invariant is the isodiameter. This is defined as 
$$isod(\tilde{A}) \ = \ \inf \frac{(diam \ \tilde{\phi}(\tilde{A}))^{2n+2}}{\mathcal{H}^{2n}(A)} $$
with respect to all $\tilde{\phi} \in Hom(H(n), vol, Lip)$.
Because the diameter of $\tilde{A}$ is greater than the square root of 
the $\infty$ height of $\tilde{A}$, it follows that 
$$isod(\tilde{A}) \ \geq \ \frac{h_{\infty}(\tilde{A})^{n+1} \ w(\tilde{A})}{ \mathcal{H}^{2n+2}(\tilde{A})}$$

All these invariants scale with dilatations like this (denote generically  by $inv$ any of the mentioned invariants): 
$$inv \ \delta_{\varepsilon} \tilde{A} \ = \ \varepsilon^{2} \tilde{A}$$

The groupoid 
$\mathcal{H}om(H(n), vol, Lip)$ with objects open sets $H(n)$ and arrows 
volume preserving locally bi-Lipschitz maps between then, can be considered instead of the group $Hom(H(n), vol, Lip)$. Then the  
heights are no longer invariants. Indeed, set $\tilde{A}$ to be the disjoint union of two isometric sets. Then we can keep one of these sets fixed and put 
on the top of it the other one, using only translations. The width will change.  We can define  other invariants though, using the heights and the number of connected components, for example. 

By using lifts of symplectic flows one can define symplectic invariants. Conversely, 
with any symplectic capacity comes an invariant of $Hom(H(n), vol, Lip)$. Indeed, let $c$ be a capacity. Then $\tilde{c}(\tilde{A}) \ = \ c(A)$ is an 
invariant. 

Is there any general definition of a $Hom(H(n), vol, Lip)$ invariant which has as particular realization heights and 
symplectic capacities? 

How about the invariants of the  group $Hom(N, vol, Lip)$, where $N$ is a general Carnot group? We leave this for further study.

\section{Sub-Riemannian Lie groups}

In this section we shall look closer to the last example of section 
\ref{exemple}.  Let $G$ be a real 
 connected Lie group with Lie algebra $\mathfrak{g}$ and $D \subset \mathfrak{g}$ 
a vector space which generates the algebra. This means that the sequence 
$$V_{1} \ = \ D \ \ , \ V^{i+1}  \ = \ V^{i} + [D,V^{i}]$$
provides a filtration of $\mathfrak{g}$:
\begin{equation} 
V^{1} \ \subset \ V^{2} \ \subset \ ... \ V^{m} \ = \ \mathfrak{g}
\label{filtra}
\end{equation}

\subsection{Nilpotentisation}
\label{nilpo}

The filtration has the straightforward property: if $x \in V^{i}$ and 
$y \in V^{j}$ then $[x,y] \in V^{i+j}$ (where $V^{k} = \mathfrak{g}$ for all $k \geq m$ and $V^{0} = \left\{ 0 \right\}$). This allows to construct the Lie algebra 
$$\mathfrak{n}(\mathfrak{g}, D) \ = \ \oplus_{i=1}^{m} V_{i} \ \ , \ 
V_{i} \ = \ V^{i}/V^{i-1}$$
with Lie bracket $ [\hat{x}, \hat{y}]_{n} \ = \ \hat{[x,y]}$, where 
$\hat{x} = x + V^{i-1}$ if $x \in V^{i}\setminus 
V^{i-1}$. 

\begin{prop}
$\mathfrak{n}(\mathfrak{g},D)$  with distinguished space $V_{1} \ = \ D$ 
is the Lie algebra of a Carnot group of step $m$. It is called the nilpotentisation of the filtration \eqref{filtra}. 
\end{prop}
Remark that $\mathfrak{n}(\mathfrak{g},D)$ and $\mathfrak{g}$ have the same dimension, therefore they are isomorphic as vector spaces.

We set the distribution induced by $D$ to be 
$$D_{x} \ = \ T L_{x} D$$
where $x \in G$ is arbitrary and $L_{x}: G \rightarrow G$, $L_{x}y \ = \ xy$ 
is the left translation by $x$. We shall use the same notation $D$ for the 
induced distribution.

Let $\left\{ X_{1}, ... , X_{p} \right\}$ be a basis of the vector space $D$. 
We shall build a basis of $\mathfrak{g}$ which will give a vector spaces isomorphism with $\mathfrak{n}(\mathfrak{g},D)$. 

A word with letters $A \ = \ \left\{X_{1}, ... , X_{p}\right\}$ is a string 
$X_{h(1)}...X_{h(s)}$ where $h: \left\{1 , ... , s \right\} \rightarrow \left\{ 1 , ... , p \right\}$. The set of words forms the dictionary $Dict(A)$, ordered lexicographically. We set the function $Bracket: Dict(A) \rightarrow 
\mathfrak{g}$ to be 
$$Bracket(X_{h(1)}...X_{h(s)}) \ = \ [X_{h(1)},[X_{h(2)},[ ... , X_{h(s)}] ... ] $$
For any $x \in Bracket(Dict(A))$ let 
$\hat{x} \in Dict(A)$ be the least word  such that $Bracket(\hat{x}) \ = \ x$. 
The collection of all these words is denoted by $\hat{g}$. 

The length $l(x) \ = \ length(\hat{x})$ is well defined. The dictionary 
$Dict(A)$ admit a filtration made by the length of words function. In the same 
way the function $l$ gives the filtration 
$$V^{1} \cap Bracket(Dict(A)) \ \subset \ V^{2} \cap Bracket(Dict(A)) \ \subset \ ... \ Bracket(Dict(A))$$
Choose now, in the lexicographic order in $\hat{g}$, a set $\hat{B}$ such that 
$B \ = \ Bracket(\hat{B})$ is a  basis for $\mathfrak{g}$. Any  element $X$ in this basis can be written as 
$$X \ = \ [X_{h(1)},[X_{h(2)},[ ... , X_{h(s)}] ... ] $$ 
such that $l(X) \ = \ s$ or equivalently $X \in V^{s} \setminus V^{s-1}$. 

It is obvious that the map 
$$X \in B \ \mapsto \tilde{X} \ = \ X + V^{l(X)-1} \in \mathfrak{n}(\mathfrak{g},D)$$
is a bijection and that $\tilde{B} \ = \ \left\{ \tilde{X_{j}} \mbox{ : } 
j \ = \ 1, ... , dim \ \mathfrak{g} \right\}$ is a basis for 
$\mathfrak{n}(\mathfrak{g},D)$. We can identify then $\mathfrak{g}$ with 
$\mathfrak{n}(\mathfrak{g},D)$ by the identification $X_{j} \ = \ \tilde{X}_{j}$. 

Equivalently we can define the nilpotent Lie bracket $[ \cdot , \cdot ]_{n}$ 
directly on $\mathfrak{g}$, with the use of the dilatations on $\mathfrak{n}(\mathfrak{g},D)$. 

Instead of the filtration \eqref{filtra} let us start with a direct sum decomposition of $\mathfrak{g}$ 
\begin{equation}
\mathfrak{g} \ = \ \oplus_{i = 1}^{m} W_{i} \ , \ \ V^{i} \ = \ \oplus_{j=1}^{i} W_{j} 
\label{deco}
\end{equation}
such that $[V^{i}, V^{j}] \subset V^{i+j}$. The chain $V^{i}$ form a filtration 
like \eqref{filtra}. 

We set 
$$\delta_{\varepsilon} (x) \ = \ \sum_{i=1}^{m} \varepsilon^{i} x_{i}$$
for any $\varepsilon > 0$ and $x \in \mathfrak{g}$, which decomposes according 
to \eqref{deco} as $x \ = \ \sum_{i=1}^{m} x_{i}$. 

\begin{prop}
The limit 
\begin{equation}
[x,y]_{n} \ = \ \lim_{\varepsilon \rightarrow 0} \delta_{\varepsilon}^{-1} 
[\delta_{\varepsilon} x , \delta_{\varepsilon} y ]
\label{dnil}
\end{equation}
exists for any $x,y \in \mathfrak{g}$ and $(\mathfrak{g}, [ \cdot , \cdot]_{n})$ is the Lie algebra of a Carnot group with dilatations $\delta_{\varepsilon}$.
\label{limbra} 
\end{prop}

\paragraph{Proof.}
Let $x \ = \ \sum_{j=1}^{i} x_{j}$ and $y \ = \ \sum_{k=1}^{l} y_{k}$. Then 
$$[\delta_{\varepsilon} x , \delta-{\varepsilon} y ] \ = \ 
\sum_{s = 2}^{i+l} \varepsilon^{s} \sum_{j+k =s} \sum_{p=1}^{s} [x_{j},y_{k}]_{p}$$
We apply $\delta_{\varepsilon}^{-1}$ to this equality and we obtain: 
$$\delta_{\varepsilon}^{-1} 
[\delta_{\varepsilon} x , \delta-{\varepsilon} y ] \ = \ 
\sum_{s = 2}^{i+l}  \sum_{j+k =s} \sum_{p=1}^{s} \varepsilon^{s-p} [x_{j},y_{k}]_{p}$$
When $\varepsilon$ tends to 0 the expression converges to the  limit
\begin{equation}
[x,y]_{n} \ = \ \sum_{s = 2}^{i+l}  \sum_{j+k =s}  [x_{j},y_{k}]_{s}
\label{expre}
\end{equation}
For any $\varepsilon > 0$ the expression 
$$[x,y]^{\varepsilon} \ = \ \delta_{\varepsilon}^{-1} 
[\delta_{\varepsilon} x , \delta-{\varepsilon} y ]$$
is a Lie bracket (bilinear, antisymmetric and it satisfies the Jacobi identity). Therefore at the limit $\varepsilon \rightarrow 0$ we get a Lie bracket. Moreover, it is straightforward to see from the definition of 
$[x,y]_{n}$ that $\delta_{\varepsilon}$ is an algebra isomorphism. We conclude that $(\mathfrak{g}, [ \cdot , \cdot]_{n})$ is the Lie algebra of a Carnot group with dilatations $\delta_{\varepsilon}$.
\quad $\blacksquare$

\begin{prop}
Let $\left\{X_{1}, ... ,X_{dim \ \mathfrak{g}} \right\}$ be a basis 
of $\mathfrak{g}$ constructed from a basis $\left\{X_{1}, ... , X_{p}\right\}$ 
of $D$. Set 
$$W_{j} \ = \ span \ \left\{ X_{i} \mbox{ : } l(X_{i}) \ = \ j \right\}$$
Then the spaces $W_{j}$ provides a direct sum decomposition \eqref{deco}. Moreover the identification $\hat{X}_{i} = X_{i}$ gives a Lie algebra  isomorphism between 
$(\mathfrak{g}, [ \cdot , \cdot ]_{n})$ and $\mathfrak{n}(\mathfrak{g}, D)$. 
\end{prop}

\paragraph{Proof.}
In the basis $\left\{X_{1}, ... ,X_{dim \ \mathfrak{g}} \right\}$ the Lie bracket on $\mathfrak{g}$ looks like this: 
$$[X_{i}, X_{j}] \ = \ \sum C_{ijk}X_{k}$$
where $c_{ijk} = 0$ if $l(X_{i}) + l(X_{j}) < l(X_{k})$. From here the first part of the proposition is straightforward. The expression of the Lie bracket generated by the decomposition \eqref{deco} is obtained from  \eqref{expre}). We have 
$$[X_{i}, X_{j}]_{n} \ = \ \sum \lambda_{ijk} C_{ijk} X_{k}$$
where $\lambda_{ikj} = 1$ if $l(X_{i}) + l(X_{j}) = l(X_{k})$ and $0$ otherwise.  The Lie algebra isomorphism follows from the expression of the Lie bracket on $\mathfrak{n}(\mathfrak{g},D)$: 
$$[\hat{X}_{i}, \hat{X}_{j}] \ = \ [X_{i}, X_{j}] \ + \ V^{l(X_{i}) + l(X_{j})-1}$$
\quad $\blacksquare$

In conclusion  the expression of the nilpotent Lie  
bracket depends on the choice of basis $B$ trough the transport of the dilatations group from $\mathfrak{n}(\mathfrak{g},D)$ to $\mathfrak{g}$.

Let $N(G,D)$ be the simply connected Lie group with Lie algebra 
$\mathfrak{n}(\mathfrak{g},D)$. As previously, we identify $N(G,D)$ with 
$\mathfrak{n}(\mathfrak{g},D)$ by the exponential map. 

\subsection{Commutative smoothness for uniform groups}
This section is just a sketch and certainly needs to be improved. It is inspired from Bella\"{\i}che \cite{bell}, last section. I think that by improving the notions explained here and by making  axioms from the conclusions of lemma \ref{lgen} and Ball-Box theorem \ref{bbt}, then it is possible to avoid the  use of the notion of Lie group. 
For example, one could use consequences of proposition \ref{vit} in order 
to construct a general manifold structure, noncommutative derivative, moment map, as in sections \ref{manif}, \ref{nocog}, \ref{noco}. 

We start with the following setting: $G$ is a topological group endowed with an uniformity such that the operation is uniformly continuous.  More specifically, 
we introduce first the double of $G$, as the group $G^{(2)} \ = \ G \times G$ with operation 
$$(x,u) (y,v) \ = \ (xy, y^{-1}uyv)$$
The operation on the group $G$, seen as the function 
$$op: G^{(2)} \rightarrow G \ , \ \ op(x,y) \ = \ xy$$
is a group morphism. Also the inclusions: 
$$i': G \rightarrow G^{(2)} \ , \ \ i'(x) \ = \ (x,e) $$
$$i": G \rightarrow G^{(2)} \ , \ \ i"(x) \ = \ (x,x^{-1}) $$
are group morphisms. 

\begin{defi}
$G$ is an uniform group if we have two uniformity structures, on $G$ and 
$G^{2}$,  such that $op$, $i'$, $i"$ are uniformly continuous.   

A local action of a uniform group $G$ on a uniform  pointed space $(X, x_{0})$ is a function  
$\phi \in W \in \mathcal{V}(e)  \mapsto \hat{\phi}: U_{\phi} \in \mathcal{V}(x_{0}) \rightarrow 
V_{\phi}  \in \mathcal{V}(x_{0})$ such that the map $(\phi, x) \mapsto \hat{\phi}(x)$ is uniformly continuous from $G \times X$ (with product uniformity) 
to to $X$. Moreover, the action has the property:  for any $\phi, \psi \in G$ there is $D \in \mathcal{V}(x_{0})$ 
such that for any $x \in D$ $\hat{\phi \psi^{-1}}(x)$ and $\hat{\phi}(\hat{\psi}^{-1}(x))$ make sense and   $\hat{\phi \psi^{-1}}(x) = \hat{\phi}(\hat{\psi}^{-1}(x))$. 

Finally, a local group is an uniform space $G$ with an operation defined 
in a neighbourhood of $(e,e) \subset G \times G$ which satisfies the uniform group axioms locally.  
\end{defi}
Remark that a local group acts locally at left (and also by conjugation) on itself. 

\begin{defi}
A conical local uniform group $N$ is a local group with a local action of 
$(0,+\infty)$ by morphisms $\delta_{\varepsilon}$ such that 
$\lim_{\varepsilon \rightarrow 0} \delta_{\varepsilon} x \ = \ e$ for any 
$x$ in a neighbourhood of the neutral element $e$. 
\end{defi}

We shall make the following  hypotheses on the local uniform group $G$: there is a local action of $(0, +\infty)$ (denoted by 
$\delta$), on $(G, e)$ such that 
\begin{enumerate}
\item[H0.] $\delta_{\varepsilon}(e) \ = \ e$ and $\lim_{\varepsilon \rightarrow 0} \delta_{\varepsilon} x \ = \ e$ uniformly with respect to $x$. 
\item[H1.] the limit 
$$\beta(x,y) \ = \ \lim_{\varepsilon \rightarrow 0} \delta_{\varepsilon}^{-1} 
\left((\delta_{\varepsilon}x) (\delta_{\varepsilon}y ) \right)$$
is well defined in a neighbourhood of $e$ and the limit is uniform. 
\item[H2.] the limit 
$$\hat{x}^{-1} \ = \ \lim_{\varepsilon \rightarrow 0} \delta_{\varepsilon}^{-1} 
\left( ( \delta_{\varepsilon}x)^{-1}\right)$$
is well defined in a neighbourhood of $e$ and the limit is uniform. 
\end{enumerate}

\begin{prop}
Under the hypotheses H0, H1, H2 $(G,\beta)$ is a conical local uniform group. 
\end{prop}

\paragraph{Proof.}
All the uniformity assumptions permit to change at will the order of taking 
limits. We shall not insist on this further and we shall concentrate on the 
algebraic aspects.
 
We have to prove the associativity, existence of neutral element, existence of inverse and the property of being conical. The proof is straightforward. 
For the associativity $\beta(x,\beta(y,z)) \ = \ \beta(\beta(x,y),z)$ we compute: 
$$\beta(x,\beta(y,z)) \ = \ \lim_{\varepsilon \rightarrow 0 , \eta \rightarrow 0} \delta_{\varepsilon}^{-1} \left\{ (\delta_{\varepsilon}x) \delta_{\varepsilon/\eta}\left( (\delta_{\eta}y) (\delta_{\eta} z) \right) \right\}$$
We take $\varepsilon = \eta$ and we get 
$$ = \beta(x,\beta(y,z)) \ = \ \lim_{\varepsilon \rightarrow 0}\left\{ 
(\delta_{\varepsilon}x) (\delta_{\varepsilon} y) (\delta_{\varepsilon} z) \right\}$$
In the same way: 
$$\beta(\beta(x,y),z) \ = \ \lim_{\varepsilon \rightarrow 0 , \eta \rightarrow 0} \delta_{\varepsilon}^{-1} \left\{ (\delta_{\varepsilon/\eta}x)\left( (\delta_{\eta}x) (\delta_{\eta} y) \right) (\delta_{\varepsilon} z) \right\}$$
and again taking $\varepsilon = \eta$ we obtain
$$\beta(\beta(x,y),z) \ = \  \lim_{\varepsilon \rightarrow 0}\left\{ 
(\delta_{\varepsilon}x) (\delta_{\varepsilon} y) (\delta_{\varepsilon} z) \right\}$$
The neutral element is $e$, from H0 (first part): $\beta(x,e) \ = \beta(e,x) \ = \ x$. The inverse of $x$ is $\hat{x}^{-1}$, by a similar argument: 
$$\beta(x, \hat{x}^{-1})  \ = \ \lim_{\varepsilon \rightarrow 0 , \eta \rightarrow 0} \delta_{\varepsilon}^{-1} \left\{ (\delta_{\varepsilon}x) 
\left( \delta_{\varepsilon/\eta}(\delta_{\eta}x)^{-1}\right) \right\}$$
and taking $\varepsilon = \eta$ we obtain 
$$\beta(x, \hat{x}^{-1})  \ = \ \lim_{\varepsilon \rightarrow 0} 
\delta_{\varepsilon}^{-1} \left( (\delta_{\varepsilon}x) (\delta_{\varepsilon}x)^{-1}\right) \ = \ \lim_{\varepsilon \rightarrow 0} \delta_{\varepsilon}^{-1}(e) \ = \ e$$
Finally, $\beta$ has the property: 
$$\beta(\delta_{\eta} x, \delta_{\eta}y) \ = \ \delta_{\eta} \beta(x,y)$$
which comes from the definition of $\beta$ and commutativity of multiplication 
in $(0,+\infty)$. This proves that $(G,\beta)$ is conical.  
\quad $\blacksquare$

We arrive at a natural realization of the tangent space to the neutral element. 
Let us denote by $[f,g] \ = \ f \circ g \circ f^{-1} \circ g^{-1}$ the commutator of two transformations. For the group we shall denote by 
$L_{x}^{G} y \ = \ xy$ the left translation and by $L^{N}_{x}y \ = \ \beta(x,y)$. The preceding proposition tells us that $(G,\beta)$ acts locally by left 
translations on $G$. 

\begin{prop}
We have the equality
\begin{equation}
\lim_{ \lambda \leq \varepsilon \rightarrow 0} \delta_{\lambda}^{-1} [\delta_{\varepsilon}, L_{(\delta_{\lambda}x)^{-1}}^{G}] \delta_{\lambda} \ = \ L_{x}^{N}
\end{equation}
and the limit is uniform with respect to $x$. 
\label{vit}
\end{prop}

\paragraph{Proof.}
We shall use the equality:
\begin{equation} 
\lim_{\varepsilon \rightarrow 0} \delta_{\varepsilon}^{-1} \left( 
\beta((\delta_{\varepsilon}x),\delta_{\varepsilon}y))\right) \ = \ 
\lim_{\varepsilon \rightarrow 0} \delta_{\varepsilon} \left( \beta(\hat{x}^{-1}, \delta_{\varepsilon}^{-1}(\beta(x,y))) \right)
\label{tem1}
\end{equation}
which is a consequence of H0 and the previous proposition. 
The left hand side of this equality equals $\beta(x,y)$. We shall look now at  the right hand side (RHS) of the equality. The uniformity assumptions and some computations lead us to the following equality, under the constraint 
$0< \max \left\{ \mu, \lambda,\eta \right\} < \varepsilon$: 
$$RHS \ = \ \lim_{(\varepsilon, \mu, \lambda, \eta) \rightarrow 0} 
\delta_{\varepsilon/\mu} \left\{ \delta_{\mu/\lambda}\left( (\delta_{\lambda}x)^{-1}\right) \delta_{\mu/(\varepsilon\eta)} \left( (\delta_{\eta}x) (\delta_{\eta}y)\right) \right\}$$
The mentioned constraint forces all expressions to make sense. 
Take now $\mu = \lambda = \eta < \varepsilon$ and get
$$RHS \ = \ \lim_{\varepsilon, \lambda \rightarrow 0} 
\delta_{\lambda}^{-1} \delta_{\varepsilon} \left( \left( (\delta_{\lambda}x)^{-1}\right) \delta_{\varepsilon}^{-1} \left( (\delta_{\lambda}x) (\delta_{\lambda}y)\right) \right)$$
as desired.
\quad $\blacksquare$

An easy corollary is that 
$$\lim_{\lambda \rightarrow 0} [L_{(\delta_{\lambda}x)^{-1}}^{G}, \delta_{\lambda}^{-1}] \ = \ L^{N}_{x}$$

\begin{defi}
The group $VT_{e}G$ formed by all transformations $L_{x}^{N}$ is called the virtual tangent space at $e$ to $G$. 
\end{defi}

The virtual tangent space $VT_{x}G$ at $x \in G$ to $G$ is obtained by translating the group operation and the dilatations from $e$ to $x$. This means: define a new operation on $G$ by 
$$y \stackrel{x}{\cdot} z \ = \ y x^{-1}z$$
The group $G$ with this operation is isomorphic to $G$ with old operation and 
the left translation $L^{G}_{x}y \ = \ xy$ is the isomorphism. The neutral element is $x$. 
Introduce also the dilatations based at $x$ by 
$$\delta_{\varepsilon}^{x} y \ = \ x \delta_{\varepsilon}(x^{-1}y)$$
Then $G^{x} \ = \ (G,\stackrel{x}{\cdot})$ with the group of dilatations $\delta_{\varepsilon}^{x}$ satisfy the axioms Ho, H1, H2. Define then the virtual tangent 
space $VT_{x}G$ to be: $VT_{x}G \ = \ VT_{x} G^{x}$. A short computation using 
proposition \ref{vit} shows that 
$$VT_{x} G \ = \ \left\{ L^{N,x}_{y} \ = \ L_{x} L^{N}_{x^{-1}y} L_{x} \mbox{ : } y \in U_{x} \in \mathcal{V}(X) \right\}$$
where 
$$L^{N,x}_{y} \ = \ \lim_{ \lambda \leq \varepsilon \rightarrow 0} \delta_{\lambda}^{-1,x} [\delta_{\varepsilon}^{x}, L_{(\delta_{\lambda}x)^{x, -1}}^{G}] \delta_{\lambda}^{x}$$

We shall introduce the notion of commutative smoothness, which contains  
a derivative resembling with Pansu derivative. 
Independently, the author introduced a general "topological" derivative 
in \cite{buli} (there are a lot of typographical errors in this reference, but the ideas behind are quite nice in my opinion). Lack of knowledge stopped the development of this notion until recently, when I have learned about Carnot groups and Pansu derivative. 
As it will be seen further, interesting information will be found when 
noncommutative smoothness is considered. 

\begin{defi}
A function $f: G_{1} \rightarrow G_{2}$ is commutative smooth, where 
$G_{1}, G_{2}$ are two groups satisfying H0, H1, H2,  if the application 
$$(x,u) \in G_{1}^{(2)} \ \mapsto \ (f(x), Df(x)u) \in G_{2}^{(2)}$$
exists and it is  continuous, where 
$$Df(x)u \ = \ \lim_{\varepsilon \rightarrow 0} \delta_{\varepsilon}^{-1} 
\left(f(x)^{-1}f(x \delta_{\varepsilon}u)\right)$$
and the convergence is uniform with respect to $(x,u)$. 
\label{fdcd}
\end{defi}

For example the left translations $L_{x}$ are commutative smooth and 
the derivative equals identity. If we want to see how the derivative moves 
the virtual tangent spaces we have to give a definition. 

\begin{defi}
Le $f: G \rightarrow G$ be a commutative smooth function. The virtual tangent to $f$ is defined by: 
$$VT f(x): VT_{x}G \rightarrow VT_{f(x)} G \ , \ \ VT f(x) L^{N,x}_{y} \ = \ 
L_{f(x)} L^{N}_{Df(x)y} L^{-1}_{f(x)}$$
\label{vcd}
\end{defi}

With this definition $L_{x}$ is commutative smooth and it's virtual tangent in any point $y$ is a 
group morphism from $VT_{y}G$ to $VT_{xy}G$. More generally, reasoning as in proposition \ref{ppansuprep}, we get: 

\begin{prop}
The virtual tangent is morphism of conical groups. 
\end{prop}

Nevertheless the right translations are not commutative smooth. This failure pushes us to consider a notion of 
noncommutative derivative. We shall meet the same failure soon, when looking to the manifold structure, in the sense of definition \ref{dcum}, of the group 
$G$.

Now that we have a model for the tangent space to $e$ at $G$, we can show that 
the operation is commutative smooth. 

\begin{prop}
Let $G$ satisfy H0, H1, H2 and $\delta_{\varepsilon}^{(2)} : G^{(2)} \rightarrow G^{(2)}$ be defined by 
$$\delta_{\varepsilon}^{(2)} (x,u) \ = \ (\delta_{\varepsilon}x, 
\delta_{\varepsilon} y)$$
Then $G^{(2)}$ satisfies H0, H1, H2, $op$ is commutative smooth  and we have the relation: 
$$D \ op \ (x,u) (y,v) \ = \ \beta(y,v)$$
\label{opsm}
\end{prop}

\paragraph{Proof.}
It is sufficient to use the morphism property of the operation. Indeed, the right hand side of the relation to be proven is
$$RHS \ = \ \lim_{\varepsilon \rightarrow 0} 
\delta_{\varepsilon}^{-1} \left( op(x,u)^{-1} op(x,u) op \left(\delta_{\varepsilon}^{(2)}(y,v)\right)\right) \ = $$
$$=  \ \lim_{\varepsilon \rightarrow 0} 
\delta_{\varepsilon}^{-1} \left( op(\delta_{\varepsilon}^{(2)}(y,v))\right) \ = \ \beta(y,v)$$
The rest is trivial. 
\quad $\blacksquare$

\subsection{Manifold structure}
\label{manif}

The notion of virtual tangent space is not based on the use 
of distances, but on the use of dilatations. In fact, any manifold has a tangent space to any of its points, not only the Riemannian manifolds. We shall prove in this section that $VT_{e} G$ is isomorphic to the nilpotentisation 
$N(G,D)$. 

Nevertheless  
$G$ does not have the structure of a $N(G,D)$ $\mathcal{C}^{1}$  manifold, in the sense of definition \ref{dcum}.

We start from Euclidean norm on $D$ and we choose an orthonormal basis of $D$. We can then extend the Euclidean norm to $\mathfrak{g}$ by 
stating that the basis of $\mathfrak{g}$ constructed, as explained, from the basis on $D$, is orthonormal. By left translating the Euclidean norm on 
$\mathfrak{g}$ we endow $G$ with a structure of Riemannian manifold. The 
induced Riemannian  distance $d_{R}$ will give an uniform structure on $G$. 
This distance is left invariant: 
$$d_{R}(xy, xz) \ = \ d_{R}(y,z)$$
for any $x,y,z \in G$. 

Any left invariant distance $d$ is uniquely determined if we set $d(x) \ = \ d(e,x)$.

The following lemma is important (compare with lemma \ref{p2.4}). 

\begin{lema}
Let $X_{1}, ... , X_{p}$ be a basis of $D$. 
Then there are $U \subset G$ and $V \subset N(G,D)$, open neighbourhoods of the neutral elements $e_{G}$, $e_{N}$ respectively,  and a surjective function $g: \left\{1, ... , M \right\} \rightarrow 
\left\{ 1, ... ,p\right\}$ such that any $x \in U$, $y \in V$ can be written as 
\begin{equation}
x \ = \ \prod_{i = 1}^{M} \exp_{G}(t_{i}X_{g(i)}) \ \ , \ y \ = \ \prod_{i = 1}^{M} \exp_{N}(\tau_{i}X_{g(i)})
\label{fpgen}
\end{equation}
\label{lgen}
\end{lema}

\paragraph{Proof.}
We shall make the proof for $G$; the proof for $N(G,D)$ will follow from the identifications explained before. 

Denote by $n$ the dimension of $\mathfrak{g}$. We start the proof with the 
remark that the function 
\begin{equation}
(t_{1}, ... , t_{n})  \ \mapsto  \ \prod_{i = 1}^{n} \exp_{G}(t_{i}X_{i})
\label{tp1}
\end{equation}
is invertible in a neighbourhood of $0 \in R^{n}$, where the $X_{i}$ are elements of a basis $B$ constructed as before. Remember that each $X_{i} \in B$
is a multi-bracket of elements from the basis of $D$. If we replace a bracket  
$\exp_{G} (t [x,y])$ in the expression \eqref{tp1} by $exp_{G} (t_{1}x) \ exp_{G} (t_{2}y) \ exp_{G} (t_{3}x)
 \ exp_{G} (t_{4}y)$ and we replace $t$ by $(t_{1}, ... , t_{4})$ then  the 
image of a neighbourhood of $0$ by the obtained function still covers a 
neighbourhood of the neutral element. We repeat this procedure a finite number of times and the thesis is proven. 
\quad $\blacksquare$

As a corollary we obtain the Chow theorem for our particular example. 

\begin{thm}
Any two points $x,y \in G$ can be joined by a horizontal curve. 
\end{thm}

Let $d_{G}$ be the Carnot-Carath\'eodory distance induced by the distribution 
$D$ and the metric. This distance is also left invariant. We obviously have $d_{R} \ \leq \ d_{G}$. We want to show that $d_{R}$ and $d_{G}$ induce the same uniformity on $G$. 

Let us introduce another left invariant distance on $G$
$$d^{1}_{G}(x) \ = \ \inf \left\{ \sum \mid t_{i} \mid \mbox{ : } x \ = \ 
\prod \exp_{G} (t_{i} Y_{i}) \ , \ Y_{i} \in D \right\}$$
and the auxiliary functions : 
$$\Delta^{1}_{G}(x) \ = \ \inf \left\{ \sum_{i=1}^{M} \mid t_{i} \mid \mbox{ : } x \ = \ 
\prod \exp_{G} (t_{i} X_{g(i)}) \right\}$$
$$\Delta^{\infty}_{G}(x) \ = \ \inf \left\{ \max \mid t_{i} \mid \mbox{ : } x \ = \ 
\prod \exp_{G} (t_{i} X_{g(i)}) \right\}$$
From theorems \ref{t441amb} and \ref{t411amb} we see that $d^{1}_{G} \ = \  d_{G}$. Indeed, it is straightforward that  $d_{G} \leq d_{G}^{1}$.  On the other part 
$d_{G}^{1}(x)$ is less equal than the variation of any Lipschitz curve joining 
$e_{G}$ with $x$. Therefore we have equality. 

The functions $\Delta^{1}_{G}$, $\Delta^{\infty}_{G}$ don't induce left invariant distances. Nevertheless they are useful, because of their equivalence: 
\begin{equation}
\Delta^{\infty}_{G}(x) \ \leq \ \Delta^{1}_{G}(x) \ \leq \ M \ \Delta_{G}^{\infty}(x)
\label{dequi}
\end{equation}
for any $x \in G$. This is a consequence of the lemma 
\ref{lgen}. 

We have therefore the chain of inequalities: 
$$d_{R} \ \leq \ d_{G} \ \leq \ \Delta^{1}_{G} \ \leq \ M \Delta^{\infty}_{G}$$
But from the proof of lemma \ref{lgen} we see that $\Delta^{\infty}_{G}$ is uniformly continuous. This proves the equivalence of the uniformities. 

Because $\exp_{G}$ does not deform much the $d_{R}$ distances near $e$, we see that the group $G$ with the dilatations
$$\tilde{\delta}_{\varepsilon} (\exp_{G} x) \ = \ \exp_{G}(\delta_{\varepsilon} 
x)$$
satisfies H0, H1, H2. 

The same conclusion is true for the local uniform group (with the uniformity induced by the Euclidean distance) $\mathfrak{g}$ with the operation: 
$$XY \ = \ \log_{G} \left( \exp_{G} (X) \exp_{G} (Y) \right)$$
for any $X,Y$ in a neighbourhood of $0 \in \mathfrak{G}$. Here the dilatations 
are $\delta_{\varepsilon}$. We shall denote this group by $\log G$These two groups are isomorphic as local uniform groups by the map $\exp_{G}$. Dilatations commute with the isomorphism. They have therefore isomorphic (by $\exp_{G}$) virtual tangent spaces. 

\begin{thm}
The virtual tangent space $VT_{e}G$ is isomorphic to $N(G,D)$. More precisely 
$N(G,D)$ is equal (as local group) to the virtual tangent space to $\log G$: 
$$N(G,D) \ = \ VT_{0} \log G$$
\label{teore}
\end{thm}

\paragraph{Proof.}
The product $XY$ in $\log G$ is given by  Baker-Campbell-Hausdorff formula 
$$X \opg Y \ = \ X + Y + \frac{1}{2} [X,Y] + ... $$
Use proposition \ref{limbra} to compute $\beta(X,Y)$ and show that 
$\beta(X,Y)$ equals the nilpotent multiplication. 
\quad $\blacksquare$

Let $\left\{ X_{1}, ... , X_{dim \ \mathfrak{g}} \right\}$ be a basis of 
$\mathfrak{g}$ constructed from a basis $\left\{ X_{1}, ... , X_{p} \right\}$, 
as explained before. Denote by 
$$b: \mathfrak{g} \rightarrow \mathfrak{n}(\mathfrak{g},D)$$
the identification of the vectorspace $\mathfrak{g}$ with $\mathfrak{n}(\mathfrak{g},D)$ using the basis $\left\{ X_{1}, ... , X_{dim \ \mathfrak{g}} \right\}$. The exponential and the logarithm of the groups 
$G$ and $N(G,D)$ will be denoted by $\exp_{G}$, $\log_{G}$ and 
$\exp_{N}$, $\log_{N}$ respectively.

We shall consider the following  $N(G,D)$ atlas for $G$: 
$$\mathcal{A} \ = \ \left\{ \psi_{x}: U_{x} \subset G \rightarrow N(D,G) 
\mbox{ : } x \in G  \ , \ \right. $$ 
$$\left.  \psi_{x}(y) \ = \  \exp_{N}\circ b \circ \log_{G} (x^{-1}y)  \right\}$$
We shall study the differentiability properties of the transition functions. 
We want to  work in $\mathfrak{g}$. For this we have to transport (in a neighbourhood $0 \in \mathfrak{g}$) the interesting operations, namely: 
$$X \opg Y \ = \ \log_{G} \left( \exp_{G}(X) \exp_{G}(Y) \right)$$
$$X \opn Y \ = \ \log_{N} \left( \exp_{N}(X) \exp_{N}(Y) \right)$$
and to use instead of the chart $\psi_{x}: U_{x} \subset G \rightarrow N(D,G)$,  
the chart $\phi_{X} : U_{x} \subset G \rightarrow \mathfrak{g}$, where 
$x \ = \ \exp_{G} X$ and 
$\phi_{X}(y) \ = \ (-X) \opg \log_{G}(y)$. 

The transition function from $\phi_{X}$ to $\phi_{0}$ is then the left 
translation by $X$, with respect to the operation $\opg$. We denote this translation by $L^{\mathfrak{g}}_{X}$. We want to know if this function is 
Pansu derivable from $(\mathfrak{g}, \opn)$ to itself. See the difference: 
the function is certainly derivable (or commutative smooth) in the sense of definition \ref{fdcd}, for the operation $\opg$, but what about the derivability with respect to the operation $\opn$? In order to answer we shall use the following trick: a $C^{\infty}$ function $f: N \rightarrow N$ is Pansu 
derivable if and only if it  preserves the horizontal distribution and the derivative in each point is a morphism from $N$ to $N$. 

The (classical) derivative of $L^{\mathfrak{g}}_{X}$ moves the distribution 
$D^{n}(Y) \ = \ D L^{\mathfrak{n}}_{Y}(e)$
into $D L^{\mathfrak{g}}_{X} D^{n}(Y) \ \subset \ T_{X \opg Y} \mathfrak{g}$. 
The horizontal distribution in $X \opg Y$, corresponding to the group operation $\opn$, is $D^{n}(X \opg Y)$. The difference between these two distributions 
is measured by one of the linear transformations:
$$A_{X,Y} \ : T_{x \opg Y} \mathfrak{g} \rightarrow T_{x \opg Y} \mathfrak{g}$$ 
\begin{equation} 
A_{X,Y} \ = \ D L^{\mathfrak{g}}_{X \opg Y} (0) \left( D L^{\mathfrak{g}}_{Y}(0) \right)^{-1} D L^{\mathfrak{n}}_{Y}(0) \left( D L^{\mathfrak{n}}_{X \opg Y} (0) \right)^{-1}
\label{axyj}
\end{equation}
$$A^{X,Y} \ : \mathfrak{g} \rightarrow \mathfrak{g}$$ 
\begin{equation} 
A^{X,Y} \ = \ \left( D L^{\mathfrak{n}}_{X \opg Y} (0) \right)^{-1} 
D L^{\mathfrak{g}}_{X \opg Y} (0) \left( D L^{\mathfrak{g}}_{Y}(0) \right)^{-1} D L^{\mathfrak{n}}_{Y}(0) 
\label{axys}
\end{equation}

Let then $\mathfrak{J}(G,D)$ be the Lie group generated by these transformations
$$\mathfrak{J}(G,D) \ = \ \langle \left\{ A^{X,Y} \mbox{ : } X,Y \in \mathfrak{g} \right\} \rangle$$
It is then straightforward to see that the algebra $\mathfrak{j}(G,D)$ of this group contains the algebra generated by all the linear transformation with the form 
\begin{equation}
a_{x} \ = \ ad^{G}_{x} \ - \  ad^{N}_{x}
\label{ax}
\end{equation}
The necessary and sufficient condition for the group $\mathcal{J}(G,D)$ to be included in the  group $End((\mathfrak{g}, \opn), D)$ is that all the elements $a_{X}$,  to be in the algebra of the mentioned linear group. But this is equivalent with one of the (equivalent) conditions: 
\begin{enumerate}
\item[(i)] $Ad^{G}_{x}$ is $[ \cdot , \cdot ]_{N}$ morphism, for any $x$ in a 
neighbourhood of the identity, 
\item[(ii)] for any $X,U,V \in \mathfrak{g}$ we have the identity: 
\begin{equation}
[[X,U]_{G},V]_{N} \ + \ [U,[X,V]_{G}]_{N} \ = \ [X,[U,V]_{N}]_{G}
\label{magic}
\end{equation}
\end{enumerate}
We collect what we have found in the following theorem. 

\begin{thm}
Set $Ad^{G}$ to be the adjoint representation of $G$ and $Ad^{N}$ the adjoint representation of$N(G,D)$, seen as group of linear transformations on $\mathfrak{g}$ (via the identification function $b$). The following are equivalent: 
\begin{enumerate}
\item[(a)] $\mathcal{J}(G,D) \ \subset \ End(\mathfrak{g},\opn)$, 
\item[(b)] $Ad^{G} \subset \ End(\mathfrak{g},\opn)$, 
\item[(c)] the relation \eqref{magic} is true. 
\end{enumerate}
If \eqref{magic} holds then $Ad^{G} Ad^{N}$ is a group and its Lie algebra 
is the adjoint representation of the algebra $\mathfrak{g} \oplus \mathfrak{g}$ 
with the bracket: 
\begin{equation}
[(X,U),(Y,V)] \ = \ ([X,Y]_{G}, [X,V]_{G} + [U,Y]_{G} + [V,U]_{N})
\label{magicbra}
\end{equation} 
\label{jstruct}
\end{thm}

\paragraph{Proof.}
The equivalence of (a), (b), (c) has been explained. Suppose now that \eqref{magic} is true. It is then straightforward to show that the space of all 
elements 
$$a_{(X,Y)} \ = \ ad^{G}_{X} \ - \ ad^{N}_{Y}$$
forms a Lie algebra with the linear commutator as bracket. Moreover, we have: 
$$[a_{(X,U)},a_{(Y,V)}] \ = \ a_{[(X,U),(Y,V)]}$$
This shows that $Ad^{G} Ad^{N}$ is a group and the property of its algebra. 
In order to finish the proof remark that $a_{(X,X)} \ = \ a_{X}$, defined at 
\eqref{ax}. 
\quad $\blacksquare$

\begin{cor}
The atlas $\mathcal{A}$
gives to $G$ a $\mathcal{C}^{1}$ $N(D,G)$ manifold structure if and only if 
\eqref{magic} is true and for any $\varepsilon>0$, $X,Y \in \mathfrak{g}$ 
\begin{equation}
\delta_{\varepsilon}^{-1} [ X , \delta_{\varepsilon}Y ]_{N} 
\ - \  \delta_{\varepsilon}^{-1} [ X , \delta_{\varepsilon}Y ]_{G} \ = \ 
[X,Y]_{N} \ - \ [X,Y]_{G}
\label{cucu}
\end{equation}
\label{nonatlas}
\end{cor}

\paragraph{Proof.}
The atlas $\mathcal{A}$
gives to $G$ a $\mathcal{C}^{1}$ $N(D,G)$ manifold structure if and only if 
$\mathcal{J}(G,D) \ \subset \ HL(\mathfrak{g},\opn)$ and any of its elements 
commute with dilatations $\delta_{\varepsilon}$. This is equivalent with 
\eqref{magic} and \eqref{cucu}. 
\quad $\blacksquare$

\begin{rk}
A Lie group is a manifold endowed with a smooth operation. In what sense is then $G$ a (sub-Riemannian) Lie group? We already have problems to assign an 
atlas with smooth transition functions to $G$. The real meaning of the 
corollary \ref{nonatlas} is that even if we succeed to give to $G$ an atlas 
with smooth transition functions then the left translations and the exponential map will not be smooth. This is a problem. If we renounce to give $G$ a manifold structure, at least the operation is smooth 
in the sense of proposition \ref{opsm}. 

Returning to the problem of the atlas, if $N(G,D)$ is a Heisenberg group, or 
equivalently the distribution has codimension one, then the commentary of 
Bella\"{\i}che \cite{bell}, page 73, shows that such an atlas exists, as a consequence of Darboux theorem on normal forms for contact differential forms. 
But such an atlas is not compatible with the operation. 
\end{rk}

We shall see now how it can be possible that the operation is commutative smooth but 
the right translations are not. 

We shall  look to the double $G^{(2)}$. This group is isomorphic 
with $G^{2} \ = \ G \times G$ with componentwise multiplication by the isomorphism
$$F: G^{(2)} \rightarrow G^{2} \ , \ \ F(x,u) \ = \ (x,xu)$$

 The Lie algebra of $G^{(2)}$ is easily then 
$\mathfrak{g}^{(2)} \ = \ \mathfrak{g} \times \mathfrak{g}$, with 
Lie bracket 
$$[(x,u), (y,v)]_{(2)} \ = \ \left( [x,u], [x+u,y+v] - [x,y]\right)$$
Take in $\mathfrak{g}^{(2)}$ the distribution $D^{(2)} \ = \ D \times D$. This distribution generates $\mathfrak{g}^{(2)}$, and the dilatations are 
$$\delta_{\varepsilon}^{(2)} (x,u) \ = \ (\delta_{\varepsilon}x, \delta_{\varepsilon}u)$$
  We  translate the distribution  at left all over $G^{(2)}$. 

We know that the operation is commutative smooth. We shall check if  the  classical derivative of the operation  transports the distribution 
$D^{(2)}$ in the distribution $D$. Straightforward computation shows that  
$$D \ op (e,e) D^{(2)} \ = \ D$$
where $D$ here means classical derivative. Indeed, $D \ op (e,e)(X,Y) \ = \ 
X+Y$ and $D$ is a vector space, therefore $(X,Y) \in D^{(2)}$ implies 
$X+Y \in D$.  

For general $(x,u) \in G^{(2)}$ this is no longer true. Indeed, we have: 
$$D \ op (x,u) (Y,V) \ = \ DL^{G}_{xu}(0) \left\{ Y + V + (Ad^{G}_{xu} - Ad^{G}_{u}) Y \right\}$$
therefore a sufficient condition for $D \ op (x,u) D^{(2)} \ = \ D$ is that 
$x \ = \ e$. 

In conclusion the fact that the operation is commutative smooth does not imply that so are the right translations. 

We can hope that by a slight modification in the definition of smoothness we shall be able to give to $G$ a manifold structure and simultaneously to 
have smooth right translations. For the manifold structure we have to "smooth" the group generated by the elements with the form 
$$\left(Ad^{G}_{X}\right)^{-1} Ad^{N}_{X}$$ 
To solve the problem of the right translations we need that the group $Ad^{G}$ 
be "smoothed".

Before doing this  we shall see that  the pure  metric point of view  does not feel these problems but in such a  precise way. 

\subsection{Metric tangent cone}
\label{metriccone}

We shall prove the result of Mitchell \cite{mit} theorem 1, that $G$ admits in any point a metric tangent cone, which is isomorphic with the nilpotentisation 
$N(G,D)$. Mitchell theorem is true for regular sub-Riemannian manifolds. 
The proof that we give here is based on the lemma \ref{lgen} and Gromov 
\cite{gromov} section 1.2. 

Because left translations are isometries, it is sufficient to prove that $G$ 
admits a metric tangent cone in identity and that the tangent cone is isometric 
with $N(G,D)$. 

For this we transport all in the algebra $\mathfrak{g}$, endowed with two 
brackets $[\cdot, \cdot]_{G}$, $[\cdot, \cdot ]_{N}$ and with two operations 
$\opg$ and $\opn$. Denote by $d_{G}$, $d_{N}$ the Carnot-Carath\'eodory distances corresponding to the $\opg$, respectively $\opn$ left invariant distributions on (a neighbourhood of $0$ in) $\mathfrak{g}$. $l_{G}$, 
$l_{N}$ are the corresponding length functionals. We shall denote by $B_{G}(x,R)$, $B_{N}(x,R)$ the balls centered in $x$ with 
radius $R$ with respect to $d_{G}$, $d_{N}$. 

We can refine lemma \ref{lgen} in order to obtain the Ball-Box theorem in this 
more general situation. 

\begin{thm}
(Ball-Box Theorem) Denote by 
$$Box^{1}_{G} (\varepsilon) \ = \ \left\{ x \in G \mbox{ : } \Delta^{1}_{G}(x) \ < \ \varepsilon \right\} $$
$$Box^{1}_{N} (\varepsilon) \ = \ \left\{ x \in N \mbox{ : } \Delta^{1}_{N}(x) \ < \ \varepsilon \right\} $$
For small $\varepsilon > 0$ there is a constant $C>1$ such that 
$$\exp_{G} \left( Box^{1}_{N} (\varepsilon) \right) \ \subset \ 
Box^{1}_{G} (C \varepsilon) \ \subset \
\exp_{G} \left( Box^{1}_{N} (C^{2} \varepsilon) \right)$$
\label{bbt}
\end{thm}

\paragraph{Proof.}
Reconsider the proof of lemma \ref{lgen}. This time, instead of the trick of 
replacing commutators $\exp_{G}[X,Y](t)$ with four letters words 
$$\exp_{G}X(t_{1}) \exp_{G} Y (t_{2}) \exp_{G}X (t_{3}) \exp_{G} Y (t_{4})$$ 
we shall use a smarter replacement (which, important! , works equally for 
the nilpotentisation $N(G,D)$). 
 
We start with a basis $\left\{ Y_{1}, ... , Y_{n} \right\}$ of the algebra 
$\mathfrak{g}$, constructed from multibrackets of elements $\left\{ X_{1}, ... , X_{p} \right\}$ which form a basis for the distribution $D$. Introduce an Euclidean distance by declaring the basis of $\mathfrak{g}$ orthonormal. Set 
$$[X(t), Y(t)]^{\circ} \ = \ (tX) \opg (tY) \opg (-tX) \opg (-tY)$$
It is known that for any set of vectors $Z_{1}, ... , Z_{q} \in \mathfrak{g}$, 
if we denote by $\alpha^{\circ}(Z_{1}(t), ... , Z_{q}(t))$ a $[ \cdot , \cdot]^{\circ}$ multibracket and by $\alpha(Z_{1}, ... , Z_{q})$ the same constructed 
$[\cdot , \cdot ]$ multibracket, then we have 
$$ \alpha^{\circ}(Z_{1}(t), ... , Z_{q}(t)) \ = \ t^{q} \alpha(Z_{1}, ... , Z_{q}) \ + \ o(t^{q})$$
with respect to the Euclidean norm.

Remark as previously that   the function 
\begin{equation}
(t_{1}, ... , t_{n})  \ \mapsto  \ \prod_{i = 1}^{n} (t_{i}Y_{i})
\label{tpu1}
\end{equation}
is invertible in a neighbourhood of $0 \in R^{n}$. Each $X_{i}$ from the 
basis of $\mathfrak{g}$ can be written as a multibracket
$$X_{i} \ = \ \alpha_{i}( X_{i_{1}}, ... , X_{j^{i}})$$ 
which has the length $l_{i} = j^{i} - 1$. If $l_{i}$ is odd then replace 
$(t_{i}X_{i})$ by 
$$\alpha^{\circ}( X_{i_{1}}(t^{1/l_{i}}), ... , X_{j^{i}}(t^{1/l_{i}}))$$
If $l_{i}$ is even then the multibracket $\alpha_{i}$ can be rewritten as 
$$\alpha_{i}( X_{i_{1}}, ... , X_{j^{i}}) \ = \ \beta_{i}(X_{i_{1}} , ... ,  
[X_{i_{k}}, X_{i_{k+1}}], ... , X_{j^{i}})$$
Replace then $(t_{i}X_{i})$ by 
$$\left\{ \begin{array}{ll}
\beta_{i}^{\circ}(X_{i_{1}}(\mid t \mid^{1/l_{i}}) , ... ,  
[X_{i_{k}}(\mid t \mid^{1/l_{i}}), X_{i_{k+1}}(\mid t \mid^{1/l_{i}})]^{\circ}, ... , X_{j^{i}}(\mid t \mid^{1/l_{i}})) & \mbox{ if } t \geq 0 \\ 
\beta_{i}^{\circ}(X_{i_{1}}(\mid t \mid^{1/l_{i}}) , ... ,  
[X_{i_{k+1}}(\mid t \mid^{1/l_{i}}), X_{i_{k}}(\mid t \mid^{1/l_{i}})]^{\circ}, ... , X_{j^{i}}(\mid t \mid^{1/l_{i}})) & \mbox{ if } t \leq 0
\end{array} \right.$$
After this replacements in the  expression \eqref{tpu1} one obtains a function 
$E_{G}$ which is still invertible in a neighbourhood of $0$. We obtain a function $E_{N}$ with the same algebraic expression as $E_{G}$, but with 
$[\cdot , \cdot ]^{\circ}_{N}$ brackets instead of $[\cdot , \cdot ]^{\circ}_{G}$ ones. Use these functions to (obviously) end the proof of the theorem.  
\quad $\blacksquare$

\begin{thm}
The Gromov-Hausdorff limit of pointed metric spaces $(\mathfrak{g}, 0, \lambda 
d_{G})$ as $\lambda \rightarrow \infty$ exists and equals $(\mathfrak{g}, d_{N})$. 
\end{thm}

\paragraph{Proof.}
We shall use the proposition \ref{pbur}. For this we shall construct  $\varepsilon$ isometries between $Box^{1}_{N}(1)$ and $Box^{1}_{G}(C \varepsilon)$. These are provided by  the function $E_{G} \circ 
E_{N}^{-1} \circ \delta_{\varepsilon}$. 

The trick consists in the definition of the nets. We shall exemplify the construction for the case of a 3 dimensional algebra $\mathfrak{g}$. The basis 
of $\mathfrak{g}$ is $X_{1}, X_{2}, X_{3} \ = \ [X_{1}, X_{2}]_{G}$. Divide 
the interval $[0,X_{1}]$ into $P$ equal parts, same for the interval $[0,X_{2}]$. The interval $[0,X_{3}]$ though will be divided into $P^{2}$ intervals. The net so obtained, seen in $N \ \equiv \ \mathfrak{g}$, turn $E_{G} \circ 
E_{N}^{-1} \circ \delta_{\varepsilon}$ into a $\varepsilon$ isometry between $Box^{1}_{N}(1)$ and $Box^{1}_{G}(C \varepsilon)$. To check this is mostly a matter of smart (but still heavy) notations.  
\quad $\blacksquare$

The proof of this theorem is basically a refinement of Proposition 3.15, Gromov \cite{gromov}, pages 85--86, mentioned in the introduction of these notes. 
Close examination shows that the theorem is a consequence of the following facts: 
\begin{enumerate}
\item[(a)] the identity map $id \ : (\mathfrak{g}, d_{G}) \rightarrow 
(\mathfrak{g}, d_{N})$ has finite dilatation in $0$ (equivalently 
$\exp_{G}: N(D,G) \rightarrow G$ has finite dilatation in $0$), 
\item[(b)] the change from the $G$-invariant distribution induced by $D$ to the 
$N$ invariant distribution induced by the same $D$ is classically smooth. 
\end{enumerate}
It is therefore natural that the result does not feel the non-derivability of 
$id$ map in $x \not = 0$.

\subsection{Noncommutative smoothness in Carnot groups}
\label{nocog}
Let $N$ be a Carnot group and 
\begin{equation}
N \ = \ V_{1} \ + \ ... \ + V_{m}
\label{dirsum}
\end{equation}
the graduation of its algebra. Let $End(N)$ be the group of endomorphisms 
of $N$. 

\begin{defi}
The group of vertical endomorphisms of $N$, noted by $VL(N)$, is the subgroup 
of $End(N)$ of all $N$ endomorphisms $A$ which admit the form $A = (A_{i,j})$
with respect to the direct sum decomposition of $N$ \eqref{dirsum}, such that 
$$A_{ii} \ = \ id_{V_{i}} \ \ ,  \ A_{ij} \ = 0 \ \ \forall i<j \in \left\{ 1, ... , m\right\}$$
\label{defver}
\end{defi}

\begin{prop}
\begin{enumerate}
\item[(a)] The elements of the linear group $HL(N)$ have diagonal form 
with respect to the direct sum decomposition \eqref{dirsum}. 
\item[(b)] Any $A \in End(N)$ admits the unique decomposition 
$$A \ = \ A^{v} A^{h} \ \ , \ A^{v} \in VL(N) \ , \ A^{h} \in HL(N)$$
\item[(c)] $Ad(N)$ is a normal subgroup of $VL(N)$. The quotient group 
is Carnot. 
\end{enumerate}
\end{prop}

\paragraph{Proof.}
(b) Proof by induction. Just write what the algebra morphism condition for 
$A \in End(N)$ means and get the inferior diagonal form of $A$. Such a matrix 
(with elements linear transformations) decomposes uniquely in a product of 
a vertical morphism $A^{v}$ and a diagonal morphism $A^{h}$. Check that the latter commutes with dilatations. 

(a) We know from (b) that any element $A \in HL(N)$ has inferior diagonal form. 
Commutation with dilatations forces $A$ to have diagonal form. 

(c) $Ad(N)$ is normal in $End(N)$ and it is a subgroup of $VL(N)$. It is therefore normal in $VL(N)$. The quotient group is nilpotent, because $VL(N)$ 
is nilpotent. It is Carnot because it has dilatations. Consider the map 
$$\delta^{\varepsilon} \ : End(N) \rightarrow End(N) \ , \ \ \delta^{\varepsilon} A  \ = \ \delta_{\varepsilon} A \delta_{\varepsilon}^{-1}$$
From (b) we get that $\delta^{\varepsilon}$ is an endomorphism of $VL(N)$. 
From the fact that $\delta_{\varepsilon} \in End(N)$ we get that 
$$\delta^{\varepsilon} Ad(N) \ = \ Ad(N)$$
Therefore $\delta^{\varepsilon}$ factorizes to an endomorphism of the 
quotient $VL(N)/Ad(N)$. The rest is trivial. 
\quad $\blacksquare$

We shall call the quotient $$R(N) \ = \ VL(N)/Ad(N)$$
the rest of $N$. 

The definition of noncommutative derivative follows. For any function $f: N \rightarrow N$ and $x \in N$ set 
$$f_{x}: N \rightarrow N \ , \ \ \ f_{x}(y) \ = \ f(x)^{-1}f(xy)$$ 
$f$ is Pansu derivable in $x$ if and only if $f_{x}$ is Pansu derivable 
in $0$; moreover the Pansu derivative of $f$ in $x$ equals the Pansu derivative of $f_{x}$ in $0$.

From the proof of Pansu's Rademacher theorem one can extract the following information: 

\begin{prop}
 if $f$ is Pansu derivable and the convergence 
$$\lim_{\varepsilon \rightarrow 0} \delta_{\varepsilon}^{-1} 
 f_{x}(x \delta_{\varepsilon}y)$$
is uniform with respect to $x$ then $f$ is classically derivable in $x$ along the horizontal directions  and  
we have the connection between Pansu derivative and classical derivative: 
$$P f(x) \ = \ \nabla f(x)\ = \ \left(DL_{f(x)}\right)^{-1}(0) Df(x) DL_{x}(0)$$
Conversely if $f$ is classically derivable along the distribution, such that 
the quantity 
$$\nabla f(x)_{|_{D}} \ = \ \left(DL_{f(x)}\right)^{-1}(0) Df(x) DL_{x}(0)_{|_{D}}$$
can be extended to an algebra morphism of $N$, an such that the convergence 
$$ \lim_{\varepsilon \rightarrow 0} \frac{1}{\varepsilon} \left( 
f(x\delta_{\varepsilon} y) - f(x)\right) $$
(with $y \in D$) is uniform with respect to $x$ then $f$ is Pansu derivable. 
\label{phlp}
\end{prop}

\begin{defi} 
A function $f : N \rightarrow N$ is noncommutative smooth if for any 
$B \in VL(N)$  there is a  continuous function $A_{B}: N \rightarrow VL(N)$ such that:  
\begin{enumerate}
\item[(a)] for any $x \in N$ the map 
$$ y \ \mapsto \ \tilde{f}_{x}^{B}(y) \ = \ A(x)^{-1} f_{x}(B y)$$
is Pansu derivable in $0$, 
\item[(b)] the map 
$$(x,y) \in N^{2} \mapsto (f(x), D\tilde{f}_{x}^{B}(0)y ) \in N^{2}$$
is continuous, 
\item[(c)]  the convergence 
$$Df_{x}(0) y \ = \ \lim_{\varepsilon \rightarrow 0} \delta_{\varepsilon}^{-1} 
 f_{x}(x \delta_{\varepsilon}y)$$
is uniform with respect to $x$. 
\end{enumerate}
The derivative of $f$ in $x$ is by definition 
$$Df(x) \ = \ A_{I}(x) D\tilde{f}_{x}^{I}(x)$$
\label{dnonc}
\end{defi}

Note  that if $f$ is invertible then we have the decomposition:  
$$Df(x)^{v} \ = \ A_{I}(x) \ , \ \ Df(x)^{h} \ = \ D\tilde{f}_{x}^{I}(x)$$
Moreover, if $f$ is invertible and classically derivable then the conclusion 
of proposition \ref{phlp} is still true, namely $D f(x) \ = \ \nabla f(x)$, where the "gradient" $\nabla f(x)$ is given by 
$$\left(DL_{f(x)}\right)^{-1}(0) Df(x) DL_{x}(0)$$

The derivative is called noncommutative because it does not commute with dilatations. 

\begin{prop}
If $f: N \rightarrow N$ is commutative smooth then it is also noncommutative smooth. 

If $f,g: N \rightarrow N$ are noncommutative smooth then so it is $f \circ g$. 
\end{prop}

\paragraph{Proof.}
Because the noncommutative derivative splits the function in a horizontal and a vertical part, it is sufficient to prove that if $f$ is commutative smooth and 
$g$ is noncommutative smooth such that $D^{h}g(x) \ = \ id$ then 
$f \circ g$ is noncommutative smooth. This is contained in the definiton of 
noncommutative smoothness.  
\quad $\blacksquare$

The chain formula is true for noncommutative derivative. But obviously not for the vertical and the horizontal parts.

With this definition of smoothness the right translations in a Carnot group become smooth. It is easy to show that 
$$D \ R_{x}  \ = \ Ad_{x}^{-1} \ , \ \ D^{v} R_{x} \ = \ Ad_{x}^{-1} \ , \ \ 
D^{h} R_{x} \ = \ id$$

This smoothness definition seems to be a good one for Carnot groups. It 
is not good enough though, as the following example shows. 

\subsection{Moment maps and the Heisenberg group}
Let $G \ = \ T^{k} \oplus N$ be the direct sum of the $k$ dimensional torus $T^{k}$ with the Carnot group $N$. The Lie algebra of $G$ is 
$\mathfrak{g} \ = \ R^{k} \oplus N$ (recall that we use the same notation for the Carnot group $N$ and its Lie algebra). 

The algebra $\mathfrak{g}$ is Carnot, with respect to the dilatations
$$\delta_{\varepsilon}(T,X) \ = \ (\varepsilon T , \delta_{\varepsilon} X)$$
The distribution is generated by $R^{k} \oplus V^{1}$ and the nilpotentisation 
of $G$ with respect to this distribution is $N' \ = \ R^{k} \oplus H$.  

The corollary \ref{nonatlas} tells  that $G$ admits a $N'$ manifold structure. 
As the construction of the noncommutative derivative is local in nature, we 
could reproduce it for $N'$ manifolds. In this case that simply means: in reasonings of local nature one can  work on $N'$ instead of $G$.  

Take now $N \ = \ H(n)$ and a toric Hamiltonian action 
$$(s_{1}, ... , s_{k}) \in T^{k} \ \mapsto \phi_{(s_{1},...,s_{k})} \in Sympl(R^{2n})$$
For any $v \in R^{k} \ = \ Lie \ T^{k}$ the infinitesimal generator is 
$$\xi_{v}(x) \ = \ \frac{d}{dt}_{|_{t=0}} \phi_{\exp tv}(x)$$ 
Let $H_{v}$ be the Hamiltonian of $t \mapsto \phi_{\exp tv}$. The moment map associated to the Hamiltonian toric action is then 
$$j_{\phi}: R^{2n} \ \rightarrow \ \left(R^{k}\right)^{*} \ , \ \ \langle j(x), v \rangle \ = \ H_{v}(x)$$

We shall associate to the toric action a function on $T^{K} \times H(n)$. 
\begin{equation}
\Phi: G \ = \ T^{k} \times H(n) \rightarrow G \ , \ \ \Phi(s, x ,\bar{x}) \  = \ (s , \tilde{\phi}_{s}(x , \bar{x}) ) 
\label{amne}
\end{equation} 
where $\tilde{\phi}_{s}$ is the lift of $(\phi_{s}, 0)$ (see theorem \ref{t1}) and it has the expression: 
$$(\phi_{s}(x), 
\bar{x} \ + \ F_{s}(x) ) $$
 
Let us compute the "gradient" 
$$\nabla \Phi (s,x \bar{x}) \ = \ \left( D L_{\Phi(s,x,\bar{x})}(0) \right)^{-1}  D \Phi(s,x,\bar{x}) \ D L_{(s,x,\bar{x})}(0)$$
In matrix form with respect to the decomposition of the algebra
$$R^{k} \oplus H(n) \ = \ R^{k} \oplus R^{2n} \oplus R$$ 
we obtain the following expression: 
$$\nabla \Phi(s,x,\bar{x}) \ =  \ \left( 
\begin{array}{ccc}
1 & 0 & 0 \\ 
\frac{\partial \phi_{s}}{\partial s} & D \phi_{s} & 0 \\ 
a(s,x) & 0 & 1 
\end{array} \right)$$
where $a(s,x)$ is the function: 
$$a(s,x) \ = \ \frac{\partial F}{\partial s} \ - \ \frac{1}{2} \omega (\phi, 
\frac{\partial \phi}{\partial s})$$
We can improve the expression of the function $a$ if we use the fact that the 
toric action is Hamiltonian. Indeed, by direct computation one can show that 
$$\frac{\partial a}{\partial x}(s,x) \ = \ - \ \frac{ \partial }{ \partial x} 
\left( j_{\phi}(\phi_{s}(x)) \right)$$
This is equivalent with 
$$a(s,x) \ = \ \lambda(s) \ - \ j_{\phi}(\phi_{s}(x))$$

Easy computation shows that $End(R^{k} \oplus H(n))$ is made by all matrices 
with the form: 
$$\left( 
\begin{array}{ccc}
A & 0 & 0 \\ 
0 & B & 0 \\ 
a & b & c 
\end{array} \right) \ \ , \ \ A \in GL(R^{k}) \ , \ 
\left( \begin{array}{cc}
B & 0 \\ 
b & c 
\end{array} \right) \ \in \ HL(H(n)) \ , \ a \in \left( R^{k} \right)^{*} \ = \ \left( Lie \ T^{k} \right)^{*} $$
It follows that $\nabla \Phi \not \in End(R^{k} \oplus H(n))$ and we already saw this in a particular form in proposition \ref{pho}. Therefore $\Phi$ is not 
even noncommutative smooth in the sense of definition \ref{dnonc}. 

Suppose that the toric action is by diffeomorphisms with compact support. Then $\lambda(s)$ can be taken equal to $0$ hence  
$a(x,s) \ = \ - j_{\phi}(x)$. 

\begin{prop}
The map $\Phi$ transforms the distribution
$$D^{N}(s,x,\bar{x}) \ = \ \left\{ (S,Y,\bar{Y}) \ \mbox{ : }  - \ \frac{1}{2} \omega(x,Y) \ + \ \bar{Y} \ = \ 0 \right\}$$
in the distribution 
$$D^{\phi}(s,x,\bar{x}) \ = \ \left\{ (S,Y,\bar{Y}) \ \mbox{ : } \langle j_{\phi}(x), S \rangle \ - \ \frac{1}{2} \omega(x,Y) \ + \ \bar{Y} \ = \ 0 \right\}$$
\label{preham}
\end{prop}

\paragraph{Proof.}
By direct computation. 
\quad $\blacksquare$

\subsection{General noncommutative smoothness}
\label{noco}
The introduction of noncommutative derivative makes right translations in Carnot groups "smooth". We shall apply the same strategy for having a manifold structure compatible with the group operation, in the case of a Lie group $G$ with 
a left invariant distribution $D$. 

The transition functions of the chosen atlas have the form  $\left(L^{\mathfrak{g}}_{X}\right)^{-1}$, with $X \in \mathfrak{g}$. These transition functions  forms a group. Note that a natural generalization is the case where the transition functions forms a groupoid. This generalisation is not meaningless since it leads to the introduction of intrinsic SR manifolds, if we identify a manifold with the groupoid of transition functions of an atlas. 

We want this family of transition functions to be "smooth". We shall proceed then as in the motivating example. 

The group $G$ is identified with the group of left translations 
$L^{G}$ acting on the Carnot group $(\mathfrak{g}, \opn)$. These transformation are not smooth, they preserve another distribution than the $N$ left invariant distribution on $G$.

We shall give the following "upgraded" variant of the definition \ref{dnonc},  concerning  the noncommutative derivative.

We shall denote by $N$ the set $\mathfrak{g}$ endowed with the nilpotent multiplication $\opn$. The same notation will be used for the Lie algebra 
$(\mathfrak{g}, [\cdot , \cdot ]_{N})$. We fix on $N$ an  Euclidean norm, given, for example, by declaring orthonormal a basis constructed from multi-brackets from a basis of $D$. The Euclidean norm will be denoted by $\| \cdot \|$.

We introduce the group of linear transformations $\Sigma(G,D)$ generated by 
all transformations $DL^{g}_{X}(0)$, $DR^{g}_{X}(0)$, $HL(N)$. In particular 
$\delta_{\varepsilon} \in \Sigma(G,D)$ for any $\varepsilon > 0$. From the 
formula: 
$$DL^{N}_{X}(0) \ = \ \lim_{\varepsilon \rightarrow 0} \delta_{\varepsilon}^{-1} \circ DL_{\delta_{\varepsilon}X}^{g}(0) \circ \delta_{\varepsilon}$$
we get that $DL^{N}_{X}(0) \in \Sigma(G,D)$. In the same way we obtain 
$DR^{N}_{X}(0) \in \Sigma(G,D)$. 

For any $F \in \Sigma(G,D)$ we define the group $N(F)$ to be the Carnot group 
$N$ with the bracket: 
$$[X,Y]^{F}_{N} \ = \ F [F^{-1} X , F^{-1}Y]_{N}$$
and with the dilatations 
$$\delta_{\varepsilon}^{F} \ = \ F \delta_{\varepsilon} F^{-1}$$
This group is endowed with the Carnot-Carath\'eodory distance generated by the Euclidean norm and the distribution $D^{N}(F) \ = \ F D$. 

We can see the class $\left\{ N(F) \mbox{ : } F \in \Sigma(G,D) \right\}$ in two ways, as an orbit. First, it is the orbit 
$$\bar{\Sigma}(G,D) \ = \ \left\{ (F \delta_{\varepsilon} F^{-1} , F \ ad^{N}_{F^{-1}X} F^{-1}) \ \mbox{ : } F \in \Sigma(G,D) \right\}$$ 
in the space $GL(N) \times gl(N,gl(N))$, under the action of $\Sigma(G,D)$ 
$$ F. (A,B) \ = \ (FAF^{-1}, F(B)) \ , \ \ F(B)(Z) \ = \ F B(F^{-1}Z) F^{-1}$$
Second, it is the right coset $\Sigma(G,D)/HL(N)$. Indeed, if $N(F) \ = \ N(P)$ then $F^{-1}P \in HL(N)$.

\begin{defi}
 A function $f : N \rightarrow N$ is $(G,D)$ noncommutative derivable in 
$x \in N$ if for any $F \in \Sigma(G,D)$ there exists $P \in \Sigma(G,D)$ 
such that $f: N(F) \rightarrow N(P)$ is Pansu derivable in $x$. The derivative 
of $f$ in $x$ is by definition: 
$$Df(x)(F,P) \ = \ P^{-1} \circ D_{Pansu} f(x) \circ F$$

The function $f$ is noncommutative smooth if 
\begin{enumerate}
\item[(a)] the map 
$$(x,y) \in N^{2} \mapsto (f(x), Df(x)y ) \in N^{2}$$
is continuous,
\item[(c)]  the convergence
$$ D_{Pansu} f(x)( y \ = \ \lim_{\varepsilon \rightarrow 0} \delta_{\varepsilon}^{-1,P}  \left( P(P^{-1}(f(x))^{-1} \ P^{-1}f(x \delta_{\varepsilon}^{F}y))\right)$$
is uniform with respect to $x$. 
\end{enumerate}
\label{dnong}
\end{defi}

If $f$ is invertible then following map is well defined: 
$$\Lambda(f,x): \bar{\Sigma}(G,D) \rightarrow \bar{\Sigma}(G,D)$$ 
$\Lambda(f,x)(F \ HL(N)) \ = \ P \ HL(N)$, where $f: N(F) \rightarrow N(P)$ is 
Pansu derivable in $x$. In this case the noncommutative derivative can be written as $Df(x,N(F))$. 

The definition can be adapted for functions $f: N \rightarrow R$, or 
$f: R \rightarrow N$. 

\begin{defi}
A function $f: N \rightarrow R$ is noncommutative derivable in $x$ if 
for any $F \in \Sigma(G,D)$ the function $f: N(F) \rightarrow R$ is Pansu derivable in $x$. In this case 
$$Df(x,N(F)) \ = \ D_{Pansu} f(x) \circ F$$
A function $f: R \rightarrow N$ is noncommutative derivable in $x \in R$ if there exists $P \in \Sigma(G,D)$ such that $f: R \rightarrow N(P)$ is Pansu derivable in $x \in R$. In this case 
$$Df(x) \ = \ P^{-1} \circ D_{Pansu} f(x)$$
The definition of noncommutative smoothness adapts obviously. 
\end{defi}

We might risk to have no noncommutative derivable functions from $N$ to $R$. On the contrary we shall have a lot of noncommutative derivable functions from 
$R$ to $N$. 

What is the connection with the previous notion of noncommutative smoothness? 
The noncommutative smoothness in the sense of definition \ref{dnonc} corresponds to the  definition \ref{dnong} if we replace $\Sigma(G,D)$ with 
$End(N)$.

The following proposition has now a straightforward proof. 

\begin{prop}
The transition functions of the  $N(G,D)$ atlas $\mathcal{A}$ are 
$\Sigma(G,D)$ smooth. 
\end{prop}

\paragraph{Proof.}
Indeed, it is sufficient to remark that $\mathcal{J}(G,D) \subset \Sigma(G,D)$. 
\quad $\blacksquare$

\subsection{Margulis \& Mostow tangent bundle}
In this section we shall apply Margulis \& Mostow \cite{marmos2} construction of the tangent bundle to a SR manifold  for the case of a group with left invariant distribution. It will turn that the tangent bundle does not have a group structure, due to the fact that, as previously, the non-smoothness of the right translations is not studied. 

The main point in the construction of a tangent bundle is to have a functorial definition of the tangent space. This is achieved by Margulis \& Mostow \cite{marmos2} in a very natural way. One of the geometrical definitions of a tangent 
vector $v$ at a point $x$, to a manifold $M$, is the following one: identify $v$ with the class of smooth curves which pass through $x$ and have tangent $v$. 
If the manifold $M$ is endowed with a distance then one can  define the equivalence relation based in $x$ by:  $c_{1} \equiv_{x} 
c_{2}$ if $c_{1}(0)  = c_{2}(0) = x$ and the distance between $c_{1}(t)$ and 
$c_{2}(t)$ is of order $t^{2}$ for small $t$.  The set of equivalence classes 
is the tangent space at $x$. One has to put then some structure on the tangent 
space (as, for example, the nilpotent multiplication). 

To put is practice this idea is not so easy though. This is achieved by the following sequence of definitions and theorems. For commodity we shall explain this construction in the case $M=G$ connected Lie group, endowed with a left 
invariant distribution $D$. The general case is the one of a regular sub-Riemannian manifold. We shall denote by $d_{G}$ the CC distance on $G$ and 
we identify $G$ with $\mathfrak{g}$, as previously. The CC distance induced by 
the distribution $D^{N}$, generated by left translations of $G$ using nilpotent 
multiplication $\opn$, will be denoted by $d_{N}$. 

\begin{defi}
A $C^{\infty}$ curve in $G$ with  $x \ = \ c(0)$ is called rectifiable at $t=0$ if 
$d_{G}(x, c(t)) \leq Ct$ as $t \rightarrow 0$. 

Two $C^{\infty}$ curves $c', c"$ with $c'(0) = x = c"(0)$ are called equivalent at $x$ if $t^{-1} d_{G}(c'(t),c"(t)) \rightarrow 0$ as $t \rightarrow 0$. 

The tangent cone to $G$ as $x$, denoted by $C_{x} G$ is the set of equivalence classes of all $C^{\infty}$ paths $c$ with $c(0) = x$, rectifiable at $t=0$. 
\end{defi}

Let $c: [-1,1] \rightarrow G$ be a $C^{\infty}$ rectifiable curve, $x = c(0)$ and 
\begin{equation}
v \ = \ \lim_{t \rightarrow 0} \delta_{t}^{-1} \left( c(0)^{-1} \opg c(t) \right) 
\label{ax1}
\end{equation}
The limit $v$ exists because the curve is rectifiable. 

Introduce the curve $c_{0}(t) \ = \ x \exp_{G}(\delta_{t}v)$. Then 
$$d(x, c_{0}(t)) \ = \ d(e, x^{-1}c(t)) < \mid v \mid t$$
 as $t \rightarrow 0$  (by the Ball-Box theorem) 
The curve  $c$  is equivalent with $c_{0}$. Indeed, we have (for $t>0$): 
$$\frac{1}{t} d_{G}(c(t), c_{0}(t)) \ = \ \frac{1}{t} d_{G}(c(t), x \opg \delta_{t} v)  \ = \ \frac{1}{t} d_{G}(\delta_{t}(v^{-1}) \opg x^{-1} \opg c(t) , 0)$$
The latter expression is equivalent (by the Ball-Box Theorem) with 
$$\frac{1}{t} d_{N} (\delta_{t}(v^{-1}) \opg x^{-1} \opg c(t) , 0) \ = \ 
d_{N} (\delta_{t}^{-1} \left(  \delta_{t}(v^{-1}) \opg \delta_{t}\left( 
\delta_{t}^{-1} \left( x^{-1} \opg c(t)\right) \right) \right) ) $$
The right hand side (RHS) converges to 
$d_{N}(v^{-1} \opn v , 0)$, as $t \rightarrow 0$, as a consequence of 
the definition of $v$ and theorem \ref{teore}.

Therefore we can identify $C_{x} G$ with the set of curves 
$t \mapsto x \exp_{G}(\delta_{t} v)$, for all $v \in \mathfrak{g}$. Remark that the equivalence relation between curves $c_{1}$, $c_{2}$, such that 
$c_{1}(0) \ = \ c_{2}(0) \ = \ x$ can be redefined as: 
\begin{equation}
\lim_{t \rightarrow 0} \delta_{t}^{-1} \left( c_{2}(t)^{-1} \opg  c_{1}(t) \right) \ = \ 0 
\label{ax2}
\end{equation}
 
In order to define the multiplication Margulis \& Mostow introduce  the families of segments rectifiable at  $t$. 

\begin{defi}
A family of segments rectifiable at $t=0$ is a $C^{\infty}$ map 
$$\mathcal{F} : U  \rightarrow G$$
where $U$ is an open neighbourhood of $G \times 0$ in $G \times R$ satisfying 
\begin{enumerate}
\item[(a)] $\mathcal{F}(\cdot, 0) \ = \ id $
\item[(b)] the curve $t \mapsto \mathcal{F}(x,t)$ is rectifiable at $t=0$ uniformly for all $x \in G$, that is for every compact  $K$  in $G$ there is a  constant $C_{K}$ and a compact neighbourhood $I$ of $0$ such that $d_{G}(y,\mathcal{F}(y,t)) < C_{K} t $ for all $(y,t) \in K \times I$. 
\end{enumerate}

Two families of segments rectifiable at $t=0$ are called equivalent if 
$t^{-1} d_{G}(\mathcal{F}_{1}(x,t), \mathcal{F}_{2}(x,t)) \rightarrow 0$ as $t\rightarrow 0$, uniformly on compact sets in the domain of definition. 
\end{defi}

Part (b) from the  definition of a family of segments rectifiable can be restated as: 
there exists the limit 
\begin{equation}
v(x) \ = \ \lim_{t \rightarrow 0} \delta_{t}^{-1} \left( x^{-1} \opg \mathcal{F}(x,t) \right)
\label{ax3}
\end{equation}
and the limit is uniform with respect to $x \in K$, $K$ arbitrary compact set. 

It follows then, as previously, that $\mathcal{F}$ is equivalent to 
$\mathcal{F}_{0}$, defined by: 
$$\mathcal{F}_{0}(x,t) \ = \ x \opg \delta_{t} v(x)$$
Also, the equivalence between families of segments rectifiable can be redefined 
as: 
\begin{equation}
\lim_{t \rightarrow 0} \delta_{t}^{-1} \left( \mathcal{F}_{2}(x,t)^{-1} \opg  \mathcal{F}_{1}(x,t) \right) \ = \ 0 
\label{ax2}
\end{equation}
uniformly with respect to $x \in K$, $K$ arbitrary compact set.

\begin{defi}
The product of two families $\mathcal{F}_{1}$, $\mathcal{F}_{2}$ of segments rectifiable at $t=0$ is defined by 
$$ \left( \mathcal{F}_{1} \circ \mathcal{F}_{2}\right) (x,t) \ = \ \mathcal{F}_{1}(\mathcal{F}_{2}(x,t),t)$$
\end{defi}

The product is well defined by Lemma 1.2 {\it op. cit.}. One of the main results is then the following theorem (5.5). 

\begin{thm}
Let $c_{1}$, $c_{2}$ be $C^{\infty}$ paths rectifiable at $t=0$, such that 
$c_{1}(0) = x_{0} = c_{2}(0)$. Let 
$\mathcal{F}_{1}$, $\mathcal{F}_{2}$ be two families of segments rectifiable 
at $t=0$ with: 
$$\mathcal{F}_{1}(x_{0}, t) \ = \ c_{1}(t) \ \ , \ \ \mathcal{F}_{2}(x_{0}, t) \ = \ c_{2}(t)$$
Then the equivalence class of 
$$t \mapsto \mathcal{F}_{1} \circ \mathcal{F}_{2}(x_{0}, t)$$ depends 
only on the equivalence classes of $c_{1}$ and $c_{2}$. This defines the product of the elements of the tangent cone $C_{x_{0}} G$. 
\end{thm}

This theorem is the straightforward consequence of the following facts (5.1(5) and 5.2 in Margulis \& Mostow \cite{marmos}).

\begin{rk}
 Maybe I misunderstood the notations, but it 
seems to me that several times the authors claim that the exponential map which they construct is bi-Lipschitz (as in 5.1(4) and Corollary 4.5). This is false, as explained before. In Bella\"{\i}che \cite{bell}, Theorem 7.32 and 
also at the beginning of section 7.6 we find that the exponential map is only 
$1/m$ H\"{o}lder continuous (where $m$ is the step of the nilpotentization). 
However, (most of) the results of Margulis \& Mostow hold true. It would be interesting to have a better written paper on the subject, even if the subject might seem trivial (which is not). 
\end{rk}

 We shall denote by $\mathcal{F} \approx \mathcal{F}'$ the equivalence relation of families of segments rectifiable; the equivalence relation of rectifiable curves based at $x$ will be denoted by $c \stackrel{x}{\approx} c'$. 

\begin{lema}
\begin{enumerate}
\item[(a)] Let $\mathcal{F}_{1} \approx \mathcal{F}_{2}$ and 
$\mathcal{G}_{1} \approx \mathcal{G}_{2}$. Then 
$\mathcal{F}_{1} \circ \mathcal{G}_{1} \approx \mathcal{F}_{2} \circ \mathcal{G}_{2}$. 
\item[(b)] The map $\mathcal{F} \ mapsto \mathcal{F}_{0}$ is constant on equivalence classes of families of segments rectifiable. 
\end{enumerate}
\label{lhelps}
\end{lema}

\paragraph{Proof.}
Let $$\mathcal{F}_{0} (x,t) \ = \ x \opg \delta_{t} w_{1}(x) \ \ , \ \ 
\mathcal{G}_{0} (x,t) \ = \ x \opg \delta_{t} w_{2}(x)$$
For the point (a) it is sufficient to prove that 
$$\mathcal{F} \circ \mathcal{G}  \approx  \mathcal{F}_{0} \circ \mathcal{G}_{0}$$
This is true by the following chain of estimates. 
$$\frac{1}{t} d_{G}(\mathcal{F} \circ \mathcal{G}(x,t) , \mathcal{F}_{0} \circ \mathcal{G}_{0}(x,t)) \ = \ \frac{1}{t} d_{G}(\delta_{t}w_{1}(\mathcal{G}_{0}(x,t))^{-1} \opg \delta_{t}w_{2}(x)^{-1} \opg x^{-1} \opg \mathcal{F}(\mathcal{G}(x,t),t) ,  0)$$
The RHS of this equality behaves like 
$$d_{N}(\delta_{t}^{-1} \left( \delta_{t}w_{1}(\mathcal{G}_{0}(x,t))^{-1} 
\opg \delta_{t}w_{2}(x)^{-1} \opg \delta_{t} \left( \delta_{t}^{-1} \left( 
x^{-1} \opg \mathcal{G}(x,t)\right) \right) \opg \delta_{t} \left( 
\delta_{t}^{-1} \left( \mathcal{G}(x,t)^{-1} \opg \mathcal{F}(\mathcal{G}(x,t),t) \right) \right) \right) , 0)$$
This quantity converges (uniformly with respect to $x \in K$, $K$ an arbitrary compact) to 
$$d_{N}(w_{1}(x)^{-1} \opn w_{2}(x)^{-1} \opn w_{2}(x) \opn w_{1}(x), 0) \ = \ 0$$

The point (b) is easier: let $\mathcal{F} \approx \mathcal{G}$ and consider 
$\mathcal{F}_{0}$, $\mathcal{G}_{0}$, as above. We want to prove that 
$\mathcal{F}_{0} \ = \ \mathcal{G}_{0}$, which is equivalent to $w_{1} \ = \ 
w_{2}$. 

Because $\approx $ is an equivalence relation all we have to prove is that 
if $\mathcal{F}_{0} \approx \mathcal{G}_{0}$ then $w_{1} \ = \ w_{2}$.  We have: 
$$\frac{1}{t} d_{G}(\mathcal{F}_{0}(x,t), \mathcal{G}_{0}(x,t)) \ = \ 
\frac{1}{t} d_{G}(x \opg \delta_{t} w_{1}(x) , x \opg \delta_{t} w_{2}(x))$$
We use the $\opg$ left invariance of $d_{G}$ and the Ball-Box theorem to 
deduce that the RHS behaves like 
$$d_{N}(\delta_{t}^{-1} \left( \delta_{t} w_{2}(x)^{-1} \opg \delta_{t} w_{1}(x)^{-1}\right) , 0)$$
which converges to $d_{N}(w_{1}(x), w_{2}(x))$ as $t$ goes to $0$. The two families are equivalent, therefore the limit equals $0$, which implies that 
$w_{1}(x) = w_{2}(x)$ for all $x$. 
\quad $\blacksquare$

We shall apply this theorem. Let $c_{i}(t) \ = \ x_{0} \exp_{G} \delta_{t} v_{i}$, for $i = 1,2$. It is easy to check that 
$\mathcal{F}_{i}(x,t) \ = \ x \exp_{G}(\delta_{t} v_{i})$ are families of segments rectifiable at $t=0$ which satisfy the hypothesis of the theorem. 
But then 
$$\left( \mathcal{F}_{1} \circ \mathcal{F}_{2}\right) (x,t) \ = \ 
x_{0} \exp_{G} (\delta_{t} v_{1}) \exp_{G}(\delta_{t} v_{2})$$
which is equivalent with 
$$\mathcal{F} \exp_{G}\left( \delta_{t}( v_{1} \opn v_{2} ) \right)$$
Therefore the tangent bundle defined by this procedure is the same as the virtual tangent bundle defined previously. It's definition is possible because 
the group operation has horizontal derivative in $(e,e)$. 

In terms of commutative derivative, theorem 10.5 Margulis \& Mostow 
\cite{marmos1} (and restricting to bi-Lipschitz maps) becomes the Rademacher theorem.   
Indeed, in the case of a Lie group $G$ endowed with a left invariant distribution $D$, with an associated Carnot-Carath\'eodory distance, the definition \ref{fdcd} of commutative derivative can be adapted in the following way:

\begin{defi}
Let $G_{1}$, $G_{2}$ two groups endowed with left invariant distributions and 
$d_{1}$, $d_{2}$ tow associated Carnot-Carath\'eodory distances. 

A function $f: G_{1} \rightarrow G_{2}$ is metrically commutative derivable in 
$x \in G_{1}$ if there is a $\varepsilon > 0$ such that the sequence 
$$\Delta_{\lambda}f(x)u \ = \ \delta_{\lambda}^{-1} 
\left(f(x)^{-1}f(x \delta_{\lambda}u)\right)$$
converges uniformly with respect to $u \in B(x,\varepsilon)$. The derivative 
is
$$Df(x)u \ = \ \lim_{\varepsilon \rightarrow 0} \delta_{\varepsilon}^{-1} 
\left(f(x)^{-1}f(x \delta_{\varepsilon}u)\right)$$
and the virtual tangent is 
$$VT f(x): VT_{x}G_{1} \rightarrow VT_{f(x)} G_{2} \ , \ \ VT f(x) L^{N_{1},x}_{y} \ = \ 
L^{G_{2}}_{f(x)} L^{N_{2}}_{Df(x)y} L^{-1, G_{2}}_{f(x)}$$
\end{defi}

\begin{thm}
Let $\phi: G \rightarrow G$ be a bi-Lipschitz map with respect to the Carnot-Carath\'eodory distance on $G$ induced by the left invariant distribution $D$. 
Then for almost all $x \in G$  the map $\phi$ is commutative derivable at $x$ and  the virtual tangent to $\phi$ is an isomorphism of conical groups. 
\end{thm}

We don't give the proof of this theorem. The reader might find interesting to 
adapt the proof of the mentioned Margulis \& Mostow theorem in our particular case.


\begin{thebibliography}{99}

\bibitem{allcock} D. Allcock, An isoperimetric inequality for the
Heisenberg groups, {\it Geom. Funct. Anal.}, {\bf 8}, (1998), 219--233



\bibitem{ambrosio} L. Ambrosio, Geometric measure theory and applications to the
calculus of variations, (1999), chapters 1--4


\bibitem{bell} A. Bella\"{\i}che, The tangent space in sub-Riemannian 
geometry, in: {\it Sub-Riemannian Geometry}, A. Bella\"{\i}che, J.-J. Risler 
eds., Progress in Mathematics, {\bf 144}, Birchäuser, (1996), 4 -- 78

\bibitem{buli} M. Buliga, The topological substratum of the derivative, 
{\it Stud. Cerc. Mat. (Mathematical Reports)}, {\bf 45}, no. 6, (1993), 
453--465



\bibitem{burago} D. Burago, Y. Burago, S. Ivanov, A Course in Metric Geometry, {\it 
Graduate Studies in Mathematics}, {\bf 33}, AMS Providence, Rhode Island, 
(2000)

\bibitem{cdkr} M. Cowling, A Dooley, A. Koranyi, F. Ricci, 
H-type groups and Iwasawa decompositions, {\it Adv. in Math.}, 
{\bf 87} (1991), 1--41

\bibitem{elia} Y. Eliashberg, L. Polterovich, Biinvariant metrics on the group of Hamiltonian diffeomorphisms, {\it International J. of Math.}, {\bf 
4}, no. 5, (1993), 727--738 

\bibitem{fostein}  G.B. Folland, E.M. Stein, Hardy spaces on homogeneous groups, 
{\it Mathematical notes}, {\bf 28}, Princeton - N.J.: Princeton University Press; Tokyo: University of Tokyo Press, (1982)

\bibitem{fraseca}  B. Franchi, R. Serapioni, F. Serra Cassano, Rectifiability and Perimeter in the  Heisenberg Group, {\it Math. Ann.}, {\bf 321}, (2001), 479--531

\bibitem{ganh} N. Garofalo, D.M. Nhieu, Isoperimetric and Sobolev inequalitites for Carnot-Carath\'eodory spaces and the existence of minimal surfaces, 
{\it Comm. Pure Appl. Math.}, {\bf 49} (1996), no. 10, 1081--1144

\bibitem{goo} R. Goodman, Nilpotent Lie groups: structure and applications to
analysis, {\it Lecture Notes in Math.} , {\bf 53} (1981) 


\bibitem{gromo} M. Gromov, Carnot-Caratheodory spaces seen from within, 
in: {\it Sub-Riemannian Geometry}, A. Bella\"{\i}che, J.-J. Risler 
eds., Progress in Mathematics, {\bf 144}, Birchäuser, (1996), 79 -- 323

\bibitem{gromov} M. Gromov, Metric structures for Riemannian and non-Riemannian
spaces, {\it Progress in Math.}, {\bf 152}, Birch\"auser (1999), 
chapter 1. 

\bibitem{hako} P. Hajlasz, P. Koskela, Sobolev met Poincar\'e, {\it Memoirs of the America Mathematical Society} (2000), {\bf 688}

\bibitem{hnp} J. Hilgert, K.-H. Neeb, W. Plank, Symplectic convexity theorems 
and coadjoint orbits, {\it Compositio Mathematica}, {\bf 94}, 129--180, (1994)

\bibitem{hozen} H. Hofer, E. Zehnder, Symplectic invariants and Hamiltonian dynamics, {\it 
 Birkhäuser Advanced Texts: Basler Lehrbücher},  Birkhäuser Verlag, Basel, (1994)

\bibitem{kore} A. Koranyi, M. Reimann, Foundations for the theory of quasi-conformal mappings of the Heisenberg group, {\it Adv. in Math.}, {\bf III:1} 
(1995), 1--87

\bibitem{kos} B. Kostant, On convexity, the Weyl group and the Iwasawa decomposition, 
{\it Ann. Scient. Ec. Norm. Sup.},  $4^{e}$ s\'erie, {\bf t. 6} (1973), 
413--455



\bibitem{magnani} V. Magnani, Differentiability and area formula on stratified Lie groups, preprint http:\\www.cvgmt.it, (2000)

\bibitem{marmos1} G.A. Margulis, G.D. Mostow, The differential of a quasi-conformal mapping of a Carnot-Carath\'eodory space, {\it Geom. Funct. Analysis}, {\bf 8} (1995), 2, 402--433

\bibitem{marmos2} G.A: Margulis, G.D. Mostow, Some remarks on the definition of tangent cones in a Carnot-Carath\'eodory space, {\it J. D'Analyse Math.}, {\bf 80} (2000), 299--317

\bibitem{mosto} G.D. Mostow, Strong rigidity of locally symmetric spaces, {\it Ann. of Math. Studies}, {\bf 78} (1973), Princeton Univ. Press, Princeton

\bibitem{mit} J. Mitchell, On Carnot-Carath\'eodory metrics, 
{\it J. Diff. Geometry}, {\bf 21} (1985), 35--45

\bibitem{montgo}R. Montgomery, Survey of singular geodesics, in: {\it Sub-Riemannian Geometry}, A. Bella\"{\i}che, J.-J. Risler 
eds., Progress in Mathematics, {\bf 144}, Birchäuser, (1996), 325--339


\bibitem{pansu} P. Pansu, M\'etriques de Carnot-Carath\'eodory at
quasiisometries des espaces symetriques de rang un, {\it Ann. of
Math.}, {\bf 129} (1989), 1--60

\bibitem{pansu1} P. Pansu, Croissance des boules et des g\'eod\'esiques ferm\'ees dans les nilvari\'et\'es, {\it Ergod. Th. \& Dynam. Sys.}, {\bf 3} (1983), 415--445

\bibitem{paninis} P. Pansu, Une in\'egalit\'e isop\'erim\'etrique sur le groupe 
d'Heisenberg, {\it C. R. Acad. Sci. Paris}, {\bf 295} (1982), 127--131


\bibitem{polte}  L. Polterovich, The geometry of the group of symplectic diffeomorphisms, 
{\it  Lectures in Mathematics ETH Zürich}, 
Birkhäuser Verlag, Basel, (2001)

\bibitem{reiric} H.M. Reimann, F. Ricci, The complexified Heisenberg group, 
{\it PRoc. on Analysis and Geometry}, Novosibirsk: Sobolev Institute Press, (2000), 
465--480


\bibitem{sleewa} P. Sleewaegen, Application moment et theoreme de convexite de Kostant, 
PhD. thesis, Univ. Libre de Bruxelles, 1999 

\bibitem{sus} H.J. Sussmann, A cornucopia of four-dimensional abnormal sub-Riemannian minimizers, in: {\it Sub-Riemannian Geometry}, A. Bella\"{\i}che, J.-J. Risler 
eds., Progress in Mathematics, {\bf 144}, Birchäuser, (1996), 342--364

\bibitem{varsaco} N. Varopoulos, L. Saloff-Coste, Th. Coulhon, Analysis and Geometry on Groups, {\it Cambridge Univ. Press}, (1993)

\bibitem{vodopis} S.K. Vodop'yanov, $\mathcal{P}$-Differentiability on Carnot Groups 
in Different Topologies and Related Topics, {\it Proc. on Analysis ans Geometry}, 
Novosibirsk: Sobolev Institute Press, (2000), 603 -- 670

\bibitem{voduk} S.K. Vodop'yanov, A.D. Ukhlov, Approximately differentiable 
transformations and change of variables on nilpotent groups, {\it Sib. Math. J.}, {\bf 37} (1996), 1, 62--78











\end{thebibliography}
\end{document}